\documentclass[11pt]{article}

\usepackage{float} 
\usepackage{amsmath}
\usepackage{amssymb}
\usepackage{lmodern}
\usepackage{array}
\usepackage{bbm}
\usepackage{makeidx}
\usepackage{hyperref}
\usepackage{cleveref}
\usepackage{caption}
\usepackage{subcaption}
\usepackage{amsthm}
\usepackage{color}
\usepackage{mathtools}
\usepackage{graphicx} 
\usepackage[LGR,T1]{fontenc}
\usepackage{mathrsfs}
\usepackage{lipsum}
\usepackage{booktabs}
\usepackage{xcolor}
\usepackage[title]{appendix}
\hypersetup{
    colorlinks,
    linkcolor={blue!80!black},
    citecolor={green!50!black},
}
\usepackage{apptools}
\usepackage{tikz}
\usetikzlibrary{arrows, positioning, automata}
\usepackage{comment}
\usepackage[shortlabels]{enumitem}

\usepackage[noend]{algorithmic}
\usepackage[ruled,vlined]{algorithm2e}

\usepackage[mathscr]{euscript}
\DeclareSymbolFont{rsfs}{U}{rsfs}{m}{n}
\DeclareSymbolFontAlphabet{\mathscrsfs}{rsfs}

\def\R{{\mathbb R}}

\def\cR{{\mathcal R}}

\def\cuG{\mathscrsfs{G}}
\def\cuD{\mathscrsfs{D}}
\def\cuC{\mathscrsfs{C}}

\def\Cccv{\cuC_{\mbox{\tiny\rm ccv}}}
\def\Ccont{\cuC_{\mbox{\tiny\rm cont}}}

\def\de{{\rm d}}

\def\sRS{\mbox{\tiny \rm RS}}

\def\uq{\underline{q}}

\def\<{\langle}
\def\>{\rangle}

\def\salg{\mbox{\tiny \rm alg}}

\def\prob{{\mathbb P}}

\def\ovp{\overline{p}}
\def\pd{\dot{p}}
\def\AC{{\rm AC}}

\usepackage{mathtools}

\newtheoremstyle{myremark} 
    {\topsep}                    
    {\topsep}                    
    {\rm}                        
    {}                           
    {\bf}                        
    {.}                          
    {.5em}                       
    {}  

\newtheorem{claim}{Claim}[section]
\newtheorem{lemma}[claim]{Lemma}

\newtheorem{theorem}{Theorem}
\newtheorem{proposition}[claim]{Proposition}
\newtheorem{corollary}[claim]{Corollary}
\newtheorem{definition}[claim]{Definition}

\theoremstyle{myremark}

\global\long\def\uu{\phi}%
\global\long\def\Qm{\calQ}%
\global\long\def\Qe{\calQ_{\rm err}}%
\global\long\def\lb{\text{\ensuremath{{\scriptscriptstyle {\rm lb}}}}}%
\global\long\def\LL{\mathbb{L}}%
\global\long\def\cL{\mathcal{L}}%
\global\long\def\kmax{k_{\max}}%
\global\long\def\alg{{\sf alg}}%
\global\long\def\r{\varnothing}%
\global\long\def\Ia{{\rm I}}%
\global\long\def\Ib{{\rm II}}%
\global\long\def\Ic{{\rm III}}%
\global\long\def\Iabar{\overline{{\rm I}}}%
\global\long\def\Ibbar{\overline{{\rm II}}}%
\global\long\def\Icbar{\overline{{\rm III}}}%
\global\long\def\cEr{\calE_{\delta,\r}}%
\global\long\def\cEm{\calE_{\delta,{\rm main}}}%
\global\long\def\cEbar{\calE_{\delta,{\rm apx}}}%

\allowdisplaybreaks
\usepackage[top=1in,bottom=1in,left=1in,right=1in]{geometry}
\usepackage[utf8]{inputenc}

\def\dd{\mathrm{d}}
\def\ed{\stackrel{{\mathrm d}}{=}}

\def\top{\intercal}

\def\bA{{\boldsymbol A}}
\def\bB{{\boldsymbol B}}

\def\hcL{\hat{\mathcal L}}
\def\xid{\xi'}
\def\xidd{\xi''}
\def\xidot{\dot{\xi}}
\def\xiddot{\ddot{\xi}}

\def\cA{{\mathcal A}}

\def\cS{{\mathcal S}}

\def\cF{\mathcal{F}}

\def\cS{\mathcal{S}}

\def\GOE{\mbox{GOE}}

\def\S{{\mathbb S}}

\def\step{{\sigma}}

\def\cF{\mathcal{F}}

\def\LL{\mathbb{L}}

\def\TT{\mathbb{T}}

\def\ZZ{\mathbb{Z}}

\def\calA{\mathcal{A}}

\def\calD{\mathcal{D}}
\def\calE{\mathcal{E}}

\def\calN{\mathcal{N}}

\def\calQ{\mathcal{Q}}

\def\bA{\mathbf{A}}
\def\bB{\mathbf{B}}

\def\bD{\mathbf{D}}

\def\bF{\mathbf{F}}
\def\bG{\mathbf{G}}

\def\bI{\mathbf{I}}

\def\bM{\mathbf{M}}

\def\bP{\mathbf{P}}

\def\bU{\mathbf{U}}

\def\bW{\mathbf{W}}

\def\bZ{\mathbf{Z}}

\def\bb{\boldsymbol{b}}

\def\be{\boldsymbol{e}}

\def\bg{\boldsymbol{g}}
\def\bh{\boldsymbol{h}}

\def\bm{\boldsymbol{m}}

\def\bs{\boldsymbol{s}}
\def\bt{\boldsymbol{t}}
\def\bu{\boldsymbol{u}}
\def\bv{\boldsymbol{v}}
\def\bw{\boldsymbol{w}}
\def\bx{\boldsymbol{x}}
\def\by{\boldsymbol{y}}
\def\bz{\boldsymbol{z}}

\def\bA{\boldsymbol{A}}
\def\bB{\boldsymbol{B}}

\def\bD{\boldsymbol{D}}

\def\bF{\boldsymbol{F}}
\def\bG{\boldsymbol{G}}

\def\bI{\boldsymbol{I}}

\def\bM{\boldsymbol{M}}

\def\bP{\boldsymbol{P}}

\def\bU{\boldsymbol{U}}

\def\bW{\boldsymbol{W}}

\def\bZ{\boldsymbol{Z}}

\def\tbg{\tilde{\boldsymbol g}}

\def\normal{{\mathsf{N}}}

\def\bzero{\boldsymbol{0}}

\def\bm{\boldsymbol{m}}

\def\GOE{{\sf GOE}}

\def\ALG{{\sf ALG}}

\def\reals{{\mathbb R}}

\def\sT{{\sf T}}

\def\Ball{{\sf B}}

\def\id{{\boldsymbol I}}

\def\ed{\stackrel{{\rm d}}{=}}
\def\op{{\mbox{\rm\tiny op}}}

\def\drift{{\mathcal F}}

\renewcommand{\P}{\mathbb{P}}
\newcommand{\E}{\mathbb{E}}

\newcommand{\N}{\mathcal{N}}

\newcommand{\eps}{\varepsilon}
\newcommand{\Var}{\operatorname{Var}}
\newcommand{\argmax}{\operatorname{argmax}}
\newcommand{\argmin}{\operatorname{argmin}}

\newcommand{\Tr}{\operatorname{Tr}}

\newcommand{\RN}[1]{%
  \textup{\uppercase\expandafter{\romannumeral#1}}%
}

\newcommand{\RNum}[1]{\uppercase\expandafter{\romannumeral #1\relax}}



\makeatletter
\newcommand*{\rom}[1]{\expandafter\@slowromancap\romannumeral #1@}
\makeatother

\begin{document}

\title{Tight Lipschitz Hardness for\\
 Solving Random Polynomial Equations}

\author{Andrea Montanari\thanks{Department of Statistics and Department of Mathematics, Stanford University} 
	\and 
	Eliran Subag\thanks{Department of Mathematics, Weizmann Institute of Science}
}
\pagenumbering{arabic}

\date{\today}

\maketitle

\begin{abstract}
We consider the problem of solving $n$ polynomial equations in $d$ real variables:
$F_1(\bx)=0, \ldots, F_n(\bx)=0$, with unknowns $\bx=(x_1, \ldots, x_d)$ satisfying
the normalization constraint $\|\bx\|=1$, when the polynomials are i.i.d. 
realization of a rotationally invariant Gaussian process.

In previous work \cite{montanari2023solving}, we developed an iterative 
algorithm that attempts to minimize
the energy function $H(\bx)=\|\bF(\bx)\|^2/2$ (with $\bF(\bx)=(F_i(\bx):i\le n)$), using first and second order information
(gradient and Hessian of $H$) at a sequence of query points $\bx_1, \ldots, \bx_k$.
Here, we prove a lower bound on the minimum energy achieved by any algorithm
that is Lipschitz continuous in the coefficients of the polynomials  $F_i$. 
We prove that the lower bound
matches the  the energy achieved by the algorithm in \cite{montanari2023solving},
hence characterizing the degree of accuracy to which the system of equations $\bF(\bx)=\bzero$ 
can be solved by Lipschitz algorithms.

Our results imply that,  for certain choices of the distribution of the polynomials $F_i$, 
solutions exist with high probability but Lipschitz algorithms
fails to find near-solutions, even in a very weak sense.
Namely, if $\bx^{\salg}$ is the output of a Lipschitz algorithm, then
$\|\bF(\bx^{\salg})\|^2\ge  \eps\E[\|\bF(\bx_0)\|^2]$ (for $\bx_0$ a 
non-random point on the unit sphere)
with high probability for some $\eps>0$. 
The proof is based on techniques developed in the context of spin glass theory.
\end{abstract}

\tableofcontents

\section{Introduction and main result}

Assume  that $n=\left\lceil \alpha d\right\rceil $ for some fixed
$\alpha\in(0,1]$. Given a polynomial $\xi(t):=\sum_{k=0}^{\kmax}\xi_{k}t^{k}$
with non-negative coefficients,
we let $\bF(\bx)=\bF_{d}(\bx):=(F_{1}(\bx),\ldots,F_{n}(\bx))$ be
a system of i.i.d. random polynomials 
\begin{align*}
	F_{i}(\bx) & :=\sum_{k=0}^{\kmax}\sqrt{\xi_{k}}\sum_{j_{1},\ldots,j_{k}=1}^{d}G_{i,j_{1},\ldots,j_{k}}^{(k)}x_{j_{1}}\cdots x_{j_{k}},
\end{align*}
where the random coefficients $\bG:=(G_{i,j_{1},\ldots,j_{k}}^{(k)})_{i\leq n,k\leq\kmax,j_{1},\ldots,j_{k}\leq d}$
are i.i.d. standard Gaussian random variables.
Equivalently, the $F_{i}(\,\cdot\, )$'s are i.i.d. centered Gaussian 
processes with covariance given by 
\[
\E F_{i}(\bx)F_{j}(\by)=\delta_{ij}\xi(\langle\bx,\by\rangle).
\]

We are interested in algorithms that 
solve the following system of nonlinear equations,
up to a small error:
\begin{align}
\bF(\bx) = \bzero \,,\;\;\;\; \bx\in \S^{d-1}\, ,
\end{align}
where $\S^{d-1} =\{\bx\in\R^d:\; \|\bx\|=1\}$.
Minimizing the $\ell_2$ norm error  is equivalent to minimizing 
the `energy' function 
\begin{equation}
	H(\bx):=\frac{1}{2}\|\bF(\bx)\|^{2}\, .\label{eq:energy}
\end{equation}
We ideally seek $\bx\in\S^{d-1}$ such that $H(\bx)=o(n)$.
Notice in contrast that for a point $\bx\in\S^{d-1}$
chosen independently of $\bF$, we have $H(\bx)=n\xi(1)/2+
o(n)$ with high probability.

In \cite{montanari2023solving}, we developed an 
algorithm for minimization of $H(\bx)$, based on
ideas from mathematical spin glass theory. The algorithm
can be viewed as a `second order method' in the language of optimization theory
\cite{nesterov2003introductory}. Namely it produces a sequence of points 
$\bx_1, \ldots, \bx_k, \dots$, where each new point $\bx_{k+1}$ is
a function of the previous ones as well as of the gradient and Hessian of $H$ 
evaluated at those points.
The analysis of \cite{montanari2023solving}, yields a sharp prediction
on the energy achieved by this algorithm. Namely, it gives an explicit
function $u_{\salg}(\alpha,\xi)$ (also presented in Section \ref{sec:SolOptControl}
below) such that, if $\bx^{\salg}$ is the output of the algorithm, then, for any $\eps>0$,
\begin{align}
\lim_{\substack{n,d\to\infty\\ n/d=\alpha}}\P\left(
\frac{1}{n}H(\bx^{\salg})\le u_{\salg}(\alpha,\xi)+\eps\right)=1\,.
\end{align}

The objective of the current paper is to show that 
this algorithm is optimal within the class of
Lipschitz algorithms, i.e. algorithms whose output 
$\bx^{\salg}$ is a Lipschitz continuous function of the coefficients of the polynomials $F_i$,
with suitably bounded Lipschitz constant.
(In fact, the bound applies to a broader class, as 
explained in Section \ref{sec:Pf_outline_optimalityLip}.)

Any algorithm $\calA_{d}$ defines a mapping $\bG\mapsto\calA_{d}(\bG)=\bx_{\alg}$.
By an abuse of notation we will use $\calA_{d}$ to denote both the algorithm
and the latter mapping. We will always assume that the output belongs
to the unit sphere $\bx_{\alg}\in\S^{d-1}$. We say that an algorithm
$\calA_d$ is $t$-Lipschitz if for any $d$, 
\[
\|\calA_{d}(\bG)-\calA_{d}(\bG')\|_{2}\leq t\|\bG-\bG'\|_{2}
\]
for all $\bG$, $\bG'$, where $\|\bG\|_{2}^{2}$ is the sum of squares
of the elements of $\bG$. We will often also refer to a sequence
$\calA=(\calA_{d})_{d\geq1}$ as an algorithm and write $\calA(\bG)=\calA_{d}(\bG)$
when $d$ is understood from the context.

\begin{definition}\label{def:ControlF}
The set of \emph{control functions} $\cuC$ 
is defined as the set of non-decreasing, continuous
functions $p:[0,1]\to[0,1]$ such that $p(0)=0$, $p(1)=1$ and $p$
is continuously differentiable on $(\theta,1)$, 
where $\theta=\theta(p):=\max\{t:\,p(t)=0\}$.

We further write $p\in \cuC_*\subseteq\cuC$ if there 
exists 
$\chi(q) = \sum_{k=0}^{\infty}\chi_kq^k$, with $\chi_k\ge 0$, absolutely converging for $|q|\le 1$,
such that, for $t\ge \theta = \chi_0$, we have
$p(t)  = \sup\{q\ge 0 :\chi(q)\le t\}$.
\end{definition}
For any $p\in\cuC$ 
 we write $\dot{p}$ (and below $\dot{\xi}$) for the derivative
of $p$ (similarly for $\xi$) and define
the \emph{drift function} 	$\drift: \reals^4\to \reals$ via
\begin{align}
	\drift(t,\uu;p_0,p_1) &:= \max_{s>0}\cL(s;p_0,p_1,\uu,t)\, ,\label{eq:F}\\
	\cL(s;p_0,p_1,\uu,t)&:=-(p_1\dot{\xi}(t)+p_0\ddot{\xi}(t))(\uu\vee0)s+\frac{1}{2}1+p_1\xi(t)-
	\frac{1}{2\alpha s}+\frac{p_0\dot{\xi}(t)}{2(1+p_0\dot{\xi}(t)s)}(1-p_1\xi(t) s)\,.\label{eq:cL}
\end{align}

We let $\uu_{\lb}(t)=\uu_{\lb}(t;p)$ be the
unique solution (see Appendix \ref{app:Properties}) 
of the ordinary differential equation
\begin{equation}
	\frac{\de\uu}{\de t}(t)=\drift(t,\uu;p(t),\dot{p}(t))
\,,\;\;\;\;\uu(\theta)=0\,.\label{eq:ulb}
\end{equation}
Finally, we define 
\begin{equation}
	\ALG=\ALG(\alpha,\xi):=\inf_{p\in\cuC_*}\uu_{\lb}(1;p)\,.\label{eq:ALG}
\end{equation}
The significance of this quantity lies in the fact
that, as shown in Section \ref{sec:SolOptControl},
it coincides with the energy that is achievable by 
the algorithm of \cite{montanari2023solving}. 
Our main result is that $\ALG$ is (asymptotically) 
a lower bound on the energy achievable by any 
Lipschitz algorithm.
\begin{theorem}
	\label{thm:optimalityLip}
	If $\calA$ is $M/\sqrt{d}$-Lipschitz for
	some fixed $M>0$, then for any $\eps>0$ there exists a constant
	$C=C(\eps,\xi,\alpha,M)>0$ such that 
	\[
	\P\left(\frac{1}{n}H(\calA(\bF))\ge\ALG-\eps\right)\ge1-e^{-Cn}/C\,.
	\]
\end{theorem}
In the next section, we will define the function $u_{\salg}(\alpha,\xi)$ 
that characterizes the energy achieved by the algorithm of \cite{montanari2023solving},
and state a theorem establishing that this matches the above lower bound,
namely $u_{\salg}(\alpha,\xi)=\ALG(\alpha,\xi)$.
In Section \ref{sec:Pf_outline_optimalityLip}, we provide an outline of the proof of Theorem \ref{thm:optimalityLip},
while the technical details are presented in Sections \ref{sec:Concentration},
\ref{sec:ProofMainBound}, and \ref{sec:MainProof}. Finally, some technical estimates and 
proofs for Theorem \ref{thm:Equivalence} below are presented in the appendices.

\subsection*{Notations}

We will denote by $\Ball^d(\bx_0;r)$ the Euclidean ball in $\R^d$ with center $\bx_0$ and
radius $r$, and write $\Ball^d(r)$ when the center is $\bx_0=\bzero$. 
Given random variables (or vectors) $X,Y$ we write  $X\ed Y$ if they have the same law.  

In the following, we will often omit the argument from
$p(t)$, $\dot{p}(t)$, $\xi(t)$, $\ddot{\xi}(t)$, 
when it is clear from the context.
We use either dot (as in $\dot{p}(t)$) or prime 
(as in $p'(t)$) to denote the derivative.

\section{Solution of the optimal control problem}
\label{sec:SolOptControl}
In this section we define the energy threshold $u_{\salg}(\alpha,\xi)$ 
that is achieved by the algorithm of \cite{montanari2023solving}
(as per Theorem 1 there).
We will then give a formal statemet that the value 
$u_{\salg}(\alpha,\xi)=\ALG(\xi,\alpha)$, where $\ALG(\xi,\alpha)$ is
the lower bound defined in the previous section.
This establishes that the algorithm of \cite{montanari2023solving}
is optimal within the class of Lipschitz algorithms.

To simplify the presentation, we assume
$\xi(0)>0$ strictly.
There is no real loss of generality in this, 
because we can always replace $\bF(\bx)$ by 
$\bF_{\eps_0}(\bx) :=\eps_0\bg+\bF(\bx)$ 
with $\eps_0>0$ and $\bg\sim \normal(\bzero,\bI_n)$.
 This only results in an $O(n\eps_0)$ change in
the energy.

Define the `replica symmetry breaking radius' 
\begin{align}\label{eq:qrsDef}
q_{\sRS}(\xi):= \argmin_{q\ge 0} \frac{\xi(q)\xi'(q)}{q} \, .
\end{align}
Here we set by convention $q_{\sRS}(\xi)=\infty$ if 
$\xi(t) =\xi_0+\xi_1 t$, and it is understood that the
function  $g(q):=\xi(q)\xi'(q)/q$ is extended by continuity 
at $q=0$ (with $g(0)=\infty$ if $\xi(0)\xi'(0)>0$).

Uniqueness of the minimizer holds because, letting
$f(q) := \xi(q)\xi'(q)/q$, we have $f(q) = c\, q^{-1}
+f_+(q)$ for $c\ge 0$ and $f_+$ a non-zero polynomial with 
strictly positive coefficients (assuming $\xi$ is not affine).

We further define $q_0=q_0(\alpha,\xi)$  as the solution 
of 
\begin{align}\label{eq:q0Def}
\frac{q_0\xi'(q_0)}{\xi(q_0)} = \alpha\, .
\end{align}
(The solution exists uniquely, unless 
$\xi(t) =\xi_0+\xi_1 t$ and $\alpha\ge 1$, in which case we set
$q_0(\alpha,\xi)=\infty$.)

Next, we define the `replica symmetric' energy as
\begin{align}\label{eq:uRS_def}
u_{\sRS}(q,\alpha,\xi) := 
\frac{1}{2}\Big(\sqrt{\xi(q)}-\sqrt{\frac{1}{\alpha}q\xi'(q)}\Big)_+^2\, ,
\end{align}
and note that $q_0$ is the first point at which 
$q\mapsto u_{\sRS}(q,\alpha,\xi)$ vanishes.

For $a,b\in\reals_{\ge 0}$, define
$z_*(\alpha,a,b)$  via
\begin{align}
 Q(m;\alpha,a,b) &:= -\frac{1}{m}+\frac{\alpha b}{1+bm}-a^2m\, ,\\
 z_*(\alpha,a,b) & := -\sup_{m>0}  Q(m;\alpha,a,b)\, .
 \end{align}
 %
We distinguish two cases.

\paragraph{Case 1: $q_{\sRS}(\xi)\ge 1$.} We set 
 $u_{\sRS}(\alpha,\xi):=u_{\sRS}(1,\alpha,\xi)$ and
\begin{align}
	u_{\salg}(\alpha,\xi) := u_{\sRS}(\alpha,\xi)=
	\frac{1}{2}\left(\sqrt{\xi(1)}-\sqrt{\frac{\xi'(1)}{\alpha}}
	\right)_+^2\, .\label{eq:REenergy}
\end{align}

\paragraph{Case 2: $q_{\sRS}(\xi)< 1$.} 
 Let $u(t) = u(t ;\alpha,\xi,q)$,  be the unique solution of 
 the ordinary differential equation 
 \begin{align}
 \frac{\de u}{\de t}(t) = -\frac{1}{2\alpha}z_*\big(\alpha; \sqrt{2\alpha u(t) \xi''(q+t)}, \xi'(q+t)\big)\,,
 \;\;\;\;\; u(0) = u_{\sRS}(q,\alpha,\xi)\, .\label{eq:U_ODE}
 \end{align}
Finally we define the algorithmic threshold energy as
\begin{align}\label{eq:AchievableEnergyDef}
 u_{\salg}(\alpha,\xi)= \begin{cases}
	u(1-q_{\sRS};\alpha,\xi,q_{\sRS}) & \mbox{if $q_{\sRS}(\xi)
	<q_0(\alpha,\xi)$,}\\
	0 & \mbox{otherwise.}
 \end{cases}
\end{align}
In particular, $q_0(\alpha,\xi)\le q_{\sRS}(\xi)$ then the above 
yields $u_{\salg}(\alpha,\xi)= 0$.
The next theorem is proven in Appendix \ref{app:Equivalence}.
\begin{theorem}\label{thm:Equivalence}
	With the above definitions, 
	\begin{align*}
	\ALG(\alpha,\xi)=u_{\salg}(\alpha,\xi)\, .
	\end{align*}
\end{theorem}

\section{\label{sec:Pf_outline_optimalityLip}Proof outline for Theorem \ref{thm:optimalityLip}}

Our proof of  Theorem \ref{thm:optimalityLip} builds on a technique
developed by Huang and Sellke in \cite{huang2025tight,huang2024optimization}
in the context of spin glass models. We introduce the main ingredients of the proof in this section.

We say that two jointly Gaussian copies $(\bF^{1},\bF^{2})$ of $\bF$ are $\rho$-correlated if 
\begin{equation}
	\E\big\{ F_{i}^{1}(\bx)F_{j}^{2}(\by)\big\}=\rho\, \E\big\{ F_{i}^{1}(\bx)F_{j}^{1}(\by)\big\}=
	\rho\, \E\big\{ F_{i}^{2}(\bx)F_{j}^{2}(\by)\big\}.\label{eq:rhocorrelation}
\end{equation}
The following notion was introduced in \cite{huang2025tight}.
\begin{definition}
	\label{def:overlap_concentration}Let $\lambda,\nu>0$. An algorithm
	$\calA=(\calA_{d})_{d\geq1}$ is said to be $(\lambda,\nu)$ overlap
	concentrated if for any $d\geq1$, $\rho\in[0,1]$ and $\rho$-correlated
	$(\bF^{1},\bF^{2})$,
	\[
	\P\left(\big|\langle\calA(\bF^{1}),\calA(\bF^{2})\rangle-\E\langle\calA(\bF^{1}),\calA(\bF^{2})\rangle\big|\geq\lambda\right)\leq\nu.
	\]
We will also use the same definition for $\calA_{d}$ for a given $d\ge1$ with the obvious modification.
\end{definition}

We will prove the following theorem about optimality of such algorithms.
\begin{theorem}\label{thm:optimalityOC}
	For any $\xi$, $\alpha$ and $\delta>0$
	there exist constants $c,\lambda>0$ such that if $\calA$ is $(\lambda,\nu)$
	overlap concentrated for some $\nu\in[0,1]$, then 
	\begin{equation}\label{eq:Thm_OC}
	\P\left(\frac{1}{n}H(\calA(\bF))\le\ALG-\delta\right)\le2e^{-cd}+2\nu^{c}\,.
	\end{equation}
\end{theorem}

By an application of Gaussian concentration, $(M/\sqrt{d})$-Lipschitz
algorithms are $(\lambda,\nu)$ overlap concentrated, for arbitrary
$\lambda$ and $\nu=2e^{-\frac{d\lambda^{2}}{8M}}$
 (see Lemma \ref{lem:overlapConc} below). We use this in 
 Section \ref{subsec:pf_optimalityLip} to conclude Theorem \ref{thm:optimalityLip} from Theorem \ref{thm:optimalityOC}, and thus focus on the latter of the two for the rest of the current section.

Let $\kappa\ge1$ be an integer and $\eps>0$ a real number. Let $0\le q<q'\le1$
and $0\le\rho<\rho'\le1$. 
For $i\ge 0$, let $\bF^{(i)}$ be jointly Gaussian polynomials with covariance 
\begin{equation}
	\E F_{\ell}^{(i)}(\bx)F_{\ell'}^{(j)}(\by)=\delta_{\ell\ell'}\big[\rho+\delta_{ij}(1-\delta_{i0})(\rho'-\rho)\big]\xi(\langle\bx,\by\rangle)\label{eq:correlation}
\end{equation}
and define $H^{(i)}(\bx):=\frac{1}{2}\|\bF^{(i)}(\bx)\|^{2}$.
Define 
\begin{align}
	T(\bx)&:=  T(\bx,\kappa,\eps,q,q'):= \Big\{(\bx^{1},\ldots,\bx^{\kappa})\in(\Ball^{d}(1))^{\kappa}\;:\;\,\forall i\neq j\le\kappa,\label{eq:T(x)}\\
	& \hphantom{:=  T(\bx,\kappa,\eps,q,q'):= \Big\{}\big|\langle\bx^{i},\bx\rangle-q\big|\vee\big|\langle\bx^{i},\bx^{j}\rangle-q\big|\le\eps,\,\langle\bx^{i},\bx^{i}\rangle=q'\Big\}\,,\nonumber \\
	\Delta(\bx)&:=  \Delta(\bx,\kappa,\eps,q,q',\rho,\rho'):= \frac{1}{\kappa}\min_{(\bx^i)_{i\le \kappa}\in T(\bx)}\sum_{i=1}^{\kappa}H^{(i)}(\bx^{i})\,.\label{eq:DeltaDef}
\end{align}
We will usually use the latter with $\bx$ such that $\|\bx\|^2_2=q$.  The key estimate we use in the proof of Theorem \ref{thm:optimalityOC} is the following
proposition, which we prove in Section \ref{sec:ProofMainBound}.


%
	\begin{proposition}
	\label{prop:main_bound}
	There exists a sequence $A_d=A_d(\alpha)\overset{d\to\infty}{\longrightarrow}0$ and a 
	constant $C_{0}=C_{0}(\xi,\alpha,C)$ such that the following holds. If $q<q'$ are values in
	 $[0,1]$ such that either $q=0$ or $q'/q\le 4$, $p:[0,1]\to[0,1]$ is a function that is continuous, 
	 concave and non-decreasing on $[q,q']$ and twice-differentiable on $(q,q')$ with 
	$\sup_{t\in(q,q')}|\ddot{p}(t)|<C$ and 
	$\phi_0\in[0,C]$, then
	letting $\rho=p(q)$ and $\rho'=p(q')$, for $\bx$ with $\|\bx\|^2=q$
	we have that (for almost every $\phi_0$) 
	\begin{align*}
		\E\left[\frac{1}{n}\Delta(\bx)\,\bigg|\,\frac{1}{n}H^{(0)}(\bx)=\phi_0\right]&\ge\phi_0+
		(q'-q) \drift(q,\phi_0;p,\dot{p})\\
		&\hphantom{\ge\phi_0}-C_{0}\Big[(q'-q)^{\frac{3}{2}}+A_d(q'-q)+\frac{1}{\sqrt{\kappa}}(q'-q)^{\frac{1}{2}}+\kappa\sqrt{\eps}\Big]\,.
	\end{align*}
\end{proposition}	
We note that above, $\uu_0$ is a real number and not a function, in which case, we define $\drift(q,\phi_0;p,\dot{p})$ through 
\eqref{eq:F}  with the constant function $\uu\equiv \uu_0$ and with $\dot{p}(q)$ taken as the derivative from the right.
Observe the connection to \eqref{eq:ulb} for small $q'-q$.

\subsection{Hierarchically correlated systems}
\label{sec:Hierarchical}

In analogy with \cite{huang2025tight,huang2024optimization} 
we will construct a collection of random polynomial
systems with a hierarchical correlation structure
(see Eq.~\eqref{eq:bG(u)}), and associate to them a 
tree of points in $\R^d$ (see Eq.~\eqref{eq:xu}). An analogous construction to \eqref{eq:xu} was used in \cite{subag2024free} 
for spherical
spin glasses, to define an embedding of the ultrametric
tree of pure and ancestral states.

The construction depends on a choice of a tree $\TT$
and a variance profile $p$. 
Here, $p=p_d$ can be any non-decreasing $[0,1]\to[0,1]$ function with $p(0)=0$. We take $\theta=\theta_{d}=\max\{t:\,p(t)=0\}$
as per Definition \ref{def:ControlF}.
Given integers $\kappa,\tau\geq1$, we let $\TT:=\mathbb{T}(\kappa,\tau)$
be the complete $\kappa$-ary tree of depth $\tau$ with root $\r$, depth
$a$ vertex set $V_{a}=[\kappa]^{a}$ and set of leaves $\LL:=\LL(\kappa,\tau)=V_{\tau}$.
For any two vertices $u,v\in\TT$ let $|u|$ denote the depth of $u$,
$\r\leftrightarrow u$ the path from the root to $u$ and $u\wedge v$
the least common ancestor. We write $u\leq v$ if $v$ is a descendant
of $u$ or $u=v$, and denote by $D(u):=\{v\in \TT:\; u\le v, |v|=|u|+1 \}$
the set of direct descendants (`children') of $u$.
We define 
\begin{align}
\< u,v\>:=\theta+\frac{(1-\theta)}{\tau}|u\wedge v|\, .\label{eq:TreeMetric}
\end{align}
Note that $\langle u,v\rangle$ implicitly depends on $p$ (and therefore
on $d$).

For any edge $e=(u,v)\in\mathbb{T}$, let $\bF^{e}=\bF_{d}^{e}$ be
an i.i.d. copy of $\bF$ and, by an abuse of notation, write $p(e):=|p(\langle v,v\rangle)-p(\langle u,u\rangle)|$.
For any vertex $u\in\TT$, define the polynomial system
\begin{equation}
	\begin{aligned}\bF^{(u)}(\bx) & :=\bF_{d}^{(u)}(\bx):=(F_{1}^{(u)}(\bx),\ldots,F_{n}^{(u)}(\bx))\end{aligned}
	\label{eq:Fu}
\end{equation}
by 
\begin{equation}
	\bF^{(u)}:=\sum_{e\in\r\leftrightarrow u}\sqrt{p(e)}\bF^{e},\label{eq:bG(u)}
\end{equation}
so that
\begin{equation}
	\E F_{i}^{(u)}(\bx)F_{j}^{(v)}(\by)=\delta_{ij}p(\langle u,v\rangle)\xi(\langle\bx,\by\rangle).\label{eq:correlation-1}
\end{equation}
Note that, in the terminology above, $\bF^{(u)}$ and $\bF^{(v)}$
are $\rho$-correlated for $\rho=p(\langle u,v\rangle)$, see Eq.~\eqref{eq:rhocorrelation}. 

The correlated polynomial systems (\ref{eq:bG(u)})
can be defined with any profile $p(t)$ as above. 
In the rest of the section, however, we
will use a specific profile that we now define. Given an algorithm $\calA$,
for $\rho\in[0,1]$ and $\rho$-correlated $(\bF^{1},\bF^{2})$ we
define 
\[
\chi(\rho)=\chi_{d}(\rho):=\E\left\langle \calA(\bF^{1}),\calA(\bF^{2})\right\rangle .
\]
If $\chi_{d}$ is constant then, since $\chi_{d}(1)=1$ by definition,
$\chi_{d}\equiv1$ and therefore $\calA_{d}$ is deterministic (i.e.,
it a.s. outputs the same fixed vector). In this case, $H(\calA(\bF))=\frac{1}{2}n\xi(1)+o(n)$
with high probability. We assume henceforth that $\chi_{d}$ is not
constant. 

By expanding $ \calA(\bF^{\ell})$ into the Hermite basis of the 
underlying  Gaussian random variables  $\bG^{\ell} = (G^{\ell(k)}_{i,j_1,\dots,j_k})$,
we obtain that 
$\chi_{d}(\rho)=\sum_{i=0}^{\infty}\alpha_{d,i}\rho^{i}$
for some non-negative coefficients $\alpha_{d,i}$ 
\cite[Proposition 3.1]{huang2024optimization}. In particular, $\chi_d$ is continuous,
convex and strictly increasing on $[0,1]$. 
For the next construction (in \eqref{eq:xu}), we  suppose
that $\theta=\theta_{d}:=\chi_{d}(0)$ and
\begin{equation}
	p(t)=p_{d}(t):=\begin{cases}
		0 & t\in[0,\theta)\\
		\chi_{d}^{-1}(t) & t\in[\theta,1],
	\end{cases}\label{eq:p_linear}
\end{equation}
so that, by definition, for any $u,v\in \TT$, 
\[
\E\big\langle\calA(\bF_{d}^{(u)}),\calA(\bF_{d}^{(v)})\big\rangle=\chi_{d}(p_{d}(\langle u,v\rangle))=\langle u,v\rangle.
\]

To any leaf of the tree $u\in\LL$, we associate the point $\bx_{u}:=\calA(\bF_{d}^{(u)})\in\S^{d-1}$.
To other vertices $u\in\TT\setminus\LL$ we associate the normalized
averages 
\begin{equation}
	\bx_{u}:=\sqrt{\langle u,u\rangle}\frac{\hat{\bx}_{u}}{\|\hat{\bx}_{u}\|},\quad\hat{\bx}_{u}:=|\LL(u)|^{-1} \sum_{v\in\LL(u)}\bx_{v},\label{eq:xu}
\end{equation}
where $\LL(u):=\{v\in\LL:u\leq v\}$ is the set of leaves in the subtree
of $u$. If $\hat{\bx}_{u}=\bzero$, to avoid division by zero, we let
$\bx_{u}$ be some arbitrary deterministic vector with norm $\|\bx_{u}\|=\sqrt{\langle u,u\rangle}$.
Note that to define $\bx_{u}$ for all vertices $t\in\TT$ we only
use the outputs for $\bF^{(u)}$ associated to leaves $u\in\LL$.

The overlap concentration property of the algorithm implies the following about 
 the geometry of the  tree $(\bx_{u})_{u\in\TT}\subset\Ball_{d}(1)$. 
\begin{lemma}
	\label{lem:tree_overlap}
	Suppose that $\lambda<1/(\kappa\tau)$, $\kappa>1$,
	$p(t)$ is defined by (\ref{eq:p_linear}) and $\calA$ is $(\lambda,\nu)$
	overlap concentrated. Then, denoting $\eta_u =\sqrt{\<u,u\>}$ and $\eta_{u,v}=\eta_{u}\wedge\eta_{v}$,
	\begin{equation}
		\P\bigg(\forall u,v\in\TT:\:\big|\left\langle \bx_{u},\bx_{v}\right\rangle -\langle u,v\rangle\big|\leq \eta_{u,v} \wedge \frac{20}{\kappa\tau\eta_{u,v}}\bigg)\geq1-\kappa^{2\tau}\nu.\label{eq:overlapBd_explicit}
	\end{equation}
	Consequently, if in addition to the above we assume that $20/(\kappa\tau\eps)\le \eps$ for some $\eps>0$, then 
	\begin{equation}
		\P\Big(\forall u,v\in\TT:\:\big|\left\langle \bx_{u},\bx_{v}\right\rangle -\langle u,v\rangle\big|\leq \eps\Big)\geq1-\kappa^{2\tau}\nu.\label{eq:overlapBd}
	\end{equation}
\end{lemma}

\begin{proof}
	Denote $A_{u,v}(t):=\langle u,v\rangle+[-t,t]$. Let $\calE$ be the
	event that for all pairs of leaves $u,v\in\LL$, $\langle\bx_{u},\bx_{v}\rangle\in A_{u,v}(\lambda)$.
	From our assumptions and a union bound, $\P(\calE)\geq1-\kappa^{2\tau}\nu$.
	To complete the proof, assuming henceforth that $\calE$ occurs, we will show that the event  in \eqref{eq:overlapBd_explicit} also occurs.

	For any $u,v\in\TT$, $\langle\hat{\bx}_{u},\hat{\bx}_{v}\rangle$
	is the average of $\langle\bx_{u'},\bx_{v'}\rangle$ over all pairs
	$(u',v')\in\LL(u)\times\LL(v)$. If neither $u\leq v$ nor $v\leq u$,
	then $\langle u',v'\rangle=\langle u,v\rangle$ for all such leaves.
	If on the other hand $u\leq v$ or $v\leq u$,
	then for any natural $i$, a proportion of $(1/\kappa)^i$ at most of the pairs $(u',v')\in\LL(u)\times\LL(v)$, $\<u',v'\>=\<u,v\>+i(1-\theta)/\tau$ and $\<{\bx}_{u'},{\bx}_{v'}\>\in A_{u,v}(\lambda+i/\tau)$. For all other pairs, that do not satisfy the latter for all $i\ge1$, we have $\<u',v'\>=\<u,v\>$ and $\<{\bx}_{u'},{\bx}_{v'}\>\in A_{u,v}(\lambda)$. Therefore,
	\[
	|\<\hat{\bx}_{u},\hat{\bx}_{v}\> -\<u,v\>|\le \lambda + \sum_{i=1}^{\infty} (i/\tau)(1/\kappa)^i = \lambda + \frac{\kappa}{\tau(\kappa-1)^2}=:B\,.
	\]

	We also have that
	\begin{align*}
		|\<\bx_{u},\bx_{v}\>-\<\hat{\bx}_{u},\hat{\bx}_{v}\>| & = \eta_u\eta_v \frac{|\<\hat{\bx}_{u},\hat{\bx}_{v}\>|}{\|\hat{\bx}_{u}\|\|\hat{\bx}_{v}\|}\Big(1-\frac{\|\hat{\bx}_{u}\|\|\hat{\bx}_{v}\|}{\eta_u\eta_v}\Big)\\
		&\le \eta_u\eta_v \sup_{|t|\leq B}\bigg|1-\Big|1+\frac{t}{\eta_u^2}\Big|^{1/2}\Big|1+\frac{t}{\eta_v^2}\Big|^{1/2}\bigg|\\
		&\le \eta_{u}\eta_{v}\Big(\frac{B}{\eta_u^2}+\frac{B}{\eta_v^2}+\frac{B^2}{\eta_u^2\eta_v^2}\Big)\le \frac{3B}{\eta_{u,v}}\,,
	\end{align*} 
where we used that $1-\sqrt{1-x}\le|x|$ for $x\ge-1$ and that $B/\eta_{u,v}\le 1$, otherwise the inequality holds trivially.
	Moreover, 
	\[
	|\langle\bx_{u},\bx_{v}\rangle|\le\eta_{u}\eta_{v}\quad\text{and}\quad|\langle u,v\rangle|\leq\eta_{u}\eta_{v}.
	\]

	Combining the above, we have that 
	\[
	|\<\bx_{u},\bx_{v}\>-\<u,v\>|\le \eta_{u}\eta_{v}\wedge \frac{4B}{\eta_{u,v}} \le \eta_{u,v} \wedge \frac{20}{\kappa\tau\eta_{u,v}}\,.\qedhere
	\]
\end{proof}

\subsection{\label{subsec:energy decrements}Energies on the tree and
	concentration}

In the previous section, we constructed the random tree $(\bx_{u})_{u\in\TT}$
from the outputs of the algorithm $\cA$ when applied to the correlated systems $(\bF^{(u)}:\; u\in\LL)$. 
On this tree,
we consider the energies $H^{(u)}(\bx_{u}):=\|\bF^{(u)}(\bx_{u})\|^{2}/2$.
Overloading the notations $T(\bx)$ and $\Delta(\bx)$ introduced in
Eqs.~\eqref{eq:T(x)}, \eqref{eq:DeltaDef}, for $u\in\TT\setminus\LL$,
$\by\in\Ball^d(1)$, 
we define 
\begin{equation}
	\begin{aligned}
	T(\by,u)&:=  T(\by,u,\kappa,\tau,\eps) := \Big\{(\by_{v})_{v\in D(u)} \in (\Ball^d(1))^{\kappa}:\,\forall v\neq v',\\
		& \hphantom{:= T(\by,u,\kappa,\tau,\eps) := \Big\{} \big|\langle\by_{v},\by\rangle-\langle u,u\rangle\big|\vee\big|\langle\by_{v},\by_{v'}\rangle-\langle u,u\rangle\big|\le\eps,\,\langle\by_{v},\by_{v}\rangle=\langle v,v\rangle\Big\}\,,\\
		\Delta(\by,u)&:= {\Delta}(\by,u,\kappa,\eps):= \frac{1}{\kappa}\inf_{ (\by_{v})_{v\in D(u)} 
		\in T(\by,u)}\sum_{v\in D(u)}H^{(v)}(\by_{v})\,.
	\end{aligned}
	\label{eq:barT+barDelta}
\end{equation}
Note that, if the
event in Eq.~\eqref{eq:overlapBd} occurs, then $(\bx_v)_{v\in D(u)}\in T(\bx_{u},u)$ and therefore
\begin{equation}
	\frac{1}{\kappa}\sum_{v\in D(u)}H^{(v)}(\bx_{v})\ge \Delta(\bx_{u},u).\label{eq:rel_to_barDelta}
\end{equation}

Recalling Eqs.~\eqref{eq:TreeMetric} and \eqref{eq:correlation-1}, we
 see that for $q=\langle u,u\rangle$,
$q'=q+(1-\theta)/\tau$, $\rho=p(q)$ and $\rho'=p(q')$, we have
that $T(\bx,u)=T(\bx)$ for any $\bx$ and 
\begin{equation}
	\left(H^{(u)}(\bx),\Delta(\bx,u)\right)\ed \left(H^{(0)}(\bx),\Delta(\bx)\right),\label{eq:Deltaequiv}
\end{equation}
as processes on $\sqrt{q}\S^{d-1}=\{\bx:\:\|\bx\|^{2}=q\}$. In Section \ref{sec:Concentration} we prove a concentration result for the deviations of $\Delta(\bx)$ from $\E\big[\Delta(\bx)\,|\,H^{(0)}(\bx)\big]$ uniformly over $\sqrt{q}\S^{d-1}$ (see Proposition \ref{prop:unif_conc}), by utilizing a key idea from \cite{subag2024free,chen2018generalized,chen2020generalized2} (used there in the  development
of a generalized Thouless-Anderson-Palmer approach 
to spin glasses). In the same section, we conclude the following.

\begin{theorem}[Uniform concentration]
	\label{thm:concentration}For any $\delta_{0},c>0$, there exist
	$\eps_{0}$, $\kappa_{0}$, $\tau_0$ and $d_{0}$ which depend on $\delta_{0}$,
	$c$ and $\xi$, such that for any non-decreasing $p:[0,1]\to[0,1]$ with $p(0)=0$, 
	$\eps\leq\eps_{0}$, $\kappa\geq\kappa_0$, $\tau\ge \tau_0$ and $d\geq d_{0}$,
	\begin{equation}
		\P\left(\max_{u\in\TT\setminus(\LL\cup\{\r\})}\max_{\bx:\|\bx\|^2=\<u,u\>}\Big|\Delta(\bx,u)-\Psi\Big(\frac{1}{n}H^{(u)}(\bx),u\Big)\Big|>\delta_{0}d\right)<\kappa^{\tau}e^{-cd},\label{eq:Delta_bd_uniform-6}
	\end{equation}
	where, with $q=\langle u,u\rangle$, $q'=q+(1-\theta)/\tau$, $\rho=p(q)$
	and $\rho'=p(q')$ and $\bx$ such that $\|\bx\|^{2}=q$, 
	\begin{equation}
		\Psi(\uu_0,u):=\E\Big[\Delta(\bx,\kappa,\eps,q,q',\rho,\rho')\,\Big|\,\frac{1}{n}H^{(0)}(\bx)=\uu_0\Big].\label{eq:Psi-1}
	\end{equation}
\end{theorem}

\subsection{Putting things together}
\label{subsec:pf of reduction}

Recall the constructions from Section \ref{sec:Hierarchical} parameterized by the tree $\TT$ of degree $\kappa$ and depth $\tau$. Also recall the definition of $\theta=\theta_{d}\in[0,1)$ and $p(t)=p_{d}(t)$
in (\ref{eq:p_linear}) and that $\chi_{d}(\rho)=\sum_{i=0}^{\infty}\alpha_{d,i}\rho^{i}$
for some $\alpha_{d,i}\ge0$. Since we assume that $\chi_{d}(\rho)$
is not a constant function, we have $p\in\cuC_*$. In particular, $p(t)$ is equal to $0$ on $[0,\theta]$
and is continuous, concave  and strictly increasing on $[\theta,1]$.
Let $\uu_{\lb}(t)=\uu_{\lb}(t;p)$ be the solution of Eq.~\eqref{eq:ulb}
on $[\theta,1]$.

We have the following straightforward adaptation of Proposition
3.6(a) of \cite{huang2024optimization}.
\begin{lemma}
	\label{lem:many_to_one_exceedence_prob}For any $t\geq0$,
	\begin{equation}
		\P\left(\forall u\in\LL:\,H^{(u)}(\calA(\bF^{(u)}))\leq t\right)\geq\P\left(H(\calA(\bF))\leq t\right)^{\kappa^{\tau}}.\label{eq:many_to_one_exceedence_prob}
	\end{equation}
\end{lemma}

\begin{proof}
	Suppose $u\leq u'$ are some vertices in $\TT$ with $|u'|=|u|+1\leq\tau$.
	Let $\bG^{(v)}$ be the set of the Gaussian coefficients of $\bF^{(v)}$.
	Then, by symmetry, 
	\begin{align*}
		& \P\left(\forall v\in\LL(u):\,H^{(v)}(\calA(\bF^{(v)}))\le t\,\Big|\,\bG^{(w)},\,\forall w: w\le u \right)\\
		& =\P\left(\forall v\in\LL(u'):\,H^{(v)}(\calA(\bF^{(v)}))\le t\,\Big|\,\bG^{(w)},\,\forall w: w\le u \right)^{\kappa}.
	\end{align*}
	By Jensen's inequality and by induction, if $u\leq u'$ and $|u'|=|u|+a\leq\tau$,
	\begin{equation}
		\P\left(\forall v\in\LL(u):\,H^{(v)}(\calA(\bF^{(v)}))\le t\right)\geq
		\P\left(\forall v\in\LL(u'):\,H^{(v)}(\calA(\bF^{(v)}))\le t\right)^{\kappa^{a}}.\label{eq:exceedence_prob_induction}
	\end{equation}
	The lemma follows, since (\ref{eq:many_to_one_exceedence_prob}) coincides
	with (\ref{eq:exceedence_prob_induction}) for $u=\r$ and
	arbitrary $u'\in\LL$.
\end{proof}

We next to explain the general strategy to prove Theorem \ref{thm:optimalityOC}:
we refer the reader to  Section \ref{sec:MainProof}
for the precise quantive bounds we 
need.
By the last lemma, it is sufficient to prove that 
 with overwhelmingly high probability,
 \begin{equation}
	\exists u\in\LL:\quad \frac{1}{n}H^{(u)}(\calA(\bF^{(u)}))\ge \ALG-\delta\,.\label{eq:bd_leaves}
\end{equation} 
In the previous subsections, we discussed several estimates that hold with high probability. To explain the basic idea, next we will see how \eqref{eq:bd_leaves} follows if we assume exact equalities instead of some those estimates. 
First,  suppose that the differences as in \eqref{eq:Delta_bd_uniform-6} are exactly zero and for simplicity assume that the equality holds also for $u=\r$,
\[
\forall u\in\TT\setminus\LL:\quad \Delta(\bx_{u},u)=\Psi\Big(\frac{1}{n}H^{(u)}(\bx_{u}),u\Big)\,.
\] 
Moreover, for $\uu_0 = \uu_{\lb}(q)$, note that neglecting the small error terms, the right-hand side in the bound of Proposition \ref{prop:main_bound} is equal to 
\begin{equation}\label{eq:uu_lb_lin_apx}
	\uu_{\lb}(q)+(q'-q)\dot\uu_{\lb}(q)\approx \uu_{\lb}(q')\,.
\end{equation}
Here, we neglect the error as suppose that
\begin{equation*}
		\frac{1}{n}\Psi\big(\uu_{\lb}(q) ,u\big)  = \E\left[\frac{1}{n}\Delta(\bx)\,\bigg|\,\frac{1}{n}H^{(0)}(\bx)=\uu_{\lb}(q)\right]\ge \uu_{\lb}(q')\,,
\end{equation*}
with $q=\<u,u\>$ and $q'=\<v,v\>$ for any $u\in\TT$ and $v\in D(u)$. Finally, suppose that the bound on overlap differences as in \eqref{eq:overlapBd} holds so that we can invoke \eqref{eq:rel_to_barDelta}.

If $\frac{1}{n}H^{(u)}(\bx_{u}) \ge \uu_{\lb}(\<u,u\>)$ then, assuming the above and using the monotoncity of $\Psi(\cdot,u)$,
\[
\frac{1}{n}\Delta(\bx_{u},u) \ge \frac{1}{n}\Psi\big(\uu_{\lb}(\<u,u\>) ,u\big) \ge  \uu_{\lb}(\<u,u\>+\sigma)\,,
\]
where we denote $\sigma=(1-\theta)/\tau$. From \eqref{eq:rel_to_barDelta}, there is at least one $v\in D(u)$ such that  $\frac{1}{n}H^{(v)}(\bx_{v})\ge \uu_{\lb}(\<u,u\>+\sigma)=\uu_{\lb}(\<v,v\>)$.
Now,  note that by definition $\frac{1}{n}H^{(\r)}(\bx_{\r}) =\uu_{\lb}(\theta)=0$. Hence, by induction, we may construct a path $\r=v_0,v_1,\ldots,v_{\tau}\in\LL$ such that $\frac{1}{n}H^{(v_i)}(\bx_{v_i})\ge \uu_{\lb}(\<v_i,v_i\>)$. In particular, \eqref{eq:bd_leaves} holds since for the leaf,
\[
\frac{1}{n}H^{(v_{\tau})}(\calA(\bF^{(v_{\tau})})) = \frac{1}{n}H^{(v_{\tau})}(\bx_{v_\tau})\ge \uu_{\lb}(1) \ge \ALG\,.
\]

\section{Uniform concentration: Proof of Theorem
	\ref{thm:concentration}}\label{sec:Concentration}

Recall the definition of $\bF^{(i)}(\bx)$ and $T(\bx)=T(\bx,\kappa,\eps,q,q')$
in (\ref{eq:correlation}) and (\ref{eq:T(x)}). We start with the
following bound on the conditional variance.
\begin{lemma}
	\label{lem:conditional_cov_bd}Suppose that $\|\bx\|^{2}=q<q'\leq1$,
	$0\leq\rho\leq\rho'\leq1$ and $(\bx^{1},\ldots,\bx^{\kappa})\in T(\bx)$.
	Then, for any  vectors $\bv^{i}\in\R^{n}$, $\bt\in\R^n$, and $s_{i}\in\R$ with
	$\|\bv^{i}\|=1$ and $\sum_{i\leq\kappa}s_{i}^{2}=1$,
	\[
	{\rm Var}\bigg(\sum_{i\leq\kappa}\langle s_{i}\bF^{(i)}(\bx^{i}),\bv^{i}\rangle\,\Big|\,\bF^{(0)}(\bx)=\bt\bigg)\leq4\xi(1)+3\kappa\eps\xi'(1).
	\]
\end{lemma}

\begin{proof}
	Below we abbreviate ${\rm Var}_{\bt}(\cdot):={\rm Var}(\,\cdot\,|\,\bF^{(0)}(\bx)=\bt)$.
	Recall that 
	\begin{align}\label{eq:FhatF}
		\bF^{(i)}(\bx)=\sqrt{\rho}\hat{\bF}^{(0)}(\bx)+(1-\delta_{i0})\sqrt{\rho'-\rho}\hat{\bF}^{(i)}(\bx)
	\end{align}
	for i.i.d. copies $\hat{\bF}^{(i)}(\,\cdot\,)$ of $\bF(\,\cdot\,)$. Hence, 
	\begin{align*}
		{\rm Var}_{\bt}\bigg(\sum_{i\leq\kappa}s_{i}\langle\bF^{(i)}(\bx^{i}),\bv^{i}\rangle\bigg) & =\rho{\rm Var}_{\bt}\bigg(\sum_{i\leq\kappa}s_{i}\langle\hat{\bF}^{(0)}(\bx^{i}),\bv^{i}\rangle\bigg)+(\rho'-\rho){\rm Var}\bigg(\sum_{i\leq\kappa}s_{i}\langle\hat{\bF}^{(i)}(\bx^{i}),\bv^{i}\rangle\bigg).
	\end{align*}
	
	By independence, the second variance above is bounded by $\xi(1)$. For the first term, we have 
	\begin{align*}
	{\rm Var}_{\bt}\bigg(\sum_{i\leq\kappa}s_{i}\langle\hat{\bF}^{(0)}(\bx^{i}),\bv^{i}\rangle\bigg)
		& \leq{\rm Var}\bigg(\sum_{i\leq\kappa}s_{i}\langle\hat{\bF}^{(0)}(\bx^{i})-\hat{\bF}^{(0)}(\bx),\bv^{i}\rangle\bigg)\\
		& =\sum_{i\leq\kappa}\sum_{i'\leq\kappa}s_{i}s_{i'}\<\bv^{i},\bv^{i'}\>\left(\xi(\langle\bx^{i},\bx^{i'}\rangle)-\xi(\langle\bx^{i},\bx\rangle)-\xi(\langle\bx^{i'},\bx\rangle)+\xi(\|\bx\|^{2})\right)\\
		& \leq4\xi(1)\sum_{i\leq\kappa}s_{i}^{2}+3\eps\xi'(1)\sum_{i\leq\kappa}|s_{i}|\sum_{i'\leq\kappa}|s_{i'}|
		\leq4\xi(1)+3\kappa\eps\xi'(1)\,,
	\end{align*}
	where the second to last inequality follows since, for  $(\bx_i:i\le \kappa)\in T(\bx)$,
	\[
	\big|\xi(\langle\bx^{i},\bx^{i'}\rangle)-\xi(\langle\bx^{i},\bx\rangle)-\xi(\langle\bx^{i'},\bx\rangle)+\xi(\|\bx\|^{2})\big|\leq\begin{cases}
		4\xi(1) & i=i'\\
		3\eps\xi'(1) & i\neq i'
	\end{cases}
	\,.\qedhere
	\]
\end{proof}

Next, we state and prove a simple lemma about concentration of min-max problems.
\begin{lemma}\label{lemma:GenericConc}
Let $\bg\sim\normal(\bzero,\bI_D)$ be a standard Gaussian vector and 
$\bh(\,\cdot\, ;j,k) = (h_i(\,\cdot\, ;j,k):\; i\le m)$,
$\bh(\,\cdot\,;j,k):\reals^D\to \reals^m$ be $L$-Lipschitz for all
$j\le J$, $k\le K$.
Then, for
\begin{align*}
&\psi(\bg):= \min_{j\le J}\max_{k\le K} \frac{1}{m}\sum_{i=1}^mh_{i}(\bg;j,k)^2\,,
\end{align*}
we have 
\begin{align*}
&\prob\Big(\big|\sqrt{\psi(\bg)}-\E \sqrt{\psi(\bg)}\big|\ge t \Big)\le 2\, e^{-mt^2/2L^2}\, .
\end{align*}
\end{lemma}
\begin{proof}
We claim that $\bg \mapsto \sqrt{\psi(\bg)}$ is $L/\sqrt{m}$ Lipschitz, whence the claim 
follows by Gaussian concentration (where, here, by abuse of notation $\bg\in \R^D$ is deterministic). Indeed, 
\begin{align*}
\Big|\sqrt{\psi(\bg_1)}-\sqrt{\psi(\bg_2)}\Big|&\le 
\frac{1}{\sqrt{m}}\max_{j\le J}\max_{k\le K} \Big|\big\|\bh(\bg_1;j,k)\big\|
-\big\|\bh(\bg_2;j,k)\big\|\Big|\\
&\le 
\frac{1}{\sqrt{m}}\max_{j\le J}\max_{k\le K} \big\|\bh(\bg_1;j,k)-
\bh(\bg_2;j,k)\big\|\\
&\le \frac{L}{\sqrt{m}} \|\bg_1-\bg_2\|\, .\qedhere
\end{align*}
\end{proof}

Using the last two lemmas, we first prove the following concentration
result for a single point $\bx$.
\begin{lemma}
	\label{lem:Delta_concentration_one_pt}Let $0\leq q<q'\leq1$, $0\leq\rho\leq\rho'\leq1$,
	$\eps>0$, $\kappa\geq1$ and suppose that $\|\bx\|^{2}=q$. For any
	$\delta>0$,
	\begin{equation}
		\P\left(\Big|\sqrt{\Delta(\bx)}-\E\big[\sqrt{\Delta(\bx)}\,|\,H^{(0)}(\bx)\big]\Big|>\delta\sqrt{d}\right)\leq2\exp\left(-\frac{d\delta^{2}}{4\xi(1)/\kappa+4\eps\xi'(1)}\right).\label{eq:Delta_deviation}
	\end{equation}
\end{lemma}

\begin{proof}
	For $\bt\in\R^{n}$, denote $\P_{\bt}(\cdot):=\P(\,\cdot\,|\,\bF^{(0)}(\bx)=\bt)$
	and let $\E_{\bt}(\cdot)$ and ${\rm Var}_{\bt}(\cdot)$ be the corresponding
	conditional expectation and variance. To prove (\ref{eq:Delta_deviation})
	it will be enough to show that, for arbitrary $\bt$,
	\begin{equation}
\P_{\bt}\left(\Big|\sqrt{\Delta(\bx)}-\E_{\bt}\big[\sqrt{\Delta(\bx)}\big]\Big|>\delta\sqrt{d}\right)
			\le 2\exp\left(-\frac{d\delta^{2}}{4\xi(1)/\kappa+4\eps\xi'(1)}\right).
		\label{eq:Delta_conditional}
	\end{equation}

Write
\begin{align*}
	\Delta(\bx) & =\frac{1}{2\kappa}
	\inf_{(\bx^i)\in T(\bx)}\sup_{(\bv^{i})\in (\S^{n-1})^{\kappa}}\sum_{i=1}^{\kappa}
	\<\bF^{(i)}(\bx^{i}),\bv^{i}\>^2\, .
\end{align*}
By a standard approximation argument, it is also enough to prove \eqref{eq:Delta_conditional}
with
\begin{align*}
	\Delta_0(\bx) & =\frac{1}{2\kappa}
	\inf_{(\bx^i)\in M_1}\sup_{(\bv^{i})\in M_2}\sum_{i=1}^{\kappa}
	\<\bF^{(i)}(\bx^{i}),\bv^{i}\>^2
\end{align*}
instead of $\Delta(\bx)$ for arbitrary finite $M_1\subset T(\bx)$ and  $M_2\subset (\S^{n-1})^{\kappa}$
(see e.g. \cite{AT,zeitouni2012branching} and
note that $\E\big[\sup_{T(\bx)}\big|\sum_{i=1}^{\kappa}\big(H^{(i)}(\bx^{i})-H^{(0)}(\bx)\big)\big|\,|\,\bF^{(0)}(\bx)=\bt\big]<\infty$
by the Borell-TIS inequality and separability). Enumerating $M_1=\{(\bx^i(j))_{i=1}^\kappa\}_{j\le J}$ and $M_2=\{(\bv^i(k))_{i=1}^\kappa \}_{k\le K}$ for $J=|M_1|$ and $K=|M_2|$, we may also write 
\begin{align*}
	\Delta_0(\bx) & =\frac{1}{2\kappa}
	\inf_{j\le J}\sup_{k\le K}\sum_{i=1}^{\kappa}
	\<\bF^{(i)}(\bx^{i}(j)),\bv^{i}(k)\>^2\,.
\end{align*}

Consider the  Gaussian random variables $\left\{ 2^{-\frac12}\<\bF^{(i)}(\bx^{i}(j)),\bv^{i}(k)\> :\,i\leq\kappa,\,j\leq J,\,k\le K\right\} $.
Of course, we can augment the probability space of $\P_{\bt}(\cdot)$
with random variables $\bg=(g_{1},\ldots,g_{JK})$ such that
$g_{i}$ are i.i.d. $\normal(0,1)$ variables and %
\begin{align}
	h_i(\bg;j,k) := \frac{1}{\sqrt2}\<\bF^{(i)}(\bx^{i}(j)),\bv^{i}(k)\>=\alpha_{0,i,j,k}+\sum_{\ell\ge 1}\alpha_{\ell,i,j,k}g_{\ell}\, ,
\end{align}
for an appropriate choice of deterministic coefficients $\alpha_{\ell,i,j,k}\in\R$
(depending on $\bt$).

By Lemma \ref{lemma:GenericConc},  to prove \eqref{eq:Delta_conditional} with $\Delta_0(\bx)$ and complete the proof,
it is sufficient to prove an appropriate bound for the Lipschitz constant of $\bw\mapsto \bh(\bw;j, k)$. 
On the other hand, this is an affine function by construction, so denoting the gradient in $\bw$ by $\nabla_{\bw}$, using Lemma \ref{lem:conditional_cov_bd}, for any $\bw\in\R^{JK}$ we have
\begin{align*}
	\big\|\nabla_{\bw} \bh(\bw;j, k)\big\|^2 &= \sup_{\|\bs\|= 1}
	\Var_{\bt}\{\<\bs,\bh(\bg;j, k)\>\} \le 2\xi(1)+ 2\eps\kappa\xi'(1) \, .\qedhere
\end{align*}
\end{proof}

The following concentration result still deals with a single point
$\bx$, but now for $\Delta(\bx)$ instead of $\sqrt{\Delta(\bx)}$.
\begin{corollary}
	\label{cor:Delta_concentration_one_pt-2}For any $\delta,c>0$, there exist $\eps_{0}$,$\kappa_{0}$ and $d_{0}$
	which depend on $\xi$, $\delta$ and $c$, such that for any $\eps\leq\eps_{0}$,
	$\kappa\geq\kappa_{0}$ and $d\geq d_{0}$ and any $0\leq q<q'\leq1$ and $\bx$ with $\|\bx\|^{2}=q$,
		\begin{align}\label{eq:sqrtDelta_deviation}
		\E\bigg[ \Big( \sqrt{\Delta(\bx)}-\E\Big[\sqrt{\Delta(\bx)}\,\Big|\,H^{(0)}(\bx)\Big] \Big)^2  \,\bigg|\, H^{(0)}(\bx) \bigg] & 	\le 16 (\xi(1)/\kappa+\eps\xi'(1))
	\end{align}
	and
	\begin{equation}
		\P\left(\Big|\Delta(\bx)-\E\big[\Delta(\bx)\,|\,H^{(0)}(\bx)\big]\Big|>\delta d\right)\leq e^{-cd}.\label{eq:Delta_deviation-2}
	\end{equation}
\end{corollary}

\begin{proof}
	Let $\delta,c>0$. Abbreviating $\E_{0}[\,\cdot\,]:=\E\big[\cdot|\,H^{(0)}(\bx)\big]$,
	write
	\[
	\Big|\Delta(\bx)-\E_{0}\big[\Delta(\bx)\big]\Big|=\Big|\sqrt{\Delta(\bx)}-\sqrt{\E_{0}\big[\Delta(\bx)\big]}\Big|\cdot\Big(\sqrt{\Delta(\bx)}+\sqrt{\E_{0}\big[\Delta(\bx)\big]}\Big).
	\]
	For $C=C(\xi,c)>0$  to be determined
	below, define the events
	\begin{align*}
		\calE_{1} & :=\left\{ \sqrt{\Delta(\bx)}+\sqrt{\E_{0}\big[\Delta(\bx)\big]}>\frac{1}{2}\sqrt{d}C\right\} ,\\
		\calE_{2} & :=\left\{ \Big|\sqrt{\E_{0}\big[\Delta(\bx)\big]}-\E_{0}\big[\sqrt{\Delta(\bx)}\big]\Big|>\sqrt{d}\delta/C\right\} ,\\
		\calE_{3} & :=\left\{ \Big|\sqrt{\Delta(\bx)}-\E_{0}\big[\sqrt{\Delta(\bx)}\big]\Big|>\sqrt{d}\delta/C\right\} .
	\end{align*}
	To prove the corollary it will be enough to show that, for $\eps$,
	$\kappa$ and $d$ as in its statement, 
	\[
	\P(\calE_{i})\leq e^{-cd}/3,\quad \forall i=1,2,3\,.
	\]
	
	The fact that $\P(\calE_{3})\leq e^{-cd}/3$ follows from Lemma \ref{lem:Delta_concentration_one_pt}. 
	Denoting $\Lambda=\sqrt{\Delta(\bx)}-\E_{0}\sqrt{\Delta(\bx)}$ and
	$\Lambda_{0}=\E_{0}\sqrt{\Delta(\bx)}$, we have
		\begin{align*}
		\E_{0}\big[\Delta(\bx)\big] & =\E_{0}\big[(\Lambda+\Lambda_{0})^{2}\big]=\E_{0}\Lambda^{2}+2\Lambda_{0}\E_{0}\Lambda+\Lambda_{0}^{2}=\E_{0}\Lambda^{2}+\Lambda_{0}^{2}\,.
	\end{align*}
	Thus, $\sqrt{\E_{0}\Delta(\bx)} \leq\sqrt{\E_{0}\Lambda^{2}}+\Lambda_{0}$ and
	\[
	0\le \sqrt{\E_{0}\Delta(\bx)}-\E_{0} \sqrt{\Delta(\bx)} \leq\sqrt{\E_{0}\Lambda^{2}}\,.
	\]
	Therefore,
	\begin{align*}
		\calE_{2} & \subset \Big\{\sqrt{\E_{0}\Lambda^{2}}>\sqrt{d}\delta/C\Big\}=\Big\{\E_{0}\Lambda^{2}>d\delta^2/C^2\Big\}\,.
	\end{align*}
	
	Note that under $\P_{\bt}(\cdot)$, $\Lambda$ and 
	$\sqrt{\Delta(\bx)}-\E_{\bt}\sqrt{\Delta(\bx)}$ have the same law. 
	Hence, by the tail formula for expectation and 
	Eq.~\eqref{eq:Delta_conditional},
	which is uniform over $\bt$, 
	\begin{align*}
\E_{0}\Lambda^{2} & \leq  \int_0^\infty 4 \exp\left(-\frac{s^{2}}{4\xi(1)/\kappa+4\eps\xi'(1)}\right) s\, \dd s
\le 4 (4\xi(1)/\kappa+4\eps\xi'(1))\, ,
\end{align*}
namely, the bound of \eqref{eq:sqrtDelta_deviation} holds.
	Hence, for $\eps$,
	$\kappa$ and $d$ as in the statement of the corollary, $\calE_2$ is empty
and in particular $\P(\calE_{2})\leq e^{-cd}/3$

	Fix some $(\by^{1},\ldots,\by^{\kappa})$ 
	in $T(\bx,\kappa,0,q,q')$. Since 
	\[
	\Delta(\bx)\leq D(\bx):=\frac{1}{\kappa}\sum_{i=1}^{\kappa}H^{(i)}(\by^{i})=\frac{1}{2\kappa}\sum_{i=1}^{\kappa}\|\bF^{(i)}(\by^{i})\|^{2},
	\]
	to prove that $\P(\calE_{1})\leq e^{-cd}/3$ it suffices to show that
	\begin{align*}
		\P\Big(\sqrt{D(\bx)}>\frac{1}{4}\sqrt{d}C\Big) \vee \P\Big(\sqrt{\E_{0}\big[D(\bx)\big]}>\frac{1}{4}\sqrt{d}C\Big) \leq e^{-cd}/6\,.
	\end{align*}
	Recalling \eqref{eq:FhatF}, these bounds easily follow for sufficiently large $C=C(\xi,c)$ and $\eps$, $\kappa$ and $d$
	as in the statement of the corollary.
\end{proof}

\begin{proposition}
	\label{prop:unif_conc}For any $\delta,c>0$
	there exist $\eps_{0}$, $\kappa_{0}$ and $d_{0}$ which depend on
	$\delta$, $c$ and $\xi$ such that for any 
	$\eps\leq\eps_{0}$, $\kappa_{0}\leq\kappa$ and $d\geq d_{0}$ and any $0\leq q<q'\leq1$ such that if $q>0$ then $q'/q\le4$,
	\begin{equation}
		\P\left(\sup_{\bx:\,\|\bx\|^{2}=q}\Big|\Delta(\bx)-\E\big[\Delta(\bx)\,|\,H^{(0)}(\bx)\big]\Big|>\delta d\right)<e^{-cd}.\label{eq:Delta_bd_uniform}
	\end{equation}
\end{proposition}	

\begin{proof}

If $q=0$, the proposition follows from Corollary \ref{cor:Delta_concentration_one_pt-2}.
Assume henceforth that $q>0$. Let $C=C(\xi,\delta,c)$ be a sufficiently large
constant to be determined below. Define the Lipschitz constants
\begin{align*}
	L_{d} & :=\sup_{\bx\neq\by\in\sqrt{q}\S^{d-1}}\frac{|\Delta(\bx)-\Delta(\by)|}{\|\bx-\by\|}\,,\quad
	\hat{L}_d :=\sup_{\bx\neq\by\in\sqrt{q}\S^{d-1}}\frac{\big|\E\big[\Delta(\bx)\,|\,H^{(0)}(\bx)\big]-\E\big[\Delta(\by)\,|\,H^{(0)}(\by)\big]\big|}{\|\bx-\by\|}\,,
\end{align*}
where $\sqrt{q}\S^{d-1}:=\{\bx:\,\|\bx\|^{2}=q\}$. We claim that
to complete the proof it suffices to show that, for $\eps$, $\kappa$ and $d$ as in the proposition,
\begin{equation}\label{eq:reqbds_Delta}
	\P(L_d\le Cd\delta/4)\geq1-e^{-2dc}\qquad\mbox{and}\qquad \P(\hat L_d\le Cd\delta/4)\geq1-e^{-2dc}\,.
\end{equation}

Indeed, by Corollary \ref{cor:Delta_concentration_one_pt-2}, for
any $c'>0$ and appropriate $\eps_{0},\kappa_{0},d_{0}$ depending
on $\xi$, $\delta$ and $c'$,
\[
\P\left(\Big|\Delta(\bx)-\E\big[\Delta(\bx)\,|\,H^{(0)}(\bx)\big]\Big|\geq\delta d/2\right)\leq\exp\left(-c'd\right)
\]
for any $\bx\in\sqrt{q}\S^{d-1}$ and $\eps\leq\eps_{0}$, $\kappa\geq\kappa_{0}$
and $d\geq d_{0}$. Let $\calN$ be some $\sqrt{q}/C$-net of $\sqrt{q}\S^{d-1}$
with $|\calN|\leq C_{0}^{d}$ for an appropriate $C_{0}\ge1$ depending
on $C$ only. By a union bound, 
\begin{equation}
	\P\left(\max_{\bx\in\calN}\Big|\Delta(\bx)-\E\big[\Delta(\bx)\,|\,H^{(0)}(\bx)\big]\Big|\leq\frac{d\delta}{2}\right)\geq1-\exp\left(-c'd+\log(C_{0})d\right)\,.\label{eq:net}
\end{equation}
Thus, by taking $c'$ large enough, the event in (\ref{eq:net}),
$\{L_d\le Cd\delta/4\}$ and $\{\hat L_d\le Cd\delta/4\}$ all occur simultaneously with probability
at least $1-e^{-dc}$ for $d\geq d_{0}$, by increasing $d_{0}$ if
needed. On the intersection of those events, for any $\bx\in\sqrt{q}\S^{d-1}$,
denoting by $\bar{\bx}\in\calN$ the closest point to $\bx$, 
\begin{align*}
	\Big|\Delta(\bx)-\E\big[\Delta(\bx)\,|\,H^{(0)}(\bx)\big]\Big| & \leq\Big|\Delta(\bar{\bx})-\E\big[\Delta(\bar{\bx})\,|\,H^{(0)}(\bar{\bx})\big]\Big|\\
	& +\Big|\Delta(\bx)-\Delta(\bar{\bx})\Big|+\Big|\E\big[\Delta(\bx)\,|\,H^{(0)}(\bx)\big]-\E\big[\Delta(\bar{\bx})\,|\,H^{(0)}(\bar{\bx})\big]\Big|\leq\delta d\,.
\end{align*}

It remains to prove the bounds of \eqref{eq:reqbds_Delta}. 
If $\|\bx\|=\|\by\|=q$ and $O$ is an orthogonal matrix which maps
$O\bx=\by$ and acts as the identity on the orthogonal space to $\bx$
and $\by$, then for $(\bx^{1},\ldots,\bx^{\kappa})\in T(\bx)$ we
have that $(O\bx^{1},\ldots,O\bx^{\kappa})\in T(\by)$ and $\|O\bx^{i}-\bx^{i}\|\leq \sqrt{q'/q} \|\bx-\by\|\le 2  \|\bx-\by\|$.
Therefore, 
if we denote by $L_{i,d}$ the Lipschitz constant of $H^{(i)}(\bx)$ on $\Ball^{d}(1)$, then 
the  Lipschitz constant of $\Delta(\bx)$ satisfies $L_d\le 2\kappa^{-1}\Sigma_{i\le \kappa}L_{i,d}$ on $\sqrt{q}\S^{d-1}$. Using H\"{o}lder's inequality and Chernoff's bound, 
\begin{equation}\label{eq:Chernoff}
\P(L_d\ge d A)\le e^{-d\theta A}\E e^{\theta L_d}\le e^{-d\theta A} \E \prod_{i=1}^{\kappa} e^{2\theta L_{i,d}/\kappa}\le e^{-d\theta A}  \prod_{i=1}^{\kappa} \E[ e^{2\theta L_{i,d}}]^{1/\kappa} = e^{-d\theta A}   \E[ e^{2\theta L_{1,d}}]\,,
\end{equation} 
for any $\theta>0$.
Thus, the required bound on $\P(L_d\le Cd\delta/4)$ follows from Lemma \ref{lem:basicestimates} by taking sufficiently small $\theta$ so that $\E[ e^{2\theta L_{1,d}}]\le e^{d\bar c(\xi)}$ for an appropriate  $\bar c(\xi)$ and large enough $A=A(\xi,c,\theta)$, and choosing  $C=C(\xi,\delta,c)$ such that
$A\le C\delta/4$.

Recalling \eqref{eq:FhatF}, for $i\ge1$ we write
\begin{equation}
	\begin{aligned}
		\bF^{(i)}(\bx^{i}) = \frac{\xi(\langle\bx,\bx^{i}\rangle)}{\xi(\langle\bx,\bx\rangle)}\bF^{(0)}(\bx) 
		+\Big(\sqrt{\rho} \hat \bF^{(0)}(\bx^{i}) - \frac{\xi(\langle\bx,\bx^{i}\rangle)}{\xi(\langle\bx,\bx\rangle)}\bF^{(0)}(\bx) \Big)
		+\sqrt{\rho'-\rho}\hat{\bF}^{(i)}(\bx^i)\,,
	\end{aligned}\label{eq:condexp1}
\end{equation}
for $\hat{\bF}^{(i)}$ are i.i.d. copies of $\bF$ and note that the three summands are independent processes of $\bx^i$. Define
\begin{equation*}
	\begin{aligned}
		 \tilde{\bF}_{\bx}^{(i)}(\bx^{i}) := \Big(\sqrt{\rho} \hat \bF^{(0)}(\bx^{i}) - \frac{\xi(\langle\bx,\bx^{i}\rangle)}{\xi(\langle\bx,\bx\rangle)}\bF^{(0)}(\bx) \Big)
		+\sqrt{\rho'-\rho}\hat{\bF}^{(i)}(\bx^i)\,.
	\end{aligned}
\end{equation*}

Fixing some arbitrary unit vector $\bu$, for $\bx,\by\in\sqrt{q}\S^{d-1}$ we have that
\begin{align*}
	& \big|\E\big[\Delta(\bx)\,|\,H^{(0)}(\bx)= t\big]-\E\big[\Delta(\by)\,|\,H^{(0)}(\by)=s\big]\big|\\
	& =\big|\E\big[\Delta(\bx)\,|\,\bF^{(0)}(\bx)=\sqrt{2t}\bu\big]-\E\big[\Delta(\bx)\,|\,\bF^{(0)}(\bx)=\sqrt{2s}\bu\big]\big|\\
	& =\frac{1}{2\kappa}\bigg|\E\min_{T(\bx)}\sum_{i=1}^{\kappa}\bigg\Vert \frac{\xi(\langle\bx,\bx^{i}\rangle)}{\xi(\langle\bx,\bx\rangle)}\sqrt{2t}\bu + \tilde{\bF}_{\bx}^{(i)}(\bx^{i})\bigg\Vert ^{2}-\E\min_{T(\bx)}\sum_{i=1}^{\kappa}\bigg\Vert \frac{\xi(\langle\bx,\bx^{i}\rangle)}{\xi(\langle\bx,\bx\rangle)}\sqrt{2s}\bu + \tilde{\bF}_{\bx}^{(i)}(\bx^{i})\bigg\Vert ^{2}\bigg|\\
	& \leq\frac{1}{\kappa}\E\max_{T(\bx)}\sum_{i=1}^{\kappa}\bigg(\bigg[ \frac{\xi(\langle\bx,\bx^{i}\rangle)}{\xi(\langle\bx,\bx\rangle)}|\sqrt{2t}-\sqrt{2s}|\bigg] \cdot\bigg[\big\Vert \tilde{\bF}_{\bx}^{(i)}(\bx^{i})\big\Vert+\sqrt{2(t\vee s)}\frac{\xi(\langle\bx,\bx^{i}\rangle)}{\xi(\langle\bx,\bx\rangle)} \bigg] \bigg)\,,
\end{align*}
where for the inequality we used $\big|\|\bx\|^{2}-\|\by\|^{2}\big|\leq\|\bx-\by\|\|\bx+\by\|$.

Note that by Lemma \ref{lem:basicestimates}, $\kappa^{-1}\E \sum_{i=1}^{\kappa} \max_{\by\in\Ball^{d}(1)} \Vert \tilde{\bF}^{(i)}(\by)\Vert\le \sqrt d c_1(\xi)$, for some constant $c_1(\xi)$. Moreover, for any $(\bx^j)_{j\le \kappa}\in T(\bx)$ and $i\le \kappa$, $\big|\frac{\xi(\langle\bx,\bx^{i}\rangle)}{\xi(\langle\bx,\bx\rangle)}\big|\le (q'/q)^{\deg(\xi)/2}\le 2^{\deg(\xi)}=:c_2(\xi)$.
Hence, on the event that $\max_{\by\in\Ball^{d}(1)}\|{\bF}^{(0)}(\by)\|\le \sqrt d M$, 
\begin{align*}
	& \Big|\E\big[\Delta(\bx)\,|\,H^{(0)}(\bx)\big]-\E\big[\Delta(\by)\,|\,H^{(0)}(\by)\big]\Big|\le \sqrt d c_2(\xi)(c_1(\xi)+Mc_2(\xi)) \Big|\|\bF^{(0)}(\bx)\|-\|\bF^{(0)}(\by)\|\Big|\,.
\end{align*}
The bound $\P(\hat L_d\le Cd\delta/4)\geq1-e^{-2dc}$ follows for $d\ge d_{0}(\xi,c)$ with appropriate $d_{0}(\xi,c)$ and $C=C(\xi,\delta,c)$ sufficiently large by Lemma \ref{lem:basicestimates}.
\end{proof}

\subsection{Proof of Theorem \ref{thm:concentration}}

For any $u\in\TT\setminus(\LL\cup \r)$, we have $[\< u,u\>+(1-\theta)/\tau]/\< u,u\>\le 2$. 
Therefore, \eqref{eq:Delta_bd_uniform-6} follows using \eqref{eq:Deltaequiv} by a union bound over $u\in\TT\setminus(\LL\cup \r)$ from Proposition \ref{prop:unif_conc}. \qed

\section{Energy difference at one level: Proof of Proposition \ref{prop:main_bound}}
\label{sec:ProofMainBound}

We start by stating several auxiliary results, which we prove in Subsections \ref{subsec:pf_Delta2_eps0} to \ref{subsec:pf_quadratic_min} below. In Subsection \ref{subsec:pf_main_bound} we  combine them to prove Proposition \ref{prop:main_bound}. 
Until Subsection \ref{subsec:pf_main_bound} our standing assumption (which will not be repeated) is that $\bx$ is an arbitrary point such that $\|\bx\|^2=q$ and that $0\leq q<q'\leq1$, $0\leq\rho\leq\rho'\leq1$, $\kappa\geq1$, $\eps>0$ and $\uu_0\ge0$ are general parameters such that either $q=0$ or $q'/q\le 4$, and if $q=\xi(0)=0$ then $\uu_0=0$ (since in this case $H^{(0)}(\bx)=0$ a.s.). 

For the convenience of the reader, we recall definitions \eqref{eq:T(x)} and 
\eqref{eq:DeltaDef}:
\begin{align*}
	T(\bx)&:=  T(\bx,\kappa,\eps,q,q'):=\Big\{(\bx^{1},\ldots,\bx^{\kappa})\in(\Ball^{d}(1))^{\kappa}:\,\forall i<j\le\kappa,\\
	& \hphantom{:= T(\bx,\kappa,\eps,q,q'):=\Big\{}\big|\langle\bx^{i},\bx\rangle-q\big|\vee\big|\langle\bx^{i},\bx^{j}\rangle-q\big|\le\eps,\,\langle\bx^{i},\bx^{i}\rangle=q'\Big\}\,,\\
	\Delta(\bx)&:=  \Delta(\bx,\kappa,\eps,q,q',\rho,\rho'):=\frac{1}{\kappa}\min_{T(\bx)}\sum_{i=1}^{\kappa}H^{(i)}(\bx^{i})\,.
\end{align*}
Here, $H^{(i)}(\bx):=\frac{1}{2}\|\bF^{(i)}(\bx)\|^{2}$ for $\bF^{(i)}$ as in \eqref{eq:correlation}, which we may write as 
\begin{equation}\label{eq:FihatFi}
\bF^{(i)}(\bx)=\sqrt{\rho}\hat{\bF}^{(0)}(\bx)+\sqrt{\rho'-\rho}\hat{\bF}^{(i)}(\bx)
\end{equation}
for $i\geq1$ and $\bF^{(0)}(\bx)=\sqrt{\rho}\hat{\bF}^{(0)}(\bx)$,
for i.i.d. copies $\hat{\bF}^{(i)}(\bx)$ of $\bF(\bx)$.

To
make the dependence in $\eps$ explicit, we denote $T(\bx,\eps)=T(\bx,\kappa,\eps,q,q')$
and $\Delta(\bx,\eps)=\Delta(\bx,\kappa,\eps,\rho,\rho',q,q')$. We prove the following lemma in Subsection \ref{subsec:pf_Delta2_eps0}. In the setting of the lemma, if $\xi(0)=0$ and $q=0$, then $\frac{1}{n}H^{(0)}(\bx)=0$ a.s. and we only allow to condition on $\frac{1}{n}H^{(0)}(\bx)=\uu_0$ for $\uu_0=0$. We implicitly assume the same for the rest of the section without repeating this remark.
\begin{lemma}
	\label{lem:Delta2_eps0}There
	exists a constant $C(\xi,\alpha)$ such that
	\begin{equation}
			\E\left[\frac{1}{n}\big|\Delta(\bx,\eps)-\Delta(\bx,0)\big|\,\Big|\,\frac{1}{n}H^{(0)}(\bx)=\uu_0\right]\leq C(\xi,\alpha)(\kappa+1)\sqrt{\eps}\,(\uu_0\vee1)\,.\label{eq:diffeps}
	\end{equation}
\end{lemma}

Note that, by \eqref{eq:FihatFi},
\begin{equation}
	\label{eq:H_diff}
	\begin{aligned}
		H^{(i)}(\bx^{i})-H^{(0)}(\bx)&=\frac{1}{2}\big\Vert \sqrt{\rho}\hat{\bF}^{(0)}(\bx^{i})+\sqrt{\rho'-\rho}\hat{\bF}^{(i)}(\bx^{i})\big\Vert ^{2}-\frac{1}{2}\big\Vert \sqrt{\rho}\hat{\bF}^{(0)}(\bx)\big\Vert ^{2}\\
		&=(\rho'-\rho)\Ia^{(i)}(\bx^{i})+\sqrt{\rho}\sqrt{\rho'-\rho}\Ib^{(i)}(\bx^{i})+\rho\Ic^{(i)}(\bx^{i})\,,
	\end{aligned}
\end{equation}
	where
	\begin{equation*}
		\begin{aligned}\Ia^{(i)}(\bx^{i}) & =\frac{1}{2}\bigl\Vert\hat{\bF}^{(i)}(\bx^{i})\bigl\Vert^{2},\\
			\Ib^{(i)}(\bx^{i}) & =\bigl\langle\hat{\bF}^{(0)}(\bx^{i}),\hat{\bF}^{(i)}(\bx^{i})\bigr\rangle,\\
			\Ic^{(i)}(\bx^{i}) & =\frac{1}{2}\left(\bigl\Vert\hat{\bF}^{(0)}(\bx^{i})\bigl\Vert^{2}-\bigl\Vert\hat{\bF}^{(0)}(\bx)\bigl\Vert^{2}\right).
		\end{aligned}
	\end{equation*}

Denote by $\bD\bF(\bx):=(\partial_{j}F_{i}(\bx))_{j,i}\in\R^{n\times d}$
the differential matrix of $\bF(\bx)$ and define 
\begin{equation}
	\begin{aligned}\Iabar^{(i)}(\bx^{i}) & :=\hat{H}^{(i)}(\bx),\\
		\Ibbar^{(i)}(\bx^{i}) & :=\left\langle \hat{\bF}^{(0)}(\bx),\bD\hat{\bF}^{(i)}(\bx)(\bx^{i}-\bx)\right\rangle +\left\langle \bD\hat{\bF}^{(0)}(\bx)(\bx^{i}-\bx),\hat{\bF}^{(i)}(\bx)\right\rangle ,\\
		\Icbar^{(i)}(\bx^{i}) & :=\frac{1}{2}\left\langle \nabla^{2}\hat{H}^{(0)}(\bx),(\bx^{i}-\bx)^{^{\otimes2}}\right\rangle .
	\end{aligned}
	\label{eq:Ibar}
\end{equation}
We prove the following in Subsection \ref{subsec:pf_I-III}.

\begin{lemma}
	\label{lem:I-III}For some constant $C(\xi,\alpha)>0$, 
	\begin{align}
		\label{eq:I_error}\E\left[\frac{1}{\kappa n}\sup_{T(\bx,0)}\sum_{i=1}^{\kappa}\Big|\Ia^{(i)}(\bx^{i})-\Iabar^{(i)}(\bx^{i})\Big|\,\bigg|\,\frac{1}{n}H^{(0)}(\bx)=\uu_0\right]&\leq C(\xi,\alpha)(q'-q)^{\frac{1}{2}}\,,\\
		\label{eq:II_error}\E\left[\frac{1}{\kappa n}\sup_{T(\bx,0)}\Big|\sum_{i=1}^{\kappa}\Ib^{(i)}(\bx^{i})-\Ibbar^{(i)}(\bx^{i})\Big|\,\bigg|\,\frac{1}{n}H^{(0)}(\bx)=\uu_0\right]&\leq C(\xi,\alpha)(\uu_0\vee1)\Big(\frac{1}{\sqrt{n\kappa\rho}}+\frac{q'-q}{\sqrt{\rho}}\Big)\,,\\
		\label{eq:III_error}\E\left[\frac{1}{\kappa n}\sup_{T(\bx,0)}\Big|\sum_{i=1}^{\kappa}\Ic^{(i)}(\bx^{i})-\Icbar^{(i)}(\bx^{i})\Big|\,\bigg|\,\frac{1}{n}H^{(0)}(\bx)=\uu_0\right]&\leq C(\xi,\alpha)(\uu_0\vee1)\\
		&\times \frac1\rho\Big((q'-q)^{\frac{3}{2}}+\frac{1}{\sqrt{\kappa}}(q'-q)^{\frac{1}{2}}\Big)\,,\nonumber
	\end{align}
where the last two inequalities assume that $\rho>0$.
\end{lemma}

Using \eqref{eq:H_diff} and Lemma \ref{lem:I-III}, we may express $\Delta(\bx,0)$ by the Gaussian variables in \eqref{eq:Ibar}. By computing calculating their joint Gaussian law we conclude the following lemma, which we prove in Subsection \ref{subsec:pf_Delta_Gaussian_vars}.
We write $\bW\sim\GOE(N)$ if $\bW=\bW^{\top}$ and $(W_{ij})_{i\le j\le N}$
are independent with $W_{ii}\sim\normal(0,2)$, $W_{ij}\sim\normal(0,1)$ for
$i<j$. We write $\bZ\sim\GOE(M,N)$ if $\bZ\in\mathbb R^{M\times N}$ has
independent entries $Z_{ij}\sim\normal(0,1)$.

\begin{lemma}
	\label{lem:Delta_Gaussian_vars}Let $\bg\sim \normal(0,\bI_{n})$, $\bw\sim \normal(0,\bI_{d})$,
	$\bZ\sim\GOE(n,d)$ and $\bW\sim\GOE(d)$ be independent
	of each other.  For some constant $C(\xi,\alpha)>0$, 
	\begin{equation}\label{eq:Delta_Gaussian_vars}
		\E\left[\frac{1}{n}\Delta(\bx,0)\,\Big|\,\frac{1}{n}H^{(0)}(\bx)=\uu_0\right]\geq \uu_0+ \Qm - C(\xi,\alpha)(\uu_0\vee1)\Qe
	\end{equation}
	where we define 
	\begin{align*}
		\Qm &= \frac{1}{n}\E\min_{\by\in\S^{d-1}}\bigg\{ \frac{1}{2}(\rho'-\rho)\xi(q)\|\bg\|^{2}+\sqrt{(\rho'-\rho)(q'-q)\xi'(q)}\left(\sqrt{2n\uu_0}\langle\bw,\by\rangle+\sqrt{\rho\xi(q)}\left\langle \bZ^{\top}\bg,\by\right\rangle \right)\\ 
		& \hphantom{\frac{1}{n}\E\min_{\by\in\S^{d-1}}\bigg\{\ }
		+\frac{1}{2}(q'-q) \sqrt{2n\uu_0\rho \xi''(q)} \left\langle \bW,\by\otimes \by\right\rangle 
		+\frac{1}{2}(q'-q)\rho\xi'(q)\left\langle \bZ^{\top}\bZ,\by\otimes \by\right\rangle \bigg\}\,,\\
		\Qe &=(\rho'-\rho)(q'-q)^{\frac{1}{2}}+\sqrt{\rho'-\rho}\Big(\frac{1}{\sqrt{n\kappa}}+q'-q\Big)+(q'-q)^{\frac{3}{2}}+\frac{1}{\sqrt{\kappa}}(q'-q)^{\frac{1}{2}}\,.
	\end{align*}
\end{lemma}

In Section \ref{subsec:pf_quadratic_min} we prove the following lemma for computing the value of $\Qm$.
\begin{lemma}
	\label{lem:quadratic_min}Let $\bg$, $\bw$,
	$\bZ$ and $\bW$ be as in Lemma \ref{lem:Delta_Gaussian_vars} and let $a_0,a_1,b_0,b_1\ge0$. Then, 
	\begin{align}
	\label{eq:quadratic_min}	&\E\min_{\by\in\S^{d-1}}\bigg\{ \frac{b_0}{\sqrt{n}} \langle\bw,\by\rangle+ \frac{b_1}{n} \left\langle \bZ^\top\bg,\by\right\rangle + \frac{a_0}{2\sqrt{n}}\left\langle \bW,\by\otimes \by\right\rangle  		+ \frac{a_1}{2n}\left\langle \bZ^{\top}\bZ,\by\otimes \by\right\rangle \Big\}\\
	\nonumber	&\hphantom{\E\min_{\by\in\S^{d-1}}} \ge \max_{s>0}\Big\{ -\frac12 (a_0^2+b_0^2)s -\frac{1}{2\alpha s} +\frac12 \frac{a_1-b_1^2 s}{1+a_1s}\Big\}  - (a_0\vee a_1\vee b_0\vee b_1) A_d\,,
	\end{align}
for some sequence $A_d=A_d(\alpha)\overset{d\to\infty}{\longrightarrow}0$.
\end{lemma}

We will need the following consequence of Lemmas \ref{lem:basicestimates}  and \ref{lem:basicestimates_condtiional}.
\begin{corollary}\label{cor:conditional_Lip_bds}
	For some constant $C(\xi,\alpha)>0$,
	\begin{align*}
	\E\left[ 	\frac{1}{n}\max_{0\le k\le 3}\sup_{\by\in\Ball^d(\sqrt{q'})}\bigl\Vert\nabla^{k}{\bF}^{(1)}(\by)\bigl\Vert_{\op}\,\Big|\,\frac{1}{n}H^{(0)}(\bx)=\uu_0\right]&\leq  C(\xi,\alpha)\cdot\sqrt{(\uu_0\vee1)/d}\,,\\
	\E\left[ 	\frac{1}{n}\max_{0\le k\le 3}\sup_{\by\in\Ball^d(\sqrt{q'})}\bigl\Vert\nabla^{k}{H}^{(1)}(\by)\bigl\Vert_{\op}\,\Big|\,\frac{1}{n}H^{(0)}(\bx)=\uu_0\right]&\leq C(\xi,\alpha)\cdot (\uu_0\vee1)\,.
	\end{align*}
\end{corollary}
\begin{proof}
	Assume that $q>0$ and recall the definition $D(\xi,r,q):=\xi(q)^{-1}\max_{0\le k\le 3} \xi^{(k)}(rq^{1/2})q^{k/2}$  from Lemma \ref{lem:basicestimates_condtiional}. We have that $D(\xi,\sqrt{q'},q)\le  c_\xi (q'/q)^{\deg(\xi)/2}\le c_\xi'$ for some constants $c_\xi,c_\xi'$ depending on $\xi$ (using our assumption that $q'/q\le 4$).
	Hence, the corollary follows by Lemma \ref{lem:basicestimates_condtiional}.
	
	If $q=0$ and $\xi(0)>0$, the same argument works  with $q'/q$ replaced by $1$. If $q=0$ and $\xi(0)=0$, then $H^{(0)}(\bx)=0$ a.s. and  $\uu_0=0$. In this case, the conditioning can be omitted and we may use Lemma \ref{lem:basicestimates} instead  to conclude the corollary.
\end{proof}

\subsection{\label{subsec:pf_Delta2_eps0}Proof of Lemma \ref{lem:Delta2_eps0}}

		Let $\bx$ with $\|\bx\|^2=q$. Denote $t=(\kappa+1)\sqrt{\eps}(1+r_0)$, where $r_0:=\sqrt{q'/q}$ if $q>0$ and $r_0:=0$ if $q=0$.  Assuming that
	\begin{equation}\label{eq:epsto0mathcing}
			\mbox{for any }(\bx^{1},\ldots,\bx^{\kappa})\in T(\bx,\eps)\mbox{ there exists }(\by^{1},\ldots,\by^{\kappa})\in T(\bx,0)\mbox{ with }\|\by^i-\bx^i\|\le t\,,
	\end{equation} 
	Lemma \ref{lem:Delta2_eps0} follows from Corollary \ref{cor:conditional_Lip_bds} since
	\begin{align*}
		\frac{1}{n}|\Delta(\bx,\eps)-\Delta(\bx,0)|
		& \leq t \frac{1}{\kappa n}\sum_{i=1}^{\kappa} \sup_{\by\in\Ball^d(\sqrt{q'})}\bigl\Vert\nabla{H}^{(i)}(\by)\bigl\Vert_{\op}\,.
	\end{align*}
	In the rest of the proof, we establish \eqref{eq:epsto0mathcing}.

	Define $\bx^0=\bx$. By a change of basis, WLOG we may assume that all but the first $\kappa+1$ entries of $\bx^i$ are zero, for $0\le i\le \kappa$. Hence, we may, and will, assume that $d=\kappa+1$. 
	Define the matrices $A=(\<\bx^i,\bx^j\>)_{i,j=0}^{\kappa}$ and  $Q=(q_{ij})_{i,j=0}^{\kappa}$ where $q_{ij}=q+\delta_{ij}(1-\delta_{i0})(q'-q)$.
	We note that $(A-Q)_{ii}=0$ and $|(A-Q)_{ij}|\le \eps$ for $i\neq j$. Denote by $A^{1/2}$ and $Q^{1/2}$ the unique positive semidefinite square root matrices. Then,
	\[
		\begin{aligned}
		\|A^{1/2}-Q^{1/2}\|_F
		&\le \sqrt{\kappa+1}\|A^{1/2}-Q^{1/2}\|_2
		\le \sqrt{\kappa+1}\|A-Q\|_2^{1/2}\\
		&\le \sqrt{\kappa+1}\|A-Q\|_F^{1/2}
		\le \kappa^{1/4}(\kappa+1)^{3/4}\sqrt{\eps}
		\le(\kappa+1)\sqrt{\eps},
		\end{aligned}
	\]
	see \cite[Sect. X.1]{bhatia} for the second inequality.
	
	Let $X$ be the matrix whose columns are $(\bx^i)_{i=0}^{\kappa}$ so that $X^\top X=A$.
		There exists an orthogonal $O$ such that $X=OA^{1/2}$. Denote the columns of $OQ^{1/2}$ by $(\bar\by^i)$ (also indexed by $0\le i\le\kappa$). We have that $(\bar\by^i)$ have Gram matrix $(OQ^{1/2})^\top OQ^{1/2}=Q$ and
		$\|\bx^i-\bar\by^i\|\le\|X-OQ^{1/2}\|_F\le(\kappa+1)\sqrt{\eps}$.
		If $q=0$, \eqref{eq:epsto0mathcing} follows with $\by^i=\bar\by^i$. Otherwise, choose an orthogonal matrix $O'$ mapping $\bar\by^0$ to $\bx^0$ and acting as the identity outside their span. Then
		$\|O'\bar\by^i-\bar\by^i\|\le \sqrt{q'/q}\,\|\bx^0-\bar\by^0\|$.
		Thus \eqref{eq:epsto0mathcing} follows with $\by^i=O'\bar\by^i$.\qed

\subsection{\label{subsec:pf_I-III}Proof of Lemma \ref{lem:I-III}}
\subsubsection{Proof of \eqref{eq:I_error}}
	Recall that $\Ia^{(i)}(\bx^{i})-\Iabar^{(i)}(\bx^{i})=\hat{H}^{(i)}(\bx^{i})-\hat{H}^{(i)}(\bx)$.
	Since $\|\bx^{i}-\bx\|=\sqrt{q'-q}$ on $T(\bx,0)$,
	\[
	\frac{1}{\kappa n}\sup_{T(\bx,0)}\sum_{i=1}^{\kappa}\Big|\hat{H}^{(i)}(\bx^{i})-\hat{H}^{(i)}(\bx)\Big|\leq\frac{\sqrt{q'-q}}{\kappa n}\sum_{i=1}^{\kappa}\sup_{\bx\neq\by\in\Ball(1)}\frac{\Big|\hat{H}^{(i)}(\by)-\hat{H}^{(i)}(\bx)\Big|}{\|\by-\bx\|}\,.
	\]
	The proof is completed by Lemma \ref{lem:basicestimates}, since the conditioning in \eqref{eq:I_error} can be dropped due to independence of $H^{(0)}(\cdot)$ and $\hat H^{(i)}(\cdot)$.\qed

\subsubsection{Proof of \eqref{eq:II_error}}

	By definition, $\Ib^{(i)}(\bx^{i})-\Ibbar^{(i)}(\bx^{i})$ is equal to 
	\begin{align*}
		\bigl\langle\hat{\bF}^{(0)}(\bx^{i}),\hat{\bF}^{(i)}(\bx^{i})\bigr\rangle-\left\langle \hat{\bF}^{(0)}(\bx),\bD\hat{\bF}^{(i)}(\bx)(\bx^{i}-\bx)\right\rangle -\left\langle \bD\hat{\bF}^{(0)}(\bx)(\bx^{i}-\bx),\hat{\bF}^{(i)}(\bx)\right\rangle \,.
	\end{align*}
		By some algebra,
	\begin{align*}
		\Ib^{(i)}(\bx^{i})-\Ibbar^{(i)}(\bx^{i}) & =\left\langle \hat{\bF}^{(0)}(\bx),\hat{\bF}^{(i)}(\bx)\right\rangle +\left\langle \bD\hat{\bF}^{(0)}(\bx)(\bx^{i}-\bx),\bD\hat{\bF}^{(i)}(\bx)(\bx^{i}-\bx)\right\rangle \\
		& +\left\langle \hat{\bF}^{(0)}(\bx^{i}),B_{i}\right\rangle +\left\langle A_{i},\hat{\bF}^{(i)}(\bx^{i})\right\rangle -\left\langle A_{i},B_{i}\right\rangle \\
		& =:R_{0}^{(i)}(\bx)+R_{1}^{(i)}(\bx^{i})+R_{2}^{(i)}(\bx^{i})+R_{3}^{(i)}(\bx^{i})-R_{4}^{(i)}(\bx^{i})\,,
	\end{align*}
	where we define 
	\begin{align*}
		A_{i} & :=\hat{\bF}^{(0)}(\bx^{i})-\left[\hat{\bF}^{(0)}(\bx)+\bD\hat{\bF}^{(0)}(\bx)(\bx^{i}-\bx)\right]\,,\\
		B_{i} & :=\hat{\bF}^{(i)}(\bx^{i})-\left[\hat{\bF}^{(i)}(\bx)+\bD\hat{\bF}^{(i)}(\bx)(\bx^{i}-\bx)\right]\,.
	\end{align*}

	For $R_{0}^{(i)}(\bx)$, which does not depend on $\bx^{i}$, for $W\sim\normal (0,1)$,
	\[
	\E\left[\frac{1}{\kappa n}\sup_{T(\bx,0)}\Big|\sum_{i=1}^{\kappa}R_{0}^{(i)}(\bx)\Big|\,\Big|\,\frac{1}{n}H^{(0)}(\bx)=\uu_0 \right]=\frac{1}{\kappa n}\E\left[\sqrt{\frac{2\kappa n\uu_0 \xi(q)}{\rho}}|W|\right]\,.
	\]

	Denote$\nabla^{k}\bF(\bx):=((\partial_{j_{1}}\cdots\partial_{j_{k}}F_{i}(\bx))_{i\leq n,j_{1},\ldots,j_{k}\leq d}\in\R^{n}\otimes(\R^{d})^{\otimes k}$
	and define its operator norm by
	\begin{equation}
		\|\nabla^{k}\bF(\bx)\|_{\op}=\max_{\|\bu\|=\|\bv\|=1}\langle\nabla^{k}\bF(\bx),\bv\otimes\bu^{\otimes k}\rangle=\max_{\|\bu\|=\|\bv\|=1}\sum v_{\ell}\partial_{\bu}^{k}F_{\ell}(\bx)\,.\label{eq:tensor_op_norm}
	\end{equation}
	Note that by a Taylor expansion,
	\begin{align*}
		\|A_{i}\| & \leq\frac{1}{2}\max_{\|\bv\|=1}\langle\nabla^{2}\hat{\bF}^{(0)}(\by),\bv\otimes(\bx^{i}-\bx)^{\otimes2}\rangle\,,\\
		\|B_{i}\| & \leq\frac{1}{2}\max_{\|\bv\|=1}\langle\nabla^{2}\hat{\bF}^{(i)}(\bz),\bv\otimes(\bx^{i}-\bx)^{\otimes2}\rangle\,,
	\end{align*}
	for some $\by,\bz\in\R^{d}$ in the interval connecting $\bx^{i}$
	and $\bx$.
	
	Hence, for $1\leq j\leq4$, 
	\[
	|R_{j}^{(i)}(\bx^{i})|\leq\|\bx^{i}-\bx\|^{2}M^{(0)}M^{(i)}=(q'-q)M^{(0)}M^{(i)},
	\]
	where 
	\[
	M^{(i)}:=\max_{\by\in\Ball^{d}(\sqrt{q'})}\left\{ \|\hat{\bF}^{(i)}(\by)\|\vee\|\bD\hat{\bF}^{(i)}(\by)\|_{\op}\vee\|\nabla^{2}\hat{\bF}^{(i)}(\by)\|_{\op}\right\} .
	\]
	
	Note that $M^{(0)}$ and $M^{(i)}$ are independent. By Corollary \ref{cor:conditional_Lip_bds} (applied with $\rho=\rho'$),  
	$$\E\Big[  M^{(0)}\,\Big|\,\frac1nH^{(0)}(\bx)=\uu_0 \Big]\leq C(\xi,\alpha)\cdot\sqrt{d(\uu_0\vee1)/\rho}\,.$$
	For $i\ge1$, $M^{(i)}$ is independent of $H^{(0)}$ and by Lemma \ref{lem:basicestimates}, $\E M^{(i)}\leq C(\xi,\alpha)\cdot\sqrt{d}$.
	Thus,
	\[
	\sup_{\kappa\geq1}\frac{1}{\kappa n}\E\left[\sup_{T(\bx,0)}\Big|\sum_{i=1}^{\kappa}R_{j}^{(i)}(\bx^{i})\Big|\,\bigg|\,\frac{1}{n}H^{(0)}(\bx)=\uu_0\right]\leq(q'-q)C(\xi,\alpha)(\uu_0\vee1)/\sqrt{\rho}.\tag*{\qedsymbol}
	\]

\subsubsection{Proof of \eqref{eq:III_error}}

	Recall that $\Ic^{(i)}(\bx^{i})=\hat{H}^{(0)}(\bx^{i})-\hat{H}^{(0)}(\bx)$. By a Taylor expansion,
	\[
	\Ic^{(i)}(\bx^{i})-\Icbar^{(i)}(\bx^{i})=\left\langle \nabla\hat{H}^{(0)}(\bx),\bx^{i}-\bx\right\rangle +\frac{1}{6}\left\langle \nabla^{3}\hat{H}^{(0)}(\by),(\bx^{i}-\bx)^{\otimes3}\right\rangle 
	\]
	for some $\by\in\R^{d}$ that belongs to the line connecting $\bx$
	and $\bx^{i}$.	
	By Cauchy-Schwarz,
	\[
	\frac{1}{\kappa}\sup_{T(\bx,0)}\Big|\sum_{i=1}^{\kappa}\left\langle \nabla\hat{H}^{(0)}(\bx),\bx^{i}-\bx\right\rangle \Big|\leq\bigl\Vert\nabla\hat{H}^{(0)}(\bx)\bigl\Vert\sup_{T(\bx,0)}\Bigl\Vert\frac{1}{\kappa}\sum_{i=1}^{\kappa}(\bx^{i}-\bx)\Bigr\Vert\leq\sqrt{\frac{q'-q}{\kappa}}\bigl\Vert\nabla\hat{H}^{(0)}(\bx)\bigl\Vert\,,
	\]
	where we used the orthogonality of $(\bx^{i}-\bx)_{i\leq\kappa}$
	and that $\|\bx^{i}-\bx\|=\sqrt{q'-q}$.
	
	By Corollary \ref{cor:conditional_Lip_bds} (with $\rho=\rho'$)
	\[
	\E\left[\bigl\Vert\nabla\hat{H}^{(0)}(\bx)\bigl\Vert\vee\bigl\Vert\nabla^{3}\hat{H}^{(0)}(\bx)\bigl\Vert_{\op}\,\Big|\,\Big|\,\frac{1}{n}H^{(0)}(\bx)=\uu_0\right]\leq dC(\xi)(\uu_0\vee1)/\rho
	\]
	and \eqref{eq:III_error} follows.\qed

\subsection{\label{subsec:pf_Delta_Gaussian_vars}Proof of Lemma \ref{lem:Delta_Gaussian_vars}}

	From the definition of $\Delta(\bx,0)$, \eqref{eq:H_diff} and Lemma \ref{lem:I-III}, we have that
	\begin{align*}
		\E\left[\frac{1}{n}\Delta(\bx,0)\,\Big|\,\frac{1}{n}H^{(0)}(\bx)=\uu_0\right] = \uu_0 +\Qm_\kappa\,,\qquad |\Qm_\kappa-\bar\Qm_\kappa|\le C(\xi,\alpha)(\uu_0\vee1)\Qe\,,
	\end{align*}
	where we define 
	\begin{align*}
		\Qm_\kappa :=\E\left[\frac{1}{\kappa n}\min_{T(\bx,0)}\sum_{i=1}^{\kappa}[(\rho'-\rho)\Ia^{(i)}(\bx^{i})+\sqrt{\rho}\sqrt{\rho'-\rho}\Ib^{(i)}(\bx^{i})+\rho\Ic^{(i)}(\bx^{i})]\,\Big|\,\frac{1}{n}H^{(0)}(\bx)=\uu_0\right]
	\end{align*}
	and similarly define $\bar\Qm_\kappa$ with $\Iabar^{(i)},\Ibbar^{(i)},\Icbar^{(i)}$ instead of $\Ia^{(i)},\Ib^{(i)},\Ic^{(i)}$. From independence of the processes $\hat\bF^{(i)}$,  $\bar\Qm_\kappa\ge \bar\Qm_1$ for any $\kappa\ge1$ (where we note that $T(\bx,0)=T(\bx,\kappa,0,q,q')$ implicitly depends on $\kappa$). Hence,
	\begin{align*}
		\E\left[\frac{1}{n}\Delta(\bx,0)\,\Big|\,\frac{1}{n}H^{(0)}(\bx)=\uu_0\right] \ge \uu_0 +\bar\Qm_1 - C(\xi,\alpha)(\uu_0\vee1)\Qe\,.
	\end{align*}
	
	Assume that $q>0$. From now on we work with $\kappa=1$, for which 
	\begin{equation}\label{eq:Tkappa=1}
		T(\bx,0)=T(\bx,1,0,q,q') = \{ \by: \<\by-\bx,\bx\>=0, \|\by-\bx\|^2=q'-q\}\,.
	\end{equation}
	Let $\bU_{\bx}\in\R^{d\times d-1}$ be an arbitrary matrix whose columns
	form an orthonormal  basis of the orthogonal space to $\bx$. We recall that by
	Lemma K.1 of \cite{montanari2023solving}, we may augment our probability space with variables $\bg_{i}\sim\normal(0,\bI_{n})$, $\bW_{i}\sim\GOE(d-1)$ and $\bZ_{i}\sim\GOE(n,d-1)$
	which are independent, between each other and for different $i$, such that for $i=0,1$,
	\begin{equation*}
		\Big(\hat{\bF}^{(i)}(\bx),\bD\hat{\bF}^{(i)}(\bx)\bU_{\bx},\bU_{\bx}^{\top}\nabla^2\hat{H}^{(i)}(\bx)\bU_{\bx}\Big) = \Big(\sqrt{\xi(q)}\bg_{i},\sqrt{\xi'(q)}\bZ_{i},\sqrt{\xi(q)\xi''(q)}\|\bg_{i}\|\bW_{i}+\xi'(q)\bZ_{i}^{\top}\bZ_{i}\Big).
	\end{equation*}
	Since
	$\hat{\bF}^{(0)}(\cdot)\overset{d}{=}\bF(\cdot)/\sqrt{\rho}$ and
	$\hat{\bF}^{(1)}(\cdot)$
		are independent, the effect of the conditioning $\frac{1}{n}H^{(0)}(\bx)=\uu_0$ on the variables above simply amounts to replacing $\bg_0$ by $\sqrt{2n\uu_0/(\rho\xi(q))}\,\bg_0/\|\bg_0\|$, (i.e., determining the norm of $\bg_0$ through the relation $\|\hat{\bF}^{(0)}(\bx)\|=\sqrt{\xi(q)}\|\bg_0\|=\sqrt{2n\uu_0/\rho}$, leaving its direction uniformly distributed) which affects one of the terms in $\Ibbar^{(1)}(\bx^{1})$ without changing the law of the other two terms $\Iabar^{(1)}(\bx^{1})$ and $\Icbar^{(1)}(\bx^{1})$ in the definition of $\bar\Qm_1$. Namely, conditionally on $\frac{1}{n}H^{(0)}(\bx)=\uu_0$,
	\[
	\begin{aligned}\Iabar^{(1)}(\bx^{1}) & =\frac{1}{2}\xi(q)\|\bg_{1}\|^{2}\,,\\
		\Ibbar^{(1)}(\bx^{1}) & =\left\langle \sqrt{\frac{2n\uu_0}{\rho}} \frac{\bg_{0}}{\|\bg_{0}\|},\sqrt{\xi'(q)}\bZ_{1}\bU_{\bx}^{\top}(\bx^{1}-\bx)\right\rangle +\left\langle \sqrt{\xi'(q)}\bZ_{0}\bU_{\bx}^{\top}(\bx^{1}-\bx),\sqrt{\xi(q)}\bg_{1}\right\rangle \\
		& =\sqrt{\frac{2n\uu_0}{\rho}\xi'(q)}\bw_0^{\top}\bU_{\bx}^{\top}(\bx^{1}-\bx)+\left\langle \sqrt{\xi'(q)}\bZ_{0}\bU_{\bx}^{\top}(\bx^{1}-\bx),\sqrt{\xi(q)}\bg_{1}\right\rangle \,,\\
		\Icbar^{(1)}(\bx^{1}) & =\frac{1}{2}\left\langle \sqrt{\frac{2n\uu_0}{\rho}\xi''(q)}\bU_{\bx}\bW_{0}\bU_{\bx}^{\top}+\xi'(q)\bU_{\bx}\bZ_{0}^{\top}\bZ_{0}\bU_{\bx}^{\top},(\bx^{1}-\bx)^{^{\otimes2}}\right\rangle \,,
	\end{aligned}
	\]
	where we assume that $\bx^{1}\in T(\bx,0)$ and use the fact that
	$\bx^{1}-\bx\perp \bx$ and where $\bg_i$, $\bZ_i$, $\bW_i$ are independent and $\bw_0:=\|\bg_0\|^{-1}\bg_0^\top\bZ_1\sim \normal(0,\bI_{d-1})$. Let 
	$\bw\sim\normal(0,\bI_{d})$, $\bW\sim\GOE(d)$ and $\bZ\sim\GOE(n,d)$ be the variables obtained from extending $\bw_0$, $\bW_0$ and $\bZ_0$ by adding independent Gaussian variables at the additional entries. Let $\hat\bU_{\bx}\in\R^{d\times d}$ be the matrix obtained from $\bU_{\bx}$ by adding $\bx/\|\bx\|$ as the last column. Then, since $\bx^{1}-\bx\perp \bx$, replacing $\bw_0^{\top}\bU_{\bx}^{\top}$, $\bZ_{0}\bU_{\bx}^{\top}$, $\bU_{\bx}\bZ_{0}^{\top}\bZ_{0}\bU_{\bx}^{\top}$ and $\bU_{\bx}\bW_{0}\bU_{\bx}^{\top}$ above by $\bw^{\top}\hat\bU_{\bx}^{\top}$, $\bZ\hat\bU_{\bx}^{\top}$, $\hat\bU_{\bx}\bZ^{\top}\bZ\hat\bU_{\bx}^{\top}$ and $\hat\bU_{\bx}\bW\hat\bU_{\bx}^{\top}$ has no effect on the value. 
	
	Combining the above with a change of variables $\by=\hat\bU_{\bx}^{\top}(\bx^1-\bx)/\sqrt{q'-q}$ we obtain that
	\begin{align*}
		\bar\Qm_1 & = \frac{1}{n}\E\min_{\substack{\by\in\S^{d-1}:\\ \<\by,\be_d\>=0}}\bigg\{ \frac{1}{2}(\rho'-\rho)\xi(q)\|\bg\|^{2}+\sqrt{(\rho'-\rho)(q'-q)\xi'(q)}\left(\sqrt{2n\uu_0}\langle\bw,\by\rangle+\sqrt{\rho\xi(q)}\left\langle \bg,\bZ\by\right\rangle \right)\\
		& \hphantom{\frac{1}{n}\E\min_{\by\in\S^{d-1}}\bigg\{\ }
		+\frac{1}{2}(q'-q) \sqrt{2n\uu_0\rho \xi''(q)} \left\langle \bW,\by\otimes \by\right\rangle 
		+\frac{1}{2}(q'-q)\rho\xi'(q)\left\langle \bZ^{\top}\bZ,\by\otimes \by\right\rangle \bigg\}\,,
	\end{align*}
	where $\be_d=(0,\ldots,0,1)\in \R^d$. Removing the constraint that $\<\by,\be_d\>=0$ only decreases the minimum, which proves \eqref{eq:Delta_Gaussian_vars}.
	
	Finally, if $q=0$ then $\bx=\bzero$ and \eqref{eq:Delta_Gaussian_vars} follows by a similar proof, without the need to restrict everywhere to the orthogonal space to $\bx$. In this case the final step of removing the constraint $\<\by,\be_d\>=0$ is also not required as $T(\bx,0)=\sqrt{q'}\S^{d-1}$ for $\kappa=1$.\qed

	\subsection{Proof of Lemma \ref{lem:quadratic_min}}
	\label{subsec:pf_quadratic_min}

    Throughout the proof we use the notations $a:=(a_0,a_1)$, $b:=(b_0,b_1)$, and 
    \begin{align}\label{eq:A_b_def}
	\bA:= \frac{a_0}{\sqrt{n}}\bW+\frac{a_1}{n}\bZ^{\sT}\bZ\,,\;\;\;\; 
	\bb := \frac{b_0}{\sqrt{n}}\bw+\frac{b_1}{n}\bZ^{\sT}\bg\, .
	\end{align}
	With these notation, the value in the expectation in \eqref{eq:quadratic_min} is given by
\begin{align}\label{eq:Rn_Ab_def}
	R_n(a,b) &:= \min_{\by\in\S^{d-1}} \Big\{\<\bb,\by\>+\frac{1}{2}\<\by,\bA\by\>\Big\}\, .
\end{align}
By convex duality \cite[Appendix B]{boyd2004convex}, we obtain
	\begin{align}
	R_n(a,b) &= \sup_{\lambda>-\lambda_{\min}(\bA)}
	L_n(\lambda;a,b)\, ,\\
	L_n(\lambda;a,b)&:= -\frac{1}{2} \lambda -\frac{1}{2}\<\bb,(\bA+\lambda\id)^{-1}\bb\>\, .
	\end{align}
We further define
\begin{align}
	L(z;a,b) & := -\frac{1}{2}(a_0^2+b_0^2) z-\frac{1}{2\alpha z}+\frac{1}{2}
\frac{a_1-b_1^2 z}{1+a_1 z}\, .
\end{align}
	Note that $R_{n}(s\,a,s\,b)=s\,R_n(a,b)$ for any $s>0$,
	by the definition of Eqs.~\eqref{eq:A_b_def},  \eqref{eq:Rn_Ab_def}.
	Similarly, $L(z/s;s\,a,s\,b) = sL(z;a,b)$ by the last display.
	Hence, we can restrict ourselves to the case $\|a\|,\|b\|\le 1$,
	and it is sufficient to prove that
	\begin{align}
		\lim_{d\to\infty}
	\sup_{\|a\|,\|b\|\le 1}\sup_{K^{-1}<n/d<K}
	\Big|\E\max_{\lambda>-\lambda_{\min}(\bA)}
	L_n(\lambda;a,b)- \sup_{z>0}L(z;a,b)\Big| =0\, .
	\end{align}
	Throughout the proof, we use $C, c,\dots$ to denote constants that might depend on $K$.
	As a preliminary remark, note that
	\begin{align}
		\sup_{\|a\|,\|b\|\le 1}\big|R_n(a,b)\big|
		&\leq \sup_{\|b\|\leq1}\|\bb\|
		+\frac12\sup_{\|a\|\leq1}\|\bA\|_{\op}\nonumber\\
		&\leq \frac{\|\bw\|}{\sqrt n}
		+\frac{\|\bZ\|_{\op}\|\bg\|}{n}
		+\frac{\|\bW\|_{\op}}{2\sqrt n}
		+\frac{\|\bZ\|_{\op}^{2}}{2n}
		=:B_n .\label{eq:Delta_uniform_integrability}
	\end{align}
	Standard Gaussian norm estimates imply that, for every fixed $p<\infty$,
	\[
		\sup_{K^{-1}<n/d<K}\E B_n^p<C_{p,K}.
	\]
	Thus the random variables on the left-hand side of
	\eqref{eq:Delta_uniform_integrability} are uniformly integrable. Moreover,
	uniformly over $\|a\|,\|b\|\leq1$ and $K^{-1}<\alpha<K$,
	\[
		-\frac{3+K}{2}\leq L(1;a,b)
		\leq\sup_{z>0}L(z;a,b)\leq\frac12.
	\]
	Consequently, convergence in probability, uniform in the parameters,
	implies convergence in $L^1$. It  is therefore sufficient to prove that, for any $t>0$,
	\begin{align}
	\lim_{d\to\infty}\sup_{\|a\|,\|b\|\le 1}\sup_{K^{-1}<n/d<K}\prob\left(
	\Big|\max_{\lambda>-\lambda_{\min}(\bA)}
	L_n(\lambda;a,b)- \sup_{z>0}L(z;a,b)\Big|\ge t\right) =0\, .\label{eq:MainBoundQuadOpt}
	\end{align}
	The rest of the proof is devoted to establishing this.
	Define 
	\begin{align}
	S_0(\lambda)&:= \frac{1}{n}\Tr\Big(\Big(\frac{a_0}{\sqrt{n}}\bW+\frac{a_1}{n}\bZ^{\sT}\bZ+\lambda\id\Big)^{-1}\Big)\, ,\\
	S_1(\lambda)&:= \frac{1}{n^2}\Tr\Big(\Big(\frac{a_0}{\sqrt{n}}\bW+\frac{a_1}{n}\bZ^{\sT}\bZ+\lambda\id\Big)^{-1}\bZ^{\sT}
	\bZ\Big)\, ,\\
	L_{*,n}(\lambda;a,b)&:=-\frac{1}{2} \lambda-\frac{1}{2}b_0^2S_0(\lambda)-\frac{1}{2}b_1^2S_1(\lambda)\, .
	\end{align}
	Note that 
	\begin{align}
	L_n(\lambda;a,b) &= -\frac{1}{2}\lambda -\frac{1}{2}\<\tbg,\bM\tbg\>\, ,\\
	\tbg &:= \left(\begin{matrix} \bw \\ \bg \end{matrix}\right)\sim\normal(0,\bI_{d+n})\, ,\\
	\bM&:= \frac{1}n\left(\begin{matrix} b_0\bI_d \\ b_1\bZ/\sqrt{n} \end{matrix}\right)
	\Big(\bA+\lambda\bI_d\Big)^{-1}
	\left(\begin{matrix} b_0\bI_d & b_1\bZ^{\sT}/\sqrt{n} \end{matrix}\right)\, .
	\end{align}
	By Hanson-Wright inequality, there exists an absolute constand $C$ 
	such that for any $t>0$,
	\begin{align}
		\P\Big(\big|L_n(\lambda;a,b)-L_{*,n}(\lambda;a,b)\big|\ge 
		\sqrt{n}\|\bM\|_{\op}t
		 \Big|\, \bZ,\bW\Big)\le  C\, e^{-(t^2\wedge t\sqrt{n})/C}\, ,\label{eq:Bound_Delta_Delta0}
		\end{align}
	Hence, on the event  $\lambda>-\lambda_{\min}(\bA) + \eps$
	(which measurable on $\bZ,\bW$),
	\begin{align}
	\P\Big(\big|L_n(\lambda;a,b)-L_{*,n}(\lambda;a,b)\big|\ge \frac{(|b_0|\vee |b_1|\|\bZ\|_{\op}/\sqrt{n})^2}{\eps\sqrt{n}}
	t \Big|\, \bZ,\bW\Big)\le  C\, e^{-(t^2\wedge t\sqrt{n})/C}\, .\label{eq:Bound_Delta_Delta0}
	\end{align}

	Define the analytic function 
	\begin{align} 
	 \ell(z):= a_0^2 z+\frac{1}{\alpha z} -\frac{a_1}{1+a_1 z}\, ,
	 \end{align}
	 and let $z_*  := \argmin_{z>0} \ell(z)$ and $\lambda_* = \ell(z_*)$.
	 For $a_0= 0$, $\alpha\le 1$, we set by convention $z_*=\infty$, $\lambda_*=0$.
	 Note that $z_*$ is the unique 
	 solution\footnote{To show uniqueness,  define $1/z^2 = x$, we can rewrite the equation $\ell'(z) = 0$
	 as $x/\alpha= a^2_0+g(x)$, $g(x) :=a_1x/(a_1+\sqrt{x})^2$, and the claim follows because $g$ 
	 is non-decreasing concave with $g(0) = 0$ and $g(x)\to a_1$ as $x\to\infty$.} of
	 $\ell'(z) = 0$ on $\reals_{>0}$ hence $\ell:\reals_{>0}\to \reals$ is monotone decreasing on
	 $(0,z_*)$, and increasing on $(z_*,\infty)$, with
	  $\lim_{z\to 0}\ell(z)= \lim_{z\to \infty}\ell(z)= \infty$.
	 
	 For $\lambda> \lambda_*$,  denote by $z_0(\lambda)$  the smallest solution of
	 $\ell(z) = \lambda$.
	 We note that $\ell:(0,z_*]\to [\lambda_*,\infty)$ is a bijection with inverse $z_0$.
	 
	 The R-transform of the asymptotic spectral distribution of $\bA$, call it $R_{\bA}$, is the sum 
	 of R-transforms of a Wigner and of a Wishart matrix. Recalling the results for the latter 
	 matrix ensembles \cite{mingo2017free}, we get
	 \begin{align}
	 R_{\bA}(z) = \frac{a_0^2}{\alpha}z +\frac{a_1}{1-a_1z/\alpha}\, .
	 \end{align}
	 The connection between this R-transform and $\ell$ above is given by 
	 $R_{\bA}(z)+z^{-1} = -\ell(-z/\alpha)$. In particular, using the relation between
	 R-transforms and Stieltjis transforms $z_0(\lambda)=\lim_{n\to\infty} S_0(\lambda)$ for
	 any $\Im(\lambda)<0$ (where it is understood that $z_0$ is extended analytically from $(\lambda_*,\infty)$). 
	 
	 By \cite{capitaine2007strong} $\lim_{n\to\infty}\lambda_{\min}(\bA)=-\lambda_*$ almost surely for 
	 $n/d\to\alpha$.
	 Note that  $\lambda_{\min}(\bA)$ is a Lipschitz function of $\bW$, $\bZ$ on the event 
	 $\|\bZ\|_{\op}\le C \sqrt{n}(1\vee 1/\sqrt{\alpha})$, which holds with probability at least $1-\exp(-c_0n)$,
	 with Lipschitz constant at most $C\|a\|(1\vee 1/\sqrt{\alpha})/\sqrt{n}$.
	 Hence, for any  $t, K>0$  there exists $c(t,K)>0$ such that
	 \begin{align}
	 &\sup_{\|a\|\le 1}\sup_{K^{-1}\le n/d\le K}
	 \prob\Big( \big|\lambda_{\min}(\bA)-\Lambda_n(a,\alpha)\big| \ge t \Big)\le 2  e^{-nc(t,K)}\, ,
	 \label{eq:Lambdamin}\\
	 &\Lambda_n(a,\alpha):=\E\lambda_{\min}(\bA)\, .
	 \end{align}

		 Further, the function $(a,\alpha)\mapsto \Lambda_n(a,\alpha)$ is Lipschitz in $a$
		 and $1/2$-H\"older in $\alpha$ over the domain $\|a\|\le 1$,
		 $K^{-1}\le \alpha\le K$, uniformly in $n$.
	 Indeed  
	 \begin{align*}
	 |\Lambda_n(a,\alpha)-\Lambda_n(a',\alpha)|\le \Big\{\frac{1}{\sqrt{n}}\E(\|\bW\|_{\op})+\frac{1}{n}
	 \E(\|\bZ\|_{\op}^2)\Big\}\|a-a'\|_2\le C_1\|a-a'\|_2\, .
	 \end{align*}
	 Further write $\bA=\bA(n,d)$ to emphasize the dimension dependence, and let
	 $d=n/\alpha$, $d'=n/\alpha'$ with $\alpha<\alpha'$, so that $d>d'$.
	 Under the natural coupling, $\bA(n,d')$ is the leading $d'\times d'$ principal
	 minor of $\bA(n,d)$. Interlacing therefore gives
	 \begin{align*}
	 	\Lambda_n(a,\alpha)\le \Lambda_n(a,\alpha')\,.
	 \end{align*}
	 Conversely, let $\bv_{\min}$ be a unit lowest eigenvector of $\bA(n,d)$ and
	 let $\bP_{\le d'}$ be the projector onto the first $d'$ coordinates. By
	 orthogonal invariance, $\bv_{\min}$ is uniformly distributed on the sphere, and
	 hence $\|(\id-\bP_{\le d'})\bv_{\min}\|$ has expectation bounded by
	 $C(K)|\alpha'-\alpha|^{1/2}$. Projecting and normalizing $\bv_{\min}$, and using
	 the uniform moment bound on $\|\bA(n,d)\|_{\op}$, yields
	 \begin{align*}
	 	\Lambda_n(a,\alpha')\le
	 	\Lambda_n(a,\alpha)+C(K)|\alpha'-\alpha|^{1/2}\,.
	 \end{align*}
	
	  Hence $\Lambda_n(a,\alpha)$ converges uniformly to $-\lambda_*(a,\alpha)$ over the same domain
	  $\|a\|\le 1$, $K^{-1}\le \alpha\le K$,
	  and therefore  Eq.~\eqref{eq:Lambdamin} implies (eventually for a smaller $c>0$)
	 \begin{align}
	 \sup_{\|a\|\le 1}\sup_{K^{-1}\le n/d\le K}
	 \prob\Big( \big|\lambda_{\min}(\bA)+\lambda_*(a,n/d)\big| \ge t \Big)\le 2  e^{-nc(t,K)}\, .
	 \label{eq:Lambdamin_bis}
	 \end{align}

	 Since (by the last result) $\lambda\mapsto S_0(\lambda)$, $\lambda\mapsto S_1(\lambda)$,
	are (eventually almost surely) analytic in   $\Re\lambda>\lambda_*+\eps$ for any $\eps>0$,
	convergence of the Stieltjes transform also holds in this domain.
	In other words
	(computation of the asymptotic Stieltjis transform of $S_1(\lambda)$ is standard
	and omitted)
	 \begin{align}
	 &S_0(\lambda) = z_0(\lambda)+E_{0,n}(\lambda) ,\;\;\;\;S_1(\lambda) = \frac{z_0(\lambda)}{1+a_1z_0(\lambda)}+E_{1,n}(\lambda)\, ,\\
	  &\lim_{n\to\infty}E_{i,n}(\lambda) =0\, ,
	 \end{align}
	 where, by using a concentration argument similar to the one for $\lambda_{\min}(\bA)$, we further obtain:
	 \begin{align}
	  \sup_{\|a\|,\|b\|\le 1}\sup_{K^{-1}\le n/d\le K} \sup_{\lambda\ge \lambda_*+\eps}
	  \P\Big(|E_{i,n}(\lambda)|\ge t \Big) & \le 2\, e^{-n c(t,K, \eps)}\, .\label{eq:S0S1_Dev_bound_0}
	 \end{align}
	 
	Putting together Eqs.~\eqref{eq:Bound_Delta_Delta0}, \eqref{eq:Lambdamin_bis}, 
	\eqref{eq:S0S1_Dev_bound_0} we obtain
	 \begin{align}
	  \sup_{\|a\|,\|b\|\le 1}\sup_{K^{-1}\le n/d\le K} \sup_{\lambda\ge \lambda_*+\eps}
	  \P\Big(|L_n(\lambda;a,b)-L(z_0(\lambda);a,b)|\ge t \Big) & \le 2\, e^{-n c(t,K, \eps)}\, .
	  \label{eq:S0S1_Dev_bound}
	 \end{align}

	Also note that, for $\lambda\ge -\lambda_{\min}(\bA)+\eps$,
	 \begin{align}
	 \frac{\partial \phantom{\lambda} }{\partial \lambda}L_n(\lambda;a,b)= 
	-\frac{1}{2}+\frac{1}{2}\<\bb,(\bA+\lambda\id)^{-2}\bb\>\le \frac{1}{2\eps^2} \|\bb\|^2\, .
	\end{align}
		We note that this is bounded on the event $\|\bw\|,\|\bg\|\le C(\sqrt{n}\vee \sqrt{d})$
	and $\|\bZ\|_{\op}\le C(1\vee\alpha^{-1/2})\sqrt{n}$. Therefore, there exist
	constants $C(K,\eps)$ and $c(K)>0$ such that
	 \begin{align}
	  \sup_{\|a\|,\|b\|\le 1}\sup_{K^{-1}\le n/d\le K}
	  \P\Big( \sup_{\lambda\ge \lambda_*+\eps}\big|\partial_{\lambda}L_n(\lambda;a,b)\big|
	  \ge C(K,\eps)\Big)\le 2e^{-c(K)n }\, .
	  \label{eq:LibDelta}
	\end{align}
	For any $M>\eps$, let $S_{\lambda}=\{\lambda_i: i\le N\}$ be $N\le Mn$ points in
	$[\lambda_*+\eps,\lambda_*+M]$ with $|\lambda_{i+1}-\lambda_i|\le 1/n$.
	Taking a union bound and using Eqs.~\eqref{eq:S0S1_Dev_bound} and \eqref{eq:LibDelta}
	we obtain 
	 \begin{align}
	   \sup_{\|a\|,\|b\|\le 1}\sup_{K^{-1}\le n/d\le K}
	  \P\Big(\sup_{\lambda\in[\lambda_*+\eps,\lambda_*+M]}|L_n(\lambda;a,b)-L(z_0(\lambda);a,b)|
	  \ge t \Big) & \le 2\,e^{-n c(t,K, \eps)}\, .\label{eq:UnifDelta_Dev_bound}
	 \end{align}

	We next temporarily assume $b_0\ge \delta>0$ for some positive constant $\delta$
	and will remove the assumption at the end of the proof.
	For convenience, we define $\cS_{\delta}:=\{a,b: \|a\|\le 1, \|b\|\le 1, b_0\ge \delta \}$.
	Let 
	\begin{align}
	z_{M} := \argmax_{z>0}L(z;a,b)\, . 
	\end{align}
	We note that $z_M\in (0,z_*)$, because $L(z;a,b) = -\ell(z)/2 + h(z)$,
	$h(z) := -b_0^2z/2 - b_1^2 z/(1+a_1z)$, $-\ell(z)/2$ is uniquely maximized at $z_*$ and
	$h(z)$ is strictly monotone decreasing. Hence taking $z_-:= z_M/2$, $z_+:= 2z_M$ with
		$L(z_M ;a,b)-L(z_-;a,b)\ge c_*(a,b)$, $L(z_M ;a,b) -L(z_+;a,b)\ge c_*(a,b)$
	for some $c_*(a,b)$ bounded away from $0$ over $(a,b)\in\cS_\delta$
	and $z_{\pm}$ bounded uniformly away from $0$ and $\infty$ over the same domain.

	Let $\lambda_+ := \ell(z_+)$, $\lambda_- := \ell(z_-)$, $\lambda_M := \ell(z_M)$.
	By the above and Eq.~\eqref{eq:UnifDelta_Dev_bound}, there exists
	 $c(\delta,K)>0$ such that
	 \begin{align}
	 \sup_{(a,b)\in \cS_{\delta}}\sup_{K^{-1}\le n/d\le K}
	  \P\Big( \max_{\lambda\in \{\lambda_+,\lambda_-\}}L_n(\lambda;a,b)>L_n(\lambda_M;a,b)
	 \Big) & \le 2\,e^{-n c(\delta,K)}\, .\label{eq:Bracket}
	 \end{align}
	Note that $\lambda\mapsto L_n(\lambda;a,b)$ is strictly concave on $(-\lambda_{\min}(\bA),\infty)$
	with $\lim_{\lambda\to -\lambda_{\min}(\bA)} L_n(\lambda;a,b)=
	\lim_{\lambda\to \infty} L_n(\lambda;a,b)=-\infty$ (the former holds for $b\neq 0$
	with probability $1$). 
	Hence, by Eq.~\eqref{eq:Bracket}, with probability at least $1-Ce^{-cn}$ the unique maximizer of $L_n(\lambda;a,b)$ is
		in $[\lambda_+,\lambda_-]$. In other words
	\begin{align}
	\sup_{(a,b)\in \cS_{\delta}}\sup_{K^{-1}\le n/d\le K}
		  \P\Big( R_n(a,b)\neq \max_{\lambda\in [\lambda_+,\lambda_-]}L_n(\lambda;a,b)\Big) \le
	   2\,e^{-n c(\delta,K)}\, .
	 \end{align}
	Using Eq.~\eqref{eq:UnifDelta_Dev_bound}, we conclude that
	\begin{align}
	\sup_{(a,b)\in \cS_{\delta}}\sup_{K^{-1}\le n/d\le K}
	  \P\Big( \big|R_n(a,b)- \max_{\lambda\in [\lambda_+,\lambda_-]}L(z_0(\lambda);a,b)\big|\ge t\Big) \le
		 2\,e^{-n c(t,\delta,K)}\, .
	 \end{align}
	Finally since  $\max_{\lambda\in [\lambda_+,\lambda_-]}L(z_0(\lambda);a,b) = 
	L(z_M;a,b)= \max_{z>0} L(z;a,b)$ we conclude that
	\begin{align}
	\sup_{(a,b)\in \cS_{\delta}}\sup_{K^{-1}\le n/d\le K}
	  \P\Big( \big|R_n(a,b)- \max_{z>0}L(z;a,b)\big|\ge t\Big) \le
	   2\,e^{-n c(t,\delta,K)}\, .
	 \end{align}

	In order to complete the proof of Eq.~\eqref{eq:MainBoundQuadOpt}, we want to replace 
	 $\cS_{\delta}$ by $ \cS_{0}$ in the last bound. To this end, we claim
		 that $(i)$~$R_n(a,b)$ is uniformly continuous in $b_0$ for $(a,b)\in \cS_{0}$
	  with probability at least $1-2\exp(-cn)$,
	 and $(ii)$~$R(a,b):=\max_{z>0}L(z;a,b)$ is continuous in $b_0$ over the same set.
	  
		For claim $(i)$ consider the event $\|\bw\|\le C\sqrt{n}$, which holds with probability at
		least $1-2\exp(-cn)$ provided $C$ is a sufficiently large constant,
		and define $H(a,b;\bx):= \<\bb,\bx\>+\<\bx,\bA\bx\>/2$. On this event, $b_0\mapsto H(a,(b_0,b_1));\bx)$
		is $C$-Lipschitz for any $\|\bx\|\le 1$, $(a,b)\in\cS_{0}$. Hence,
		$b_0\mapsto R_n(a,b) =\min_{\|\bx\|=1}H(a,(b_0,b_1);\bx)$ is $C$-Lipschitz as well.
	
	Finally, for claim $(ii)$, rewrite  $L_*(a,b) = G_*((a_0^2+b_0^2)/2)$ where
	\begin{align}
	  G_*(u) = \sup_{z>0}\big\{-u z+ g(z)\big\}\, ,\;\;\; g(z) := -\frac{1}{2\alpha z}
	  +\frac{1}{2}\frac{a_1-b_1^2 z}{1+a_1 z}\, .
	\end{align}
	Since $g$ does not depend on $b_0$, it is sufficient to prove that 
	$u\mapsto G_*(u)$ is continuous for $u\ge 0$. 
	For that, we note that, for any $u>0$, $G'_*(u) = -z_*(u)$
	where $z_*(u)\in (0,\infty)$ is the unique maximizer $-u z+ g(z)$ which is bounded for
	$u>0$ because $\lim_{z\to\infty} g(z)$ is finite. it is therefore sufficient to 
	show that $G_*$ is continuous at $u=0$.  Since $G_*$ is non increasing,
	we have $G_*(0)- G_*(u)\ge 0$ for any $u>0$. On the other hand, letting 
		$G_{M}(u) = \sup_{z\in (0,M]}\big\{-u z+ g(z)\big\}\le G_*(u)$, we have that $G_M$ 
	is $M$-Lipschitz, and $\lim_{M\to\infty} G_M(u) = G_*(0)$, whence
	\begin{align}
	G_*(0)-G_*(u) \le G_*(0) - G_M(u) \le G_*(0)-G_M(0) +Mu \, ,
	\end{align}
	whence $\lim\sup_{u\to 0+}[G_*(0)-G_*(u)]\le  G_*(0)-G_M(0)$ and taking $M\to\infty$
	the claim follows. 
	
	\subsection{Proof of Proposition \ref{prop:main_bound}}
	\label{subsec:pf_main_bound}

	By combining Lemmas \ref{lem:Delta2_eps0} and \ref{lem:Delta_Gaussian_vars}, we have that
	\begin{align*}
			\E\left[\frac{1}{n}\Delta(\bx,\eps)\,\Big|\,\frac{1}{n}H^{(0)}(\bx)=\uu_0\right]&\geq \uu_0+ \Qm - C_0(\xi,\alpha,C)\big(\Qe +\kappa\sqrt{\eps} \big)\,,
	\end{align*}
	where $\Qm$ and $\Qe$ are as defined in Lemma \ref{lem:Delta_Gaussian_vars} and $C_0(\xi,\alpha,C)$ is an appropriate constant.
	
	Since $\rho'-\rho=\int_q^{q'} \dot{p}(t)\dd t$, from concavity of $p$ and the bound $\sup_{t\in(q,q')}|\ddot{p}(t)|<C$,
	\[
	0\le (q'-q)\dot{p}(q)-(\rho'-\rho)\le C(q'-q)^2\,.
	\]
	Define $\Qm^+$ and $\Qe^+$ similarly to $\Qm$ and $\Qe$, only with $\rho'-\rho$ replaced everywhere by $(q'-q)\dot{p}(q)$. 
	$\Qm$ is defined as the expectation of the minimum of some function over $\by\in\S^{d-1}$. If $\by_\star\in\S^{d-1}$ is a minimizer of this function, then the linear part of the same function is non-positive at $\by_\star$, since flipping $\by_\star\mapsto-\by_\star$ does not affect the value of the quadratic part. Hence,
	\[
	\Qm^+\le \Qm - \frac{1}{n}\E\bigg\{  \frac{1}{2}(\rho'-\rho)\xi(q)\|\bg\|^{2}\bigg\} +\frac{1}{n}\E \bigg\{ \frac{1}{2}\dot{p}(q)(q'-q)\xi(q)\|\bg\|^{2}\bigg\}\le \Qm  + \frac{C}{2}(q'-q)^2\xi(q)\,.
	\] 
	Moreover,
	\begin{align*}
		|\Qe-\Qe^+| &\le C(q'-q)^{5/2}+(q'-q)\sqrt{C }\Big(\frac{1}{\sqrt{n\kappa}}+q'-q\Big)\,.
	\end{align*}

	Hence, by adjusting the constant $C_0(\xi,\alpha,C)$,
	\begin{align*}
			\E\left[\frac{1}{n}\Delta(\bx,\eps)\,\Big|\,\frac{1}{n}H^{(0)}(\bx)=\uu_0\right]&\geq \uu_0+ \Qm^+ - C_0(\xi,\alpha,C)\Big(\Qe^+ +\kappa\sqrt{\eps} + \frac{q'-q}{\sqrt{n\kappa}}+(q'-q)^2\Big)\,.
	\end{align*}
	It remains to identify the leading term. Set $h=q'-q$ and, in Lemma
	\ref{lem:quadratic_min}, take
	\begin{align*}
		a_0&=h\sqrt{2\uu_0p(q)\xi''(q)},
		& a_1&=hp(q)\xi'(q),\\
		b_0&=h\sqrt{2\uu_0\dot p(q)\xi'(q)},
		& b_1&=h\sqrt{\dot p(q)p(q)\xi(q)\xi'(q)}.
	\end{align*}
	The term in $\Qm^+$ that does not depend on $\by$ has expectation
	$h\dot p(q)\xi(q)/2$. Therefore, Lemma \ref{lem:quadratic_min} gives
	\begin{align*}
	\Qm^+\ge{}&\frac h2\dot p(q)\xi(q)
	+\max_{r>0}\bigg\{-h^2\uu_0\big(\dot p(q)\xi'(q)+p(q)\xi''(q)\big)r
	-\frac{1}{2\alpha r}\\
	&\hspace{42mm}
	+\frac12\frac{hp(q)\xi'(q)-h^2\dot p(q)p(q)\xi(q)\xi'(q)r}
	{1+hp(q)\xi'(q)r}\bigg\}-C_1hA_d,
	\end{align*}
	where the last error bound uses $\uu_0\le C$ and
	$|\dot p(q)|\le C$ (which in turns follows from concavity using
	$\sup_{t\in(q,q')}|\ddot p(t)|<C$). Making the change of variables $r=s/h$ in the
	maximum and recalling \eqref{eq:F}--\eqref{eq:cL}, we obtain
	\begin{equation*}
		\Qm^+\ge h\,\drift(q,\uu_0;p,\dot p)-C_1hA_d.
	\end{equation*}
	Finally, since $0<h\le1$ and $0\le\dot p(q)\le C$,
	\begin{align*}
		\Qe^+ +\frac{h}{\sqrt{n\kappa}}+h^2
		\le C_2\left(h^{3/2}+\frac{h^{1/2}}{\sqrt\kappa}\right).
	\end{align*}
	Substituting the last two displays into the preceding lower bound and
	absorbing $C_1,C_2$ into $C_0(\xi,\alpha,C)$ proves the proposition.\qed

\section{Proof of Theorems \ref{thm:optimalityLip} and \ref{thm:optimalityOC}}
\label{sec:MainProof}

Most of this section is devoted to proving Theorem \ref{thm:optimalityOC}. 
In Subsection \ref{subsec:pf_optimalityLip} we will show that Theorem \ref{thm:optimalityLip}
follows from Theorem \ref{thm:optimalityOC}.
The main auxiliary results that we use in the proof were introduced in Section \ref{sec:Pf_outline_optimalityLip}. The reader is advised to recall the (simplified) outline of the proof in Subsection \ref{subsec:pf of reduction}.

Throughout this section $p=p_d\in\cuC_*$ is fixed (potentially dependent on $d$). We suppose  that $\theta=\theta_{d}\in[0,1)$ and $p(t)=p_{d}(t)$ are as defined
in (\ref{eq:p_linear}). Let $\kappa,\tau\ge1$. 
We consider the correlated systems $\bF^{(u)}$ as in \eqref{eq:bG(u)} for 
$u\in\TT(\kappa,\tau)$. Our main argument will concern their values over the tree of 
points $\bx_u$ for $u\in\TT(\kappa,\tau)$. In the current section, however, 
we view the tree $\TT(\kappa,\tau)$ as a subgraph of the larger 
tree $\TT(\kappa^3,\tau)$ (obtained by augmenting the degree at
each vertex of $\TT(\kappa,\tau)$ and adding the corresponding subtrees).
We assume that the functions $\bF^{(u)}$ and points $\bx_u\in \Ball^d(1)$ are defined through 
Eq.~\eqref{eq:xu} on the vertices of $\TT(\kappa^3,\tau)\setminus \TT(\kappa,\tau)$
as well. 
This will allow us to conclude a stronger estimate on the inner-products when invoking Lemma \ref{lem:tree_overlap}. We emphasize that this larger tree is only used for the definition of $\bx_u$ and nowhere else. We will regularly use the abbreviations $\TT=\TT(\kappa,\tau)$ and $\LL=\LL(\kappa,\tau)$.

One of the arguments we will use requires us to control the approximation 
\eqref{eq:uu_lb_lin_apx}, for which we appeal to estimates on solutions to the differential equation \eqref{eq:ulb} from Appendix \ref{app:Properties}. Those estimates require $q$ to be bounded away from $0$ and $p(q)$ to be bounded away from $1$ (which does not follow from $q$  being away from $1$). To deal with the former, we will control $H^{(u)}(\bx_u)$ for all $u$ close to the root by a different argument. For the latter, we will relate the correlated Hamiltonians $H^{(u)}$ to modified ones $\bar H^{(u)}$, defined by another correlation profile $\bar p(q)$ such that $\max_{q\in[0,1]} \bar p(q)<1$. This will also entail using a modified energy profile instead of $\uu_{\lb}$ as the bound for the energy on the root-to-leaf path that we produce (analogously to the discussion in Subsection \ref{subsec:pf of reduction}).

Let $\delta>0$ and let $\eta>0$ be a parameter to be determined below.	
Throughout the section, we use the notation $\step:=(1-\theta)/\tau$. We define 
\begin{equation}
	\label{eq:qlqr}
	\begin{aligned}
		q_{l}&:=\max\big\{\theta+i\step:\,0\leq i\step\leq\eta,\; i\in\ZZ_{\ge 0}\big\}\, ,\\
		q_{r}&:=q_l\vee \min\big\{\theta+i\step:\,p(\theta+i\step)\geq1-\eta,\; i\in\ZZ_{\ge 0}\big\}\, .
	\end{aligned}
\end{equation}
Define $\bar{p}(t)=p(t)\wedge p(q_{r})$ and let $\bar{\bF}^{(u)}$
be defined by (\ref{eq:bG(u)}) with $\bar{p}(t)$ instead of $p(t)$,
using the same independent copies $\bF^{e}$ used to define $\bF^{(u)}$.
Similarly to (\ref{eq:correlation-1}), we have that 
\[
\E\big\{\bar{F}_{i}^{(u)}(\bx_1)\bar{F}_{j}^{(v)}(\bx_2)\big\}=\delta_{ij}\bar{p}(\langle u,v\rangle)\xi(\langle\bx_1,\bx_2\rangle).
\]
Denote the corresponding energies by $\bar{H}^{(u)}(\bx):=\frac{1}{2}\|\bar{\bF}^{(u)}(\bx)\|^{2}$.

Similarly to $\uu_{\lb}(t)$,
define $\varphi:[\theta,1]\to\R$ as the solution\footnote{Although $\bar{p}\notin\cuC$ because of the discontinuity at
	$q_{r}$, the solution can be defined in the obvious way. On $[\theta,q_{r}]$,
	it coincides with $\uu_{\lb}$. On $[q_{r},1]$, we proceed by using
	the differential equation (\ref{eq:ulb}) with initial conditions
	$\varphi(q_{r})=\uu_{\lb}(q_{r})$.} of (\ref{eq:ulb}) with $p$ replaced by $\bar{p}$. Define an approximation
of $\varphi(t)$ by 
\begin{equation}\label{eq:phi_delta}
	\varphi_{\delta}(t)=\varphi(t)-\frac{\delta}{4}-\frac{\delta}{4}t\,.
\end{equation}
Recall $\cL(s;p,\dot{p},\uu,t)$ from \eqref{eq:cL}. Note that if $\varphi(t)\leq\uu_{\lb}(t)$ for some $t$, 
since $p(t)\ge \bar{p}(t)$ and   $\dot{p}(t)\ge\dot{\bar{p}}(t)$, for any $s>0$ we have
\begin{align*}
	\cL(s;p,\dot{p},\uu_{\lb},t) & \le \cL(s;\bar{p},\dot{\bar{p}},\varphi,t)+ \frac{1}{2}(\dot{p}(t)-\dot{\bar{p}}(t))\xi(t)+
	\frac{p(t)\dot{\xi}(t)}{2(1+p(t)\dot{\xi}(t)s)}-\frac{\bar{p}(t)\dot{\xi}(t)}{2(1+\bar{p}(t)\dot{\xi}(t)s)}	\,.
\end{align*}
Since $\frac{\de}{\de x}\frac{x}{1+sx}=(1+xs)^{-2}\in(0,1]$ for $x\ge0$,
\begin{align*}
	\frac{p(t)\dot{\xi}(t)}{2(1+p(t)\dot{\xi}(t)s)}-\frac{\bar{p}(t)\dot{\xi}(t)}{2(1+\bar{p}(t)\dot{\xi}(t)s)}	\le \frac{1}{2}(p(t)-\bar{p}(t))\dot{\xi}(t)
\end{align*}
and therefore
\begin{align*}
	\dot{\uu}_{\lb}(t) 
	& \leq\dot{\varphi}(t)+\frac{1}{2}(\dot{p}(t)-\dot{\bar{p}}(t))\xi(t)+\frac{1}{2}(p(t)-\bar{p}(t))\dot{\xi}(t)\\
	&=\dot{\varphi}(t)+\frac{1}{2}\frac{\de\phantom{t}}{\de t}\Big[(p(t)-\bar{p}(t))\xi(t)\Big]\,.
\end{align*}
We conclude that
\[
\uu_{\lb}(1)\leq\varphi(1)+\frac{1}{2}\Big[(p(1)-\bar{p}(1))\xi(1)-(p(0)-\bar{p}(0))\xi(0)\Big]\leq\varphi(1)+\frac{1}{2}\eta\xi(1)\,.
\]
Since $\varphi_{\delta}(1)\geq\varphi(1)-\frac{1}{2}\delta$, assuming that $\eta$ satisfies $\eta\xi(1)\leq\delta/2$,
we have that
\begin{equation}
	\varphi_{\delta}(1)\geq\uu_{\lb}(1)-3\delta/4\ge\ALG-3\delta/4\,.\label{eq:claim1}
\end{equation}

Define the events
\begin{align*}
	\cEr & :=\left\{ \forall u\in\TT\text{ s.t. }\langle u,u\rangle=q_{l}:\,\frac{1}{n}\bar{H}^{(u)}(\bx_{u})<\varphi_{\delta}(q_{l})\right\} \,,\\
	\cEm & :=\Big\{\exists u\in\TT\setminus\LL:\,\langle u,u\rangle\geq q_{l},\\
	& \qquad\frac{1}{n}\bar{H}^{(u)}(\bx_{u})\ge\varphi_{\delta}(\langle u,u\rangle)\quad\text{and}\quad\frac{1}{n\kappa}\sum_{v\in D(u)}\bar{H}^{(v)}(\bx_{v})<\varphi_{\delta}(\langle u,u\rangle+\sigma)\Big\}\,,\\
	\cEbar & :=\Big\{\exists u\in\LL:\,\frac{1}{n}\sup_{\bx\in\Ball^{d}(1)}\big|\bar{H}^{(u)}(\bx)-H^{(u)}(\bx)\big|>\delta/4\Big\}\,.
\end{align*}
We emphasize that the points $\bx_u$ are determined through the algorithm $\calA$ by the construction of Subsection \ref{sec:Hierarchical} with the original $p(q)$ of \eqref{eq:p_linear}.

\begin{lemma}
	\label{lem:3events}Provided that $\eta\xi(1)\leq\delta/2$, 
	\begin{equation}
		\begin{aligned}
			\P\left(H(\calA(\bF))\leq \ALG-\delta\right)^{\kappa^{\tau}}&\le\P\left(\forall u\in\LL:\,\frac{1}{n}H^{(u)}(\calA(\bF^{(u)}))<\ALG-\delta\right)\\
			&\le \P(\cEr)+\P(\cEm)+\P(\cEbar)\,.
		\end{aligned}
		\label{eq:3events}
	\end{equation}
\end{lemma}

\begin{proof}
	The first inequality follows from Lemma \ref{lem:many_to_one_exceedence_prob}. 
	The second follows since if $\cEr^{c}\cap\cEm^c\cap\cEbar^{c}$
	occurs, then $\frac1nH(\calA(\bF^{(u)}))\ge\ALG-\delta$ for some $u\in\LL$.
	Indeed, letting $u_0$ be the vertex $u$ with $\<u,u\>=q_l$ that maximizes $\bar{H}^{(u)}(\bx_{u})$, on $\cEr^{c}$ we have $\frac{1}{n}\bar{H}^{(u_0)}(\bx_{u_0})\ge\varphi_{\delta}(\langle u_0,u_0\rangle)$.
	If $\cEm^{c}$
	also occurs, by defining iteratively $u_{i+1}$ to be the vertex 
	$u\in D(u_i)$ which maximizes $\bar{H}^{(u)}(\bx_{u})$, we construct a
	path $u_{0},u_{1},\ldots,u_{\ell}\in\LL$ such that $\frac{1}{n}\bar{H}^{(u_{i})}(\bx_{u_{i}})\ge\varphi_{\delta}(\langle u_{i},u_{i}\rangle)$.
	In particular, the leaf $u_{\ell}$ satisfies the same inequality.
	Finally, assuming that $\cEbar^{c}$ also occurs, using \eqref{eq:claim1},
	\[
	\frac{1}{n}H^{(u_\ell)}(\calA(\bF^{(u_\ell)})) = \frac{1}{n}H^{(u_{\ell})}(\bx_{u_{\ell}})\ge\varphi_{\delta}(1)-\delta/4\geq\ALG-\delta\,.\qedhere
	\]
\end{proof}

For technical reasons, we treat separately
the case of $\theta$ close to $1$. The next lemma coincides with Theorem \ref{thm:optimalityOC}, except for the additional assumption that $\theta\ge \theta_0$ for $\theta_0$ close to $1$ (which greatly simplifies the proof). We defer its proof to Subsection \ref{subsec:proof_lemma_large_theta}.

\begin{lemma}
	\label{lem:large_theta}
	For any $\xi$, $\alpha$ and $\delta>0$
	there exist constants $c,\lambda>0$ and $\theta_{0}<1$ such that if $\theta\ge\theta_{0}$ and $\calA_d$ is $(\lambda,\nu)$
	overlap concentrated for some $d\ge0$ and $\nu\in[0,1]$, then 
	\[
	\P\left(\frac{1}{n}H(\calA(\bF))\le\ALG-\delta\right)\le2e^{-cd}+2\nu^{c}\,.
	\]
\end{lemma}

In the next subsections, we will prove the following bounds on the probability of the events
above.

\begin{lemma}
	\label{lem:E1}
	There exist  constants $c_{0}=c_{0}(\delta,\xi,\alpha)$
	and $c_{1}=c_{1}(\delta,\xi,\alpha)$ such that if $\lambda,\eta,1/\tau,1/\kappa<c_{0}$ and $\calA$ is $(\lambda,\nu)$
	overlap concentrated, then $\P(\cEr)\leq5 e^{-c_{1}d}+3\kappa^{6\tau}\nu$.
\end{lemma}

\begin{lemma}
	\label{lem:E3}For any $c_{1}>0$ there is a constant $c_{0}=c_{0}(\delta,\xi,\alpha,c_{1})$
	such that if $\eta<c_{0}$, then $\P(\cEbar)\leq\kappa^{\tau}e^{-c_{1}d}$.
\end{lemma}

\begin{lemma}
	\label{lem:E2}
	Assume that  $\theta\leq\theta_{0}$
	for $\theta_{0}$ as defined in Lemma \ref{lem:large_theta}.
	Then, for any $\xi, \alpha, \delta, c_0, c_1>0$, there exist $\lambda,\eta,\tau,\kappa,d_0>0$ such that  $\lambda,\eta,1/\tau,1/\kappa\le c_0$ and, for any $d\ge d_0$, if $\calA_d$
	is $(\lambda,\nu)$ overlap concentrated,
	\begin{align*}
			\P(\cEm)\leq 2\kappa^{\tau}e^{-c_{1}d}+\kappa^{6\tau}\nu\, .
	\end{align*}
\end{lemma}

Before tending to prove the lemmas above, we first show how Theorem \ref{thm:optimalityOC} follows from them. 
Let $\theta_0=\theta_0(\xi,\alpha,
\delta)$ be as defined in Lemma \ref{lem:large_theta}. By Lemma \ref{lem:large_theta}, if $\theta\ge\theta_0$ then
\eqref{eq:Thm_OC} with some $c'$ and $\lambda'$. It remains to prove \eqref{eq:Thm_OC} for any $\theta\le\theta_{0}$ with some $c$ and $\lambda$ since then it will hold for general $\theta$ and $c\wedge c'$ and $\lambda\vee \lambda'$.

Assume that $\theta\le \theta_0$. 
By Lemma
\ref{lem:3events}, in the setting as above,
we have 
\[
\P\left(H(\calA(\bF))\leq\ALG-\delta\right)\le \Big(
\P(\cEr)+\P(\cEm)+\P(\cEbar) \Big)^{1/\kappa^{\tau}}\, .
\]
By Lemma \ref{lem:E1}, \ref{lem:E3}, \ref{lem:E2}, there exists a constant  $c_1=c_1(\delta,\xi\,\alpha)>0$ and $\lambda,\eta,\tau,\kappa,d_0>0$
such that, if $d\ge d_0$ and $\calA_d$ is $(\lambda,\nu)$ overlap concentrated, then 
\[
\P\left(H(\calA(\bF))\leq\ALG-\delta\right)\leq\left((3\kappa^{\tau}+5)e^{-c_1d}+4\kappa^{6\tau}\nu\right)^{1/\kappa^{\tau}}\,.
\]
In fact, those three lemmas allow us to choose $\kappa$ and $\tau$ such that $\kappa^{\tau}$ is as large as we wish. The required bound follows by with $c=c(\delta,\xi,\alpha)$ sufficiently small.

Below we prove the remaining lemmas, and by that complete the proof of Theorem \ref{thm:optimalityOC}. After completing the latter proofs, we prove Theorem \ref{thm:optimalityLip} in Subsection \ref{subsec:pf_optimalityLip}.

\subsection{An auxiliary result}\label{subsec:m0}
We will use the following simple result in the proof of Lemmas \ref{lem:large_theta} and \ref{lem:E1}.
Let $(\tilde\bF^{(u)})_{u\in\TT}$ be an i.i.d. copy of the correlated systems 
$(\bF^{(u)})_{u\in \TT}$. For any leaf $v\in\LL$
define $\by_{v}:=\calA(\tilde\bF^{(v)})$, and for any $u\in\TT\setminus\LL$ define
$\hat\by_{u}$ and $\by_{u}$ similarly to (\ref{eq:xu}), using the points $\by_{v}$ instead of $\bx_v$. Finally, define $\bm_0:=\by_{\varnothing}$ to be the root.

\begin{lemma}
	\label{lem:dist_to_x0}For any $\eps_0>0$, there exists $c=c(\eps_0)>0$
	such that if $\lambda,\eta,1/\kappa\le c$ and $\calA$ is $(\lambda,\nu)$
	overlap concentrated, then
	\begin{align}
		\label{eq:xvx0}
		\P\left(\vphantom{X^2}
		\forall u\in\TT:\:\big|\< \bx_{u},\bm_{0}\> -\theta \big|\leq \eps_0
			\right) &\geq1-3\kappa^{6\tau}\nu\,,\\
			\P\left(\vphantom{X^2}\forall u\in\TT \mbox{\ s.t.\ }\<u,u\>=q_l:\:\|\bx_u-\bm_0\|\leq \eps_0\right)&\geq1-3\kappa^{6\tau}\nu\,.\label{eq:xum0}
	\end{align}
\end{lemma}
\begin{proof}
		Recall that we assume that $\bx_u$ are defined through Eq.~\eqref{eq:xu} using the tree  $\TT(\kappa^3,\tau)$. As in the proof of Lemma \ref{lem:tree_overlap}, denote  $A_{u,v}(t):=\langle u,v\rangle+[-t,t]$. 
	Let $\calE$ be the
		event that for all pairs of leaves $u,v\in\LL(\kappa^3,\tau)$, $\langle\bx_{u},\bx_{v}\rangle\in A_{u,v}(\lambda)$, $\langle\by_{u},\by_{v}\rangle\in A_{u,v}(\lambda)$ and $\langle\bx_{u},\by_{v}\rangle\in A_{\varnothing,\varnothing}(\lambda)$.
	From our assumptions and a union bound, $\P(\calE)\geq1-3\kappa^{6\tau}\nu$.
	Henceforth, we assume that $\calE$ occurs. 
	
	By the same argument and using the same notation as in the proof of Lemma \ref{lem:tree_overlap},
	\begin{equation}
		\forall u,v\in\TT:\:\big|\left\langle \bx_{u},\by_{v}\right\rangle -\theta\big|\leq \eps:=\max_{u,v\in\TT}\Big(\eta_{u,v} \wedge \frac{4B}{\eta_{u,v}}\Big)\,.
	\end{equation}
	In particular, for $v=\varnothing$ we 
	obtain Eq.~\eqref{eq:xvx0} for sufficiently small $c$.

	For small enough  $c$, \eqref{eq:xum0} follows since $\eta\le c$, and for any $u\in\TT$ with
	$\|\bx_{u}\|^2=q_{l}\le\theta+c$,
	\[
	\|\bx_u-\bm_0\|^2 \leq \theta+c  +\theta -2\left\langle \bx_{u},\bm_{0}\right\rangle \leq c+2\eps\,. \qedhere
	\]
	
\end{proof}

\subsection{Proof of Lemma \ref{lem:large_theta}}
\label{subsec:proof_lemma_large_theta}

Let $\bm_0$ be as defined in Subsection \ref{subsec:m0} and let $\bx:=\calA(\bF)$ for $\bF$ independent of $\bm_0$. Since $\|\bx-\bm_{0}\|^{2}= 1+\theta-2\<\bx,\bm_{0}\>$, by \eqref{eq:xvx0},
\begin{align*}
	\P\left(\vphantom{X^2}
	\|\bx -\bm_{0}\|^2 \leq 1-\theta +2\eps_0
		\right) &\geq1-3\kappa^{6\tau}\nu
\end{align*}
for $\eps_0$ of our choice, provided that $\lambda,\eta,1/\kappa<c$ for an appropriate $c=c(\eps_0)$ and $(\lambda,\nu)$
overlap concentrated.

Let $\delta>0$. 
Since $\|\bm_{0}\|^2=\theta$ deterministically and $\bm_0$ is independent of $\bF$, for some $c=c(\xi,\delta)$, 
\[
\P\left(\Big|\frac{1}{n}H(\bm_{0})-\frac{1}{2}\xi(\theta)\Big|\leq\frac{\delta}{2}\right)\geq 1-e^{-cn}.
\]
From Lemma \ref{lem:basicestimates}, taking  $t=1-\theta +2\eps_0$, for appropriate $c_{\xi}$, $c_{\xi}'>0$,
\[
\P\Big(\sup_{\|\by-\bm_{0}\|\leq \sqrt t}\|\bF(\bm_{0})-\bF(\by)\|\le \sqrt{dt}c_{\xi}\Big)\geq 1-e^{-d c_{\xi}'}\,.
\]

Choosing $\eps_0$ sufficiently small and $\theta_0$ sufficiently 
close to $1$ (depending on $\xi,\alpha,\delta$) and assuming that 
$\theta\ge \theta_0$, on the intersection of the three events above, 
$\frac{1}{n}H(\bx)\geq\frac{1}{2}\xi(1)-\delta$. 
Since $\frac{1}{2}\xi(1)\geq\ALG$
by Lemma \ref{lem:U-continuous}, 
\begin{equation}\label{eq:3t}
	\P\left(\frac{1}{n}H(\bx)<\ALG-\delta\right)\leq \P\left(\frac{1}{n}H(\bx)<\frac12\xi(1)-\delta\right)\leq e^{-c(\xi,\delta)n}+e^{-d c_{\xi}'}+3\kappa^{6\tau}\nu\,.
\end{equation}
Denote the right-hand side of \eqref{eq:3t} by $A$. The same bound holds if we replace $A$ by $A^s$ for some $s\in(0,1)$. Indeed, if $A\ge1$ this is trivial and if $A<1$ then $A\le A^s$. Using that $(x+y+z)^s \le x^s+y^s+z^s$ for $x,y,z>0$ and taking $s$ sufficiently small, the lemma follows.
\qed

\subsection{Proof of Lemma \ref{lem:E1}}

As in the statement of the lemma, $c_{0}$ and $c_{1}$ are constants which depend only on $\delta,\xi$
and $\alpha$ and will be assumed to be sufficiently small whenever
needed. 
Let $\TT_{q}:=\left\{ u\in\TT:\,\langle u,u\rangle=q\right\} $.
Since
$\varphi_{\delta}(t)= \varphi(t)-\delta(1+t)/4$
by Eq.~\eqref{eq:phi_delta},
$\varphi(t)\le \uu_{\lb}(t)$, and 
$\uu_{\lb}(t)\le p(t)\xi(t)/2$ by Lemma \ref{lem:U-continuous},
we have 
\[
\varphi_{\delta}(t) \le
\frac{1}{2}p(t)\xi(t)-\frac{\delta}{4}
\;\;\;\;\forall t\in[\theta,1]\, .
\]
Hence, to complete the proof it will be enough to show that 
\begin{equation}
		\P\left(\exists u\in\TT_{q_l}:\,\frac{1}{n}\bar{H}^{(u)}(\bx_{u})\geq\frac{1}{2}p(q_{l})\xi(q_{l})-\frac{\delta}{2}\right)\geq1-(5e^{-c_{1}d}+3\kappa^{6\tau}\nu).\label{eq:H_close_to_root}
\end{equation}

Let $\bm_0$ be as defined in Subsection \ref{subsec:m0}. Fix an arbitrary $u\in\TT_{q_l}$. Since $\bar{\bF}^{(u)}\ed \sqrt{\bar{p}(\langle u,u\rangle)}\bF$,
by Lemma \ref{lem:basicestimates} there exist constants $C_{\xi}$, $C_{\xi}'$ such that 
\begin{equation}
	\P\left(\max_{\bx\neq\by\in\Ball^{d}(1)}\frac{|\bar{H}^{(u)}(\bx)-\bar{H}^{(u)}(\by)|}{\|\bx-\by\|}\leq2dC_{\xi}\right)\geq1-2\exp(-dC_{\xi}')\,.\label{eq:ev2}
\end{equation}
Since $\bar{H}^{(u)}(\bm_{0})\ed\frac{1}{2}\sum_{i=1}^{n}X_{i}^{2}$
for i.i.d. $X_{i}\sim\N(0,\bar{p}(\langle u,u\rangle)\xi(\theta))$,
there exists $c=c(\xi,\delta)$ such that
\begin{equation}
	\P\left(\frac{1}{n}\bar{H}^{(u)}(\bm_{0})\ge\frac{1}{2}\bar{p}(\langle u,u\rangle)\xi(\theta)-\frac{\delta}{8}\right)\ge1-e^{-cn}.\label{eq:ev1}
\end{equation}

On the intersection of the events in Eqs.~(\ref{eq:xum0}), (\ref{eq:ev2})
and (\ref{eq:ev1}),
\begin{align*}
	\frac{1}{n}\bar{H}^{(u)}(\bx_{u}) & \ge\frac{1}{n}\bar{H}^{(u)}(\bm_{0})-\frac{1}{n}2dC_{\xi}\cdot \eps_0\\
	& \geq\frac{1}{2}\bar{p}(q_l)\xi(\theta)-\frac{\delta}{8}-\frac{2}{\alpha}C_{\xi}\eps_0\\
	& \ge\frac{1}{2}p(q_{l})\xi(q_{l})-\frac{\delta}{2},
\end{align*}
where the last inequality follows provided that $\eta\leq \delta/(16 \xi'(1))$ and $2C_{\xi}\eps_0/\alpha\leq \delta/16$. This proves that \eqref{eq:H_close_to_root} holds in the setting of Lemma \ref{lem:E1}.\qed

\subsection{Proof of Lemma \ref{lem:E3}}

For $u\in\LL$, write
\[
H^{(u)}(\bx)-\bar{H}^{(u)}(\bx)=\frac{1}{2}\|\bF^{(u)}(\bx)-\bar{\bF}^{(u)}(\bx)\|^{2}+\langle \bF^{(u)}(\bx)-\bar{\bF}^{(u)}(\bx),\bar{\bF}^{(u)}(\bx)\rangle\,.
\]
Since $p(q_{r})\geq1-\eta$, the process $F_{i}^{(u)}-\bar{F}_{i}^{(u)}$
has covariance function 
\[
\E\big\{[F_{i}^{(u)}(\bx_1)-\bar{F}_{i}^{(u)}(\bx_1)] [F_{i}^{(u)}(\bx_2)-\bar{F}_{i}^{(u)}(\bx_2)]\big\}
= \eta'\xi(\<\bx_1,\bx_2\>)\, ,
\]
for some
$\eta'\leq\eta$. Lemma \ref{lem:E3} thus follows from Lemma \ref{lem:basicestimates}.\qed

\subsection{Proof of Lemma \ref{lem:E2}}

Recall the notation $\step:=(1-\theta)/\tau$. In the current proof, we
let $\bar{\Delta}(\bx_{u},u)$ be as defined in Eq.~\eqref{eq:barT+barDelta},
but using $\bar{H}^{(v)}$ instead of $H^{(v)}$. Accordingly, we
redefine $\Psi(\uu_0,u)$ by 
\begin{equation}
	\Psi(\uu_0,u):=\E\Big[\bar{\Delta}(\bx,\kappa,\eps,q,q',\rho,\rho')\,\Big|\,\frac{1}{n}\bar H^{(0)}(\bx)=\uu_0\Big]\label{eq:Psi}
\end{equation}
as in (\ref{eq:Psi-1}), but here with $q=\langle u,u\rangle$, $q'=q+\step$,
$\rho=\bar{p}(q)$ and $\rho'=\bar{p}(q')$ and $\bx$ such that $\|\bx\|^{2}=q$
(i.e., we replace $p$ by $\bar{p}$ in (\ref{eq:Psi-1})).

Denote $\TT^{q}:=\left\{ \forall u\in\TT(\kappa,\tau):\,\langle u,u\rangle\ge q\right\} $ and define the events
\begin{equation}
	\label{eq:Di}
	\begin{aligned}
		\calD_{1}(C') &=\Big\{\max_{u\in\TT}\max_{\bx\in \Ball^{d}(1)}\frac{1}{n}\bar{H}^{(u)}(\bx)\leq C'\Big\}\,,\\
		\calD_{2}(\eps) &=\left\{ \forall u,v\in\TT^{q_l}:\:\Big|\left\langle \bx_{u},\bx_{v}\right\rangle -\langle u,v\rangle\Big|\leq\eps\right\}\,,\\
		\calD_{3}(\delta_0) &=\left\{ \max_{u\in\TT\setminus(\LL\cup \{\r\}) }\Big|\bar{\Delta}(\bx_{u},u)-\Psi\Big(\frac{1}{n}\bar{H}^{(u)}(\bx_{u}),u\Big)\Big|\leq\delta_{0}d\right\} \,.
	\end{aligned}
\end{equation}

Fix some $\xi, \alpha, \delta, c_0, c_1>0$ and fix an arbitrary $\eta\le c_0\vee 1$. Recall that $p$ and $\theta$  are determined by the algorithm $\calA_d$  and that in the current setting we assume that $\theta\le\theta_0$ for $\theta_{0}=\theta_{0}(\xi,\alpha,\delta)$ as in Lemma \ref{lem:large_theta}. 
Let $C_1=C_1(\xi,C',\eta)\ge0$ be a number to be determined below. We assume the following (which we justify after this list):
\begin{enumerate}
	\item\label{enum:1} We have  
	\begin{equation}\label{eq:sigma_bd}
		q_{l}\in[\theta+\eta/2,\theta+\eta]\qquad\mbox{and}\qquad\sigma\le 1/\tau\le \eta^2/4\le \eta/2\,.
	\end{equation}
	\item\label{enum:1'} Letting $C_{0}=C_{0}(\xi,\alpha,C_1)$
	be the constant as in Proposition \ref{prop:main_bound},
	\begin{equation}\label{eq:1st_delta_bd}
			C_{0}\sigma^{\frac{3}{2}}\vee C_{0}\frac{1}{\sqrt{\kappa}}\sigma^{\frac{1}{2}} \vee C_{0}\kappa\sqrt{\eps}\le \frac{\delta}{60}\sigma\,.
	\end{equation}
	\item
	For any $q\in\{q_{l},q_{l}+\step,\ldots,1-\step\}$ and $x\in[0,\step]$,
	\begin{align}
		|\dot{\varphi}(q+x)-\dot{\varphi}(q)|&<\frac{\delta}{20}\,.\label{eq:dotphi_cont}
	\end{align}
	\item Finally, we choose $\delta_{0}>0$ such that 
	\begin{equation}
		\frac{\delta_{0}}{\alpha}\leq\frac{\delta}{20}\frac{1-\theta_{0}}{\tau}\leq\frac{\delta}{20}\frac{1-\theta}{\tau}=\frac{\delta\sigma}{20}\,.\label{eq:reduction_bd3}
	\end{equation}
\end{enumerate}

We now explain how we choose the parameters so that both the inequalities above and suitable bounds on the probabilities of the events of \eqref{eq:Di} hold. Recall that $\eta$ is already fixed. Take $C'=C'(\xi,\alpha)>0$ sufficiently large so that using Lemma \ref{lem:basicestimates} and a union bound, $\P(\calD_{1}(C')^c)\le \kappa^\tau e^{-c_1d}$ for any $d\ge d_0$. 

We choose $\tau$ sufficiently large so 	
that $1/\tau\le\eta^2/4$ and thus \eqref{eq:sigma_bd} holds. Then, by the definition of $q_l$ and $q_r$ in \eqref{eq:qlqr}, if $q_l<q_r$ then $p(q_r-\sigma)<1-\eta\le p(q_r)$ and, using concavity of $p(t)$,
\[
\frac{p(q_r)-p(q_r-\sigma)}{\sigma}\le \frac{1}{q_l-\theta}\le \frac{1}{\eta/2}\implies p(q_r)-p(q_r-\sigma)\le \frac{2\sigma}{\eta}\le \frac{\eta}{2}\implies p(q_r)\le 1-\frac{\eta}{2}\,.
\]
Otherwise, if $q_l=q_r$, then $\bar p(t)= p(q_l)\ge 1-\eta$ for any $t\in[q_l,1]$.

By Lemma \ref{lem:SecondDerivative}, applied with $\eta/2$, there is
$L(\eta)<\infty$ such that $|\ddot p(t)|\le L(\eta)$ for
$t\in[q_l,q_r)$ when $q_r>q_l$. On $(q_r,1]$, $\bar p$
is constant, 
whence, $|\ddot{\bar p}(t)|\le L(\eta)$ for every
$t\in[q_l,1]$ except possibly $t=q_r$ when $q_r>q_l$. Moreover,
concavity and $q_l-\theta\ge\eta/2$ give
$0\le\dot{\bar p}(t)\le 2/\eta$ on the same intervals (and
$\dot{\bar p}(t)=0$ for $t>q_r$). We set
$C_1 = C'(\xi,\alpha)\vee L(\eta)\vee(2/\eta)$.
If needed, increase $\tau$ so that following hold. First, $C_{0}((1-\theta_0)/\tau)^{1/2}\le \delta/60$
so that the bound on the first term in \eqref{eq:1st_delta_bd} holds. 
Second, using Lemmas 
\ref{lem:phidot-Lip} and \ref{lem:phidot-Lip-constp},  Eq.~\eqref{eq:dotphi_cont}
holds. 

Next, choose $\delta_0$ that satisfies \eqref{eq:reduction_bd3}. We define $\eps = 20/(\kappa^3\tau\sqrt{q_l})$ and take $\kappa$ sufficiently large so that 
\[
\frac{C_{0}}{\sqrt{\kappa}} \le \frac{\delta}{60}\Big(\frac{1-\theta_0}{\tau}\Big)^{1/2}\qquad\mbox{and}\qquad C_{0}\kappa\sqrt{\eps}\le \frac{\delta}{60}\Big(\frac{1-\theta_0}{\tau}\Big)\,,
\]
which is possible because $q_l\ge\eta/2$ and
$\kappa\sqrt{\eps}=\sqrt{20}/(\sqrt{\kappa\tau}\,q_l^{1/4})$.
The bounds on the second and third terms in \eqref{eq:1st_delta_bd} follow. Finally, we take $\lambda<1/(\kappa^3\tau)\wedge c_0$ so that (recalling that we use $p(t)$ and $\TT(\kappa^3,\tau)$ to define the points $\bx_u$) the explicit bound in Lemma \ref{lem:tree_overlap} gives $\P(\calD_{2}(\eps)^c)\le \kappa^{6\tau}\nu$.

Of course, we can also require the parameters to satisfy the conditions of Theorem \ref{thm:concentration} with $\delta_0$, $\bar p$, and $c=c_1$ so that, for $d\ge d_0$, $\P(\calD_{3}(\delta_0)^c)\le \kappa^{\tau}e^{-c_1d}$. Combining this bound with those on the other two events from \eqref{eq:Di}, we have that
\begin{align*}
	\P(\cEm) & \le\P(\cEm\cap\calD_{1}(C')\cap\calD_{2}(\eps)\cap\calD_{3}(\delta_0))+\P(\calD_{1}(C')^{c})+\P(\calD_{2}(\eps)^{c})+\P(\calD_{3}(\delta_0)^{c})\\
		& \leq\P(\cEm\cap\calD_{1}(C')\cap\calD_{2}(\eps)\cap\calD_{3}(\delta_0))+2\kappa^\tau e^{-c_1d}+\kappa^{6\tau}\nu\,.
\end{align*}
To complete the proof, we will show that the intersection above is empty. 

Assume towards contradiction that the intersection is not empty and assume that it occurs. Since $\calD_{2}(\eps)$ occurs, for any $u\in\TT^{q_l}\setminus\LL$, $(\bx_{u}: v\in D(u)) \in T(\bx_{u},u)$ and thus
\begin{equation*}
	\frac{1}{\kappa}\sum_{v\in D(u)}\bar H^{(v)}(\bx_{v})\ge \bar \Delta(\bx_{u},u).
\end{equation*}
Since $\cEm$ occurs, for some $u$ with $q:=\<u,u\>\in \{q_{l},q_{l}+\step,\ldots,1-\step\}$,
\[
\frac{1}{n}\bar{H}^{(u)}(\bx_{u})\ge\varphi_{\delta}(q)\quad\text{and}\quad\frac{1}{n}\bar{\Delta}(\bx_{u},u)<\varphi_{\delta}(q+\sigma)\,.
\]
Because $\calD_{1}(C')$ and $\calD_{3}(\delta_0)$ occur, for the same $u$ and some $\uu_0\in [\varphi_{\delta}(q)\vee0,C']$, 
\[
\frac{1}{n}\Psi(\uu_0,u)=\frac{1}{n}\Psi\big(\frac{1}{n}\bar{H}^{(u)}(\bx_{u}),u\big)<\varphi_\delta(q+\sigma)+\frac{\delta_{0}}{\alpha}\le \varphi_{\delta}(q+\sigma)+\frac{\delta\sigma}{20}\,,
\]
where for the last inequality we used \eqref{eq:reduction_bd3}.

We have that
\begin{align*}
	\frac1n\Psi(\uu_0,u)&\overset{(a)}{\ge}\phi_0+
		\sigma \drift(q,\phi_0;\bar p,\dot{\bar p})-C_{0}\Big[\sigma^{\frac{3}{2}}+A_d\sigma+\frac{1}{\sqrt{\kappa}}\sigma^{\frac{1}{2}}+\kappa\sqrt{\eps}\Big]\\
	&\overset{(b)}{\ge} \phi_0+
	\sigma \drift(q,\phi_0;\bar p,\dot{\bar p})-\frac{\delta\sigma}{10}\\
	&\overset{(c)}{\ge} \varphi_{\delta}(q)+
	\sigma \drift(q,\varphi(q);\bar p,\dot{\bar p})-\frac{\delta\sigma}{10}\,,
\end{align*}
where  $(b)$ holds for $d\ge d_0$ using \eqref{eq:1st_delta_bd} and for $(c)$  we used the monotonicity of $\drift$ in its second argument.
By Proposition \ref{prop:main_bound}, $(a)$ follows since we are working with $q\in \{q_{l},q_{l}+\step,\ldots,1-\step\}$ and therefore: $q'=q+\sigma$ and $q'/q\le (\eta/2+\sigma)/(\eta/2)\le 2$ by \eqref{eq:sigma_bd}, $\uu_0\le C'$, $|\dot{\bar p}(t)|\le 2/\eta$, and $|\ddot{\bar{p}}(t)|\le L(\eta)$ on $(q,q')$.

By combining the two inequalities above and using \eqref{eq:phi_delta}, we have that
\begin{align*}
	\sigma \drift(q,\varphi(q);\bar p,\dot{\bar p}) \le \varphi(q+\sigma) - \varphi(q) -\frac{\delta\sigma}{4} +\frac{\delta\sigma}{20}+\frac{\delta\sigma}{10} \le \dot \varphi(q)\sigma - \frac{\delta\sigma}{20}
	\,,
\end{align*}
where the second inequality follows by \eqref{eq:dotphi_cont}.
This contradicts the definition of $\varphi$ (which implies 
$ \dot{\varphi} = \drift(q,\varphi;\bar{p},\dot{\bar{p}})$), thus
the proof is completed.	 \qed

\subsection{\label{subsec:pf_optimalityLip}Proof of Theorem \ref{thm:optimalityLip}}

Assume that $\calA$ is $M/\sqrt{d}$-Lipschitz for some fixed $M>0$
and any $d$. Let $\delta>0$. By Lemma \ref{lem:overlapConc} below, for
any $\lambda>0$, $\calA$  is $(\lambda,2e^{-\frac{d\lambda^{2}}{8M}})$
overlap concentrated. By Theorem \ref{thm:optimalityOC}, there exist
$c,\lambda>0$ (depending on $\xi$, $\alpha$ and $\delta$) such
that 
\[
\P\left(\frac{1}{n}H(\calA(\bF))\le\ALG-\delta\right)\le2e^{-cd}+2^{1+c}e^{-c\frac{d\lambda^{2}}{8M}}\,.
\]
Theorem \ref{thm:optimalityLip} follows for appropriate $C=C(\delta,\xi\,\alpha,M)$, assuming the correctness of Lemma \ref{lem:overlapConc}.
\qed

\begin{lemma}
	\label{lem:overlapConc}If an algorithm $\calA=(\calA_{d})_{d\geq1}$
	is $(M/\sqrt{d})$-Lipschitz for some $M>0$, then for any $\lambda>0$
	it is $(\lambda,2e^{-\frac{d\lambda^{2}}{8M}})$ overlap concentrated.
\end{lemma}
Up to the fact that in our setting the input comes from multiple functions
$F_{i}(\bx)$ instead of a single function, this is essentially the
content of Proposition 8.2 of \cite{huang2024optimization}.
\begin{proof}
	Let $(\bF^{1},\bF^{2})$ be $\rho$-correlated copies of $\bF$ and
	let $(\bG^{1},\bG^{2})$ be the corresponding Gaussian coefficients
	of the polynomial systems. Note that in law 
	\[
	(\bG^{1},\bG^{2})\ed (\sqrt{\rho}\bA^{0}+\sqrt{1-\rho}\bA^{1},\sqrt{\rho}\bA^{0}+\sqrt{1-\rho}\bA^{2}),
	\]
	for i.i.d. copies $\bA^{i}$ of $\bG$. Hence, the lemma will follow
	by Gaussian concentration if we can show that 
	\begin{equation}
		\Psi(\bA^{0},\bA^{1},\bA^{2}):=\langle\calA(\sqrt{\rho}\bA^{0}+\sqrt{1-\rho}\bA^{1}),\calA(\sqrt{\rho}\bA^{0}+\sqrt{1-\rho}\bA^{2})\rangle\label{eq:Psi_pcorrelated}
	\end{equation}
	is $(2M/\sqrt{d})$-Lipschitz, where here and below
	$\bA^{i}$ (and $\bB^{i}$) are treated as deterministic 
	arrays\footnote{Note that concentration of Lipschitz function holds
		for infinite-dimensional Gaussian measures, see e.g. discussion above 
		\cite[Theorem 2.7]{Ledoux}.}.

	Indeed, since for $\by,\bx'\in\Ball^d(1)$ and $\bx,\by'\in\R^d$, 
	\begin{align*}
		|\langle\bx,\by\rangle-\langle\bx',\by'\rangle| & \leq|\langle\bx,\by\rangle-\langle\bx',\by\rangle|-|\langle\bx',\by\rangle-\langle\bx',\by'\rangle|\leq\|\bx-\bx'\|+\|\by-\by'\|,
	\end{align*}
	we have that
	\begin{align*}
		& \big|\Psi(\bA^{0},\bA^{1},\bA^{2})-\Psi(\bB^{0},\bB^{1},\bB^{2})\big|\\
		& \quad\leq\sum_{i=1,2}\|\calA(\sqrt{\rho}\bA^{0}+\sqrt{1-\rho}\bA^{i})-\calA(\sqrt{\rho}\bB^{0}+\sqrt{1-\rho}\bB^{i})\|\\
		& \quad\leq\frac{M}{\sqrt{d}}\sum_{i=1,2}\left\Vert \sqrt{\rho}(\bA^{0}-\bB^{0})+\sqrt{1-\rho}(\bA^{i}-\bB^{i})\right\Vert .\\
		& \quad\leq\frac{M}{\sqrt{d}}\Big(2\left\Vert (\bA^{0}-\bB^{0})\right\Vert +\left\Vert (\bA^{1}-\bB^{1})\right\Vert +\left\Vert (\bA^{2}-\bB^{2})\right\Vert \Big).\qedhere
	\end{align*}
\end{proof}

\section*{Statement of AI use}

GPT-5.5 was used to prove correctness of
 the solution of the ODE in 
Theorem \ref{thm:MatchingSol}.
GPT-5.6 Sol was used to prove Lemmas \ref{thm:RS-lower-bound},
\ref{lemma:U-achieved-D}, \ref{lemma:U-achieved}, 
\ref{lem:U-continuous-in-D}, and
\ref{lemma:U-strict-inequality-gen}. 
 The authors verified and edited the
 output of these steps.

\section*{Acknowledgments}

AM was partially supported by the NSF Awards MFAI-2501597 and MFAI-2601871.
ES was supported by the Israel Science
Foundation (Grant Agreement No. 2055/21) and European Union (ERC grant PolySpin, Grant Agreement No. 101165541).

\appendix

\section{Smoothness estimates}

Define $\|\nabla^{k}H(\bx)\|_{\op}=\max_{\|\bu\|=1}\langle\nabla^{k}H(\bx),\bu^{\otimes k}\rangle$ and let 
$\partial_{\bu}^{k}H(\bx)=\langle\nabla^{k}H(\bx),\bu^{\otimes k}\rangle$  denote the directional derivative for $\bu\in\R^{d}$.
For $k\geq1$, recall that $\nabla^{k}\bF(\bx):=((\partial_{j_{1}}\cdots\partial_{j_{k}}F_{i}(\bx))_{i\leq n,j_{1},\ldots,j_{k}\leq d}\in\R^{n}\otimes(\R^{d})^{\otimes k}$
and that its operator norm is defined by (\ref{eq:tensor_op_norm}). For $k=0$, set  $\nabla^{0}\bF(\bx)=\bF(\bx)$ and $\nabla^{0}H(\bx)=H(\bx)$.
The following estimates from  \cite{montanari2024smale}
are useful to us. 
\begin{lemma}
	\label{lem:basicestimates}For any mixture $\xi$ and $\alpha\in(0,1]$
	there exist constants $C(\xi)$ and $C_{*}(\xi)$ such that, for any
	$t>0$, with probability at least $1-2\exp(-dt^{2}/C_{*}(\xi))$ we
	have 
	\begin{align}
		\forall 0\le k\leq3:\,\,\max_{\bx\in\Ball^{d}(1)}\|\nabla^{k}\bF(\bx)\|_{\op} & \leq\sqrt{d}(C(\xi)+t)\,,\label{eq:nablaF}\\
		\forall 0\le k\leq3:\,\,\max_{\bx\in\Ball^{d}(1)}\|\nabla^{k}H(\bx)\|_{\op} & \leq 16 d(C(\xi)^2+t^{2})\,.\label{eq:nablaH}
	\end{align}
\end{lemma}

\begin{proof}
	We have that 
	\begin{align}\nonumber
	\|\nabla^{k}H(\bx)\|_{\op} &\le 2^k\max_{0\le i\le k}\| \nabla^{i}\bF(\bx)\|_{\op}\| \nabla^{k-i}\bF(\bx)\|_{\op}\,,\\
	\max_{\bx\in\Ball^d(1)}	\max_{0\le k\le 3}\|\nabla^{k}H(\bx)\|_{\op} &\le 8 \Big(	\max_{\bx\in\Ball^d(1)}\max_{0\le k\le 3}\| \nabla^{k}\bF(\bx)\|_{\op}\Big)^2=:8M^2\,.\label{eq:nablaHnablaF}
	\end{align}
Hence, the bound of \eqref{eq:nablaF} implies that of \eqref{eq:nablaH}, and we only need to prove the required bound on the probability of the former.

Note that for unit vectors $\bu,\bv$, $\langle\nabla^{k}\bF(\bx),\bv\otimes\bu^{\otimes k}\rangle$
has the same law as $\partial_{\bu}^{k}F_{1}(\bx)$.
Using 
\[
\E\left\{ \frac{\partial}{\partial x_{i_{1}}}\cdots\frac{\partial}{\partial x_{i_{k}}}F_{1}(\bx)\frac{\partial}{\partial y_{j_{1}}}\cdots\frac{\partial}{\partial y_{j_{m}}}F_{1}(\by)\right\} =\frac{\partial}{\partial x_{i_{1}}}\cdots\frac{\partial}{\partial x_{i_{k}}}\frac{\partial}{\partial y_{j_{1}}}\cdots\frac{\partial}{\partial y_{j_{m}}}\E\left\{ F_{1}(\bx)F_{1}(\by)\right\} \,,
\]
one easily sees that the variance of $\partial_{\bu}^{k}F_{1}(\bx)$ is bounded by $c_{k}\sum_{i=1}^{2k}\xi^{(i)}(1)$ for some absolute $c_k$ and any $\bx\in\Ball^d(1)$.
Hence, $M$ is the maximum of a Gaussian field with variance bounded by $C_*(\xi)/2$ for some constant $C_*(\xi)$. By the Borell-TIS inequality, to complete the proof we only need to show that $\E M\le \sqrt d C(\xi)$ for appropriate $C(\xi)$. This follows from Lemma 2.1 of 
\cite{montanari2024smale}.
\end{proof}

\begin{lemma}
	\label{lem:basicestimates_condtiional}For any mixture $\xi$, $\alpha\in(0,1]$ and $0\le \rho \le \rho'\le1$, with $\bF^{(0)}$ and $\bF^{(1)}$ defined by \eqref{eq:correlation} and the same constants $C(\xi)$ and $C_{*}(\xi)$ as in Lemma \ref{lem:basicestimates}, we have the following.
	Let $\bx_0\in\Ball^d(1)$ and $\uu_0\ge0$. Conditional on $\bF^{(0)}(\bx_0)=\bv$ for some $\bv$ with $\|\bv\|=\sqrt{2n\uu_0}$, for any
	$t>0$, with probability at least $1-2\exp(-dt^{2}/C_{*}(\xi))$ we
	have 
	\begin{align}
		\forall 0\le k\leq3:\,\,\max_{\bx\in\Ball^{d}(r)}\|\nabla^{k}\bF^{(1)}(\bx)\|_{\op} & \leq\sqrt{d}(C(\xi)+t + \sqrt{2\uu_0} D(\xi,r,\|\bx\|^2))\,,\label{eq:nablaF2}\\
		\forall 0\le k\leq3:\,\,\max_{\bx\in\Ball^{d}(r)}\|\nabla^{k}H^{(1)}(\bx)\|_{\op} & \leq 24 d(C(\xi)^2+t^{2}+2\uu_0 D(\xi,r,\|\bx\|^2)^2)\,,\label{eq:nablaH2}
	\end{align}
	where $D(\xi,r,q):=\xi(q)^{-1}\max_{0\le k\le 3} \xi^{(k)}(rq^{1/2})q^{k/2}$.\footnote{If $q=0$, we implicitly assume that $\xi(0)>0$. Otherwise, $H^{(0)}=0$ a.s. and one can use the unconditional Lemma \ref{lem:basicestimates} instead.}
\end{lemma}
\begin{proof}
	Let $\bv$ and $\uu_0$ be as in the lemma. Recall that conditional on $\bF^{(0)}(\bx_0)=\bv$, 
	the process $\bF^{(1)}(\bx)$ has the law of $\tilde{\bF}^{(1)}(\bx)+\Phi(\langle\bx,\bx_0\rangle)\bv$
	for some Gaussian process $\tilde{\bF}^{(1)}(\bx)$ and $\Phi(r)=\xi(r)/\xi(\|\bx_0\|^2)$
	(see \eqref{eq:condexp1}). By Anderson's inequality, $\sup_{\bx\in\Ball(r)}\|\tilde{\bF}^{(i)}(\bx)\|$
	and $\max_{\by\in\Ball^{d}(r)}\|\nabla^{k}\tilde{\bF}^{(i)}(\bx)\|_{\op}$
	are stochastically dominated by the same expressions with the unconditional
	processes $\bF^{(i)}(\bx)$. Moreover, 
	\begin{align*}
	\max_{0\le k\le 3}\sup_{\bx\in\Ball(r)}\big\|\nabla_{\bx}^{k}\big[\Phi(\langle\bx,\bx_0\rangle)\bv\big]\big\|_{\op} &= \max_{0\le k\le 3}\sup_{\by\in\Ball(r)}\big\|\frac{\xi^{(k)}(\<\bx,\bx_0\>)}{\xi(\|\bx_0\|^2)}\bv\otimes \bx_0^{\otimes k}\big]\big\|_{\op}  = \sqrt{2n\uu_0} D(\xi,r,\|\bx\|^2)\,.
	\end{align*}

Thus, the required bound on the probability of \eqref{eq:nablaF2} follows from \eqref{eq:nablaF}. Using that $(x+y+z)^2\le 3(x^2+y^2+z^2)$ for non-negative numbers and \eqref{eq:nablaHnablaF}, \eqref{eq:nablaF2} implies \eqref{eq:nablaH2}.
\end{proof}

\section{Properties of the optimal control problem}
\label{app:Properties}

Throughout this appendix we assume $\xi(0)>0$, and $t\mapsto \xi(t)$
is non-constant. We write either $\dot{\xi}$ or $\xi'$ for its derivatives. 

\subsection{Useful identities}
\label{app:UsefulIdentities}
Recall that  $\drift: \reals^4\to\reals$ 
and $\cL:\reals^4\to\reals$ are  defined by
\begin{align}
\drift(t,u;p,r)
&:=\sup_{s>0}\cL(s;p,r,u,t),\label{eq:drift-sup}\\
\cL(s;p,r,u,t)
&:=-(r\dot\xi(t)+p\ddot\xi(t))(u\vee0)s+
\frac12 r\xi(t)-\frac{1}{2\alpha s}
+\frac{p\dot\xi(t)}{2(1+p\dot\xi(t) s)}
 (1-\dot p\xi(t) s).\label{eq:cL-app}
\end{align}
For $p:[0,1]\to[0,1]$, we introduce the shorthand
\begin{align}
a=a(t):=p(t)\xi'(t),\qquad
\qquad
b=b(t):=(p\xi)'(t) \, ,
\label{eq:ab-def}
\end{align}
so that
\begin{align}
\cL(s;p(t),\dot p(t),u,t)
&:=-a'(t)(u\vee0)s+
\frac12b(t)-\frac{1}{2\alpha s}
-\frac{a(t)b(t)s}{2(1+a(t) s)}\, .\label{eq:cL-L-identity}
\end{align}
Define
\begin{align}
L(p(t),\dot p(t),u,t)
&:=\min_{\omega\geq0}\mathcal R(\omega;p(t),\dot p(t),u,t),\label{eq:L-omega}\\
\mathcal R(\omega;p(t),\dot p(t),u,t)
&:=\frac{1}{\sqrt\alpha}\sqrt{2(u\vee 0)a'(t)+a(t)\omega^2}
+\frac12\big(\omega-\sqrt{b(t)}\big)^2.\label{eq:R-omega}
\end{align}
Then the two variational representations are related by
\begin{align}
\drift(t,u;p(t),\dot p(t))
=\frac12(p\xi)'(t)-L(t,u;p(t),\dot p(t)).
\label{eq:drift-L}
\end{align}
Indeed, the elementary identity
\begin{align}
\frac{1}{\sqrt\alpha}\sqrt{x}
=\inf_{s>0}\left\{\frac{s x}{2}+\frac{1}{2\alpha s}\right\},
\qquad x\geq0,
\label{eq:sqrt-variational}
\end{align}
gives (without loss of generality $u\ge 0$)
\begin{align}
L(p(t),\dot p(t),u,t)
&=\inf_{s>0,\,\omega\geq0}
\left\{ua'(t)s+\frac{a(t)s}{2}\omega^2+\frac{1}{2\alpha s}
+\frac12(\omega-\sqrt{b(t)})^2\right\}\nonumber\\
&=\inf_{s>0}
\left\{ua'(t)s+\frac{1}{2\alpha s}
+\frac{a(t)b(t)s}{2(1+a(t)s)}\right\},
\label{eq:L-s}
\end{align}
where in the second line we minimized over $\omega$, setting
$\omega=\sqrt b/(1+as)$.  Comparing with 
Eq.~\eqref{eq:cL-L-identity} yields
\begin{align}
\frac b2-\cL(s;p(t),\dot p(t),u,t)
=ua'(t)s+\frac{1}{2\alpha s}+\frac{a(t)b(t)s}{2(1+a(t)s)}.
\end{align}
In what follows, we will often omit the arguments
of $p(t),\dot p(t)$, $\xi(t)$, and so on, when they are clear from the context.
\subsection{Definitions}

For a non-decreasing function $p:[0,1]\to[0,1]$ with $p(1)=1$, let
\begin{align}
\theta(p):=\inf\big\{t\in[0,1]:p(t)>0\big\}.
\label{eq:theta-p-relaxed}
\end{align}
We define
\begin{align*}
\cuD&:=\Big\{p:[0,1]\to[0,1]:\;p\mbox{ is non-decreasing and right-continuous}, p(1)=1\Big\},\\
\Ccont&:= \big\{ p\in\cuD :\; 
p\big|_{[\theta(p),1]}\in C([\theta(p),1])\big\}\, ,
\end{align*}
and
\begin{align*}
\Cccv:=\Big\{p\in\Ccont:\;
p\big|_{[\theta(p),1]}\mbox{ is concave}\Big\}.
\end{align*}
Thus a function in $\Ccont$ may have a jump from zero to
$p(\theta(p))>0$. 
For any $p\in\cuD$, $\dot p \ge 0$ denotes its almost-everywhere derivative,
with an arbitrary value, say zero, assigned on the remaining null set.
Note that
$\int_0^1\dot p(t)\,\de t\le1$.
We denote by $\de p$ the measure on $[0,1]$ with distribution function $p$,
and therefore $\de p=\dot p\de t+\de p_{\rm s}$, where $\de p_{\rm s}$ is
 singular with respect to the Lebesgue measure.

We note of course that 
\begin{align}
&\cuC_*\subseteq \Cccv\subseteq \Ccont \subseteq \cuD\, .
\end{align}

For $p\in\cuD$, set $\theta=\theta(p)$. We define $u(t;p)=0$ for
$t<\theta$ and, on $[\theta,1]$, define $u(t;p)$ as the solution of
the following integral equation:
\begin{align}\label{eq:ODE-integral}
u(t;p) = \frac{1}{2}p(t)\xi(t)-\int_\theta^t 
L(p(s),\dot{p}(s),u(s;p),s) \, \de s\, .
\end{align}
We refer to Lemma \ref{lem:U-continuous} for existence and uniqueness of the solution.
Equation \eqref{eq:ODE-integral} is also 
equivalent to the measure identity
\begin{align}
\de u(t)=\drift(t,u(t);p(t),\dot p(t))\,\de t
+\frac{\xi(t)}2\,\de p_{\rm s}(t),
\qquad u(\theta)=\frac12p(\theta)\xi(\theta).
\label{eq:ODE-measure}
\end{align}
By construction $u(t)-p_{\rm s}(t)\xi(t)/2$ is continuous 
and of bounded variation on $[\theta,1]$.
If $p|_{[\theta,1]}$ is absolutely continuous, the singular term vanishes
and we recover the ordinary differential equation
\begin{align}
\dot{u}(t)=\drift(t,u(t);p(t),\dot p(t))\, ,\;\;\;
u(t_0) =u_0\, ,\label{eq:ODE-repeated}
\end{align}
with $t_0=\theta$ and $u_0=\frac{1}{2}p(\theta)\xi(\theta)$, understood
in the Carath\'eodory sense \cite{walter2013ordinary}.

Finally, for any subset of functions $\cuG\subseteq \cuD$, we let
\begin{align}
U_*(\cuG,\alpha)&:=\inf\big\{u(1;p)  :\; p\in \cuG\big\}\, , \;\;\;\;\;
U(p,\alpha) = u(1;p)\, ,
\end{align}
where we emphasized the dependence of $U,U_*$ on $\alpha$:
we will omit this argument whenever not needed.

\subsection{Basic properties}

\begin{lemma}\label{lemma:ConcavityLB}
Let $p\in\Cccv$, $\theta=\theta(p)$.
If $\theta<1$, then, for every $t\in[\theta,1]$,
\begin{align}
	p(\theta)+(1-p(\theta))\frac{t-\theta}{1-\theta}
	\le p(t)\le 1\,.
\end{align}
For $t\in(\theta,1)$, its one-sided derivatives also satisfy
\begin{align}
 \frac{1-p(t)}{1-t}\le \pd(t+)\le \pd(t-)
 \le \frac{p(t)-p(\theta)}{t-\theta}
 \le \frac{1-p(\theta)}{t-\theta}\le\frac{1}{t-\theta}\, .
\label{eq:ConcavityLB}
\end{align}
\end{lemma}
\begin{proof}
These inequalities follow immediately from concavity, using
the endpoints $(\theta,p(\theta))$ and $(1,1)$. 
\end{proof}

\begin{lemma}\label{lem:SecondDerivative}
For any $\eta>0$ there exists a positive constant $L(\eta)<\infty$ such that
the following holds. 
Let $p\in\cuC_*$ and $\theta=\theta(p)$. For any $\eta\in(0,1)$ and
$t\in(\theta,1)$ such that $t-\theta\ge\eta$ and $p(t)\le 1-\eta$,
we have
\begin{align}
-L(\eta)\le \ddot p(t)\le 0\, .
\label{eq:SecondDerivativeBound}
\end{align}
\end{lemma}
\begin{proof}
Let $\chi(q)=\sum_{k\ge0}\chi_kq^k$ be as in Definition
\ref{def:ControlF} and set $q=p(t)$. The function $\chi$ is nonconstant
 strictly increasing and convex on
$[0,1]$. Since $p(t)=\chi^{-1}(t)$ for $t-\theta\ge\eta$ and $p(t)\le 1-\eta$,
it follows that $p$ is concave, whence the upper bound in Eq.~\eqref{eq:SecondDerivativeBound} holds.

To prove the lower bound note that, by the implicit function theorem, we have
\begin{align}
	\dot p(t)=\frac{1}{\chi'(q)},
	\qquad
	\ddot p(t)=-\frac{\chi''(q)}{\chi'(q)^3},\;\;\; \mbox{ for } q=p(t)\in (0,1)\, 
	\label{eq:p-inverse-derivatives}
\end{align}
Moreover, $\dot p(t)\le \eta^{-1}$, and therefore 
$\chi'(q)\ge \eta$ by Lemma \ref{lemma:ConcavityLB}.
It is therefore sufficient to prove $\chi''(q)\le \tilde{L}(\eta)$.
Using
$q\le1-\eta$ and $\sum_{k\ge0}\chi_k=\chi(1)\le1$, we obtain
\begin{align*}
\chi''(q)
	&=\sum_{k\ge2}k(k-1)\chi_kq^{k-2}\\
	&\le \sup_{k\ge2}k(k-1)(1-\eta)^{k-2}=:\tilde{L}(\eta).
\end{align*}
It is immediate to check that $\tilde{L}(\eta)<\infty$ as claimed.
\end{proof}
\begin{lemma}\label{lemma:drift-continuous}
\begin{enumerate}
\item The map $(t,u,p,r)\mapsto\drift(t,u;p,r)$ is $C^{\infty}$
at every point of
$[0,1]\times \reals_{>0}\times [0,1]\times \reals_{\ge 0}$ at which
\begin{align}
	r\xi'(t)+p\xi''(t)>0.
\end{align}
For every fixed $(t,p,r)$, the map
$u\mapsto\drift(t,u;p,r)$ is $C^{\infty}$ on $\reals_{>0}$.
\item If $p\in\cuD$, and 
 $\pd$ exists in a neighborhood of $t\in [0,1]$
 nd is continuous at
$t$, and $u\ge 0$,
 then the map
 $(v,s)\mapsto\drift(s,v;p(s),\dot p(s))$
 is continuous at $(u,t)$.
 \item  If $p\in \cuD$ and $\pd(t)$ exists finite,
 then
 $u\mapsto\drift(t,u;p(t),\dot p(t))$
 is continuous at $0$.
Further,
\begin{align}\label{eq:drift-0}
	\drift(t,0;p(t),\dot p(t))
	= \frac{1}{2}\left(\sqrt{(p\xi)'(t)}-\sqrt{\frac{p(t)\xi'(t)}{\alpha}}\right)_+^2\, .
	\end{align} 
\item Assume $p\in \Ccont$ and define 
	\begin{align}
		J_{++}(\alpha,p) &:=
 \big\{t\in (0,1]: \;\pd(t)<\infty\, ,\;\;
 (p\xi')(t)/\alpha
< (p\xi)'(t)\big\}\, ,\\
			J_{-}(\alpha,p) &:= \big\{t\in (0,1]: \; \pd(t)<\infty\, ,\;\;
			(p\xi')(t)/\alpha
		> (p\xi)'(t)\big\}\, .
	\end{align}
		Then: 
		\begin{enumerate}
		\item[$(a)$] 
		If $t\in J_{++}(\alpha,p)$ and $p(t)\xi'(t)>0$,
		then $\drift(t,0;p(t),\dot p(t))> 0$ and
		 $u\mapsto\drift(t,u;p(t),\dot p(t))$ is continuously differentiable on
		$\reals_{\ge 0}$.
		\item[$(b)$] For $t\in J_{-}(\alpha,p)$,
		$u\mapsto\drift(t,u;p(t),\dot p(t))$, we have
		$\drift(t,u;p(t),\dot p(t))
		=-C_*(t;p)\sqrt{u}+o(\sqrt{u})$ as $u\to 0$,
	for some $C_*(t,p)\ge 0$. Moreover, $C_*(t,p)>0$ if
	$(p\xi')'(t)>0$, while $\drift(t,u;p(t),\dot p(t))=0$ for all
	$u\geq0$ if $(p\xi')'(t)=0$.
\end{enumerate}
\item  Assume  $p\in\Ccont$ with $\lim\inf_{t\to 0+} \pd(t)>0$ and either $p(0)=0$ or $\alpha>1$. 
Letting  
$t_+(\alpha,p) := \sup\{t\ge 0: \; (0,t)\subseteq
J_{++}(\alpha,p)\}$,
we have $t_+(\alpha,p)>0$ strictly.
\item Assume $p\in \cuD$. Then there exists a constant $C=C(u_1,\xi,\alpha)$
such that, for all $0<u_1<u_2$, at every differentiability point of
$p$ we have
\[
\sup_{u\in [u_1,u_2]}
\left|\partial_u\drift(t,u;p(t),\dot p(t))\right|
\le
C\big(1+\sqrt{\pd(t)}\big).
\]
\end{enumerate}
\end{lemma}
\begin{proof}

\noindent{\bf Point 1.}
For fixed $(t,p,r)$, set
\begin{align}
a=p\xi'(t),\qquad c=r\xi'(t)+p\xi''(t),\qquad
b=r\xi(t)+p\xi'(t).
\end{align}
By Eq.~\eqref{eq:L-s}, for $u>0$,
\begin{align}
L(p,r,u,t)=\inf_{s>0}\Phi(s;t,u,p,r),\qquad
\Phi(s;t,u,p,r)
=ucs+\frac{1}{2\alpha s}+\frac{abs}{2(1+as)}.
\label{eq:Phi-s-smoothness}
\end{align}
If $c>0$, the minimizer $s_*$ is unique and is characterized by
\begin{align}
\frac{1}{\alpha}
=2uc s_*^2+\frac{ab s_*^2}{(1+as_*)^2}.
\label{eq:s-star-smoothness}
\end{align}
Indeed, the right-hand side is strictly increasing from zero to
infinity. Its derivative with respect to $s_*$ is strictly positive,
so the implicit function theorem shows that $s_*$, and hence the
value in Eq.~\eqref{eq:Phi-s-smoothness}, depends smoothly on
$(t,u,p,r)$. Equation~\eqref{eq:drift-L} proves the first claim.
It also proves smoothness in $u>0$ for fixed $(t,p,r)$ when $c>0$.
When $c=0$, the variational formula is independent of $u$, which
proves the remaining claim.

\noindent{\bf Point 2.}
For $u\in \reals_{\ge 0}$, the minimizer $\omega_*$ of $\mathcal R$
is unique (by strong convexity) and belongs to
$[0,\sqrt{(p\xi)'(t)}]$.  At a continuity point of $\dot p$, all the
coefficients in Eq.~\eqref{eq:R-omega} vary continuously and the
compact interval containing the minimizer is locally bounded.  The
maximum theorem therefore gives the claimed continuity.

\noindent{\bf Point 3.}
For a fixed $t$, all the coefficients in
Eq.~\eqref{eq:R-omega} are finite.  The same compactness argument as
in Point 2 shows continuity as $u\downarrow0$.  At $u=0$ the
variational functional becomes
\begin{align}
\mathcal R(\omega;p,\dot p,0,t)
=\sqrt{\frac a\alpha}\,\omega
+\frac12(\omega-\sqrt b)^2.
\end{align}
Its minimizer is
$\omega_*=(\sqrt b-\sqrt{a/\alpha})_+$, and direct substitution into
Eq.~\eqref{eq:drift-L} gives Eq.~\eqref{eq:drift-0}.

\noindent{\bf Point 4.}
Smoothness on $\reals_{>0}$ follows from Point 1.  For part $(a)$,
use the notation $a,b,c$ from Point 1.  The assumptions give
$b>a/\alpha$ and $a>0$, so the minimizer at $u=0$ is the strictly
positive number $\omega_0=\sqrt b-\sqrt{a/\alpha}$.  In a
neighborhood of $(\omega_0,0)$, the stationarity condition is
\begin{align}
\omega_*
=\sqrt b-\frac{1}{\sqrt\alpha}
\frac{a\omega_*}{\sqrt{2cu+a\omega_*^2}}.
\label{eq:FOC-omega}
\end{align}
Its derivative with respect to $\omega_*$ at $(\omega_0,0)$ is one.
The implicit function theorem therefore extends $\omega_*$ smoothly
to $u=0$.  Equation~\eqref{eq:drift-0} gives strict positivity of the
drift there.

For part $(b)$, first suppose that $c=(p\xi')'(t)>0$.  For each
fixed $M<\infty$, uniformly over $z\in[0,M]$,
\begin{align}
\lim_{u\downarrow0}\frac{1}{\sqrt u}
\left\{\mathcal R(z\sqrt u;p,\dot p,u,t)-\frac b2\right\}
&=G(z),\\
G(z)&=\frac{1}{\sqrt\alpha}\sqrt{2c+az^2}-\sqrt b\,z.
\end{align}
Since $b<a/\alpha$, the function $G$ tends to infinity as
$z\to\infty$, satisfies $G'(0)<0$, and has a unique minimizer
$z_*>0$.  If $\omega_u$ minimizes $\mathcal R$ at $u$, comparison
with $\omega=0$ and the bound
\begin{align}
\mathcal R(\omega;u,t)-\frac b2
\ge \left(\sqrt{\frac a\alpha}-\sqrt b\right)\omega
+\frac{\omega^2}{2}
\end{align}
show that $\omega_u=O(\sqrt u)$.  Thus the locally uniform
convergence applies at the minimizer and yields
\begin{align}
\drift(t,u;p(t),\dot p(t))
=-C_*(t;p)\sqrt u+o(\sqrt u),
\qquad C_*(t;p)=G(z_*).
\end{align}
Moreover,
\begin{align}
G(z)&\ge\frac{1}{\sqrt\alpha}
\left(\sqrt{2c+az^2}-\sqrt a\,z\right)\\
&=\frac{1}{\sqrt\alpha}
\frac{2c}{\sqrt{2c+az^2}+\sqrt a\,z}>0,
\end{align}
and hence $C_*(t;p)>0$.  If instead $c=0$, the functional
$\mathcal R$ is independent of $u$.  Since $b<a/\alpha$, its
minimizer is $\omega_*=0$ and $L=b/2$, so the drift vanishes
identically.  This proves part $(b)$.

\noindent{\bf Point 5.} We have
$(p\xi)'-\alpha^{-1}p\xi'
=(1-\alpha^{-1})p\xi' +p'\xi$.
If $p(0)=0$, the possibly negative first term is $o(1)$, whereas the
liminf of the second one is strictly positive. If $\alpha>1$, both
terms are non-negative and the second again has strictly positive
liminf. Thus $(p\xi)'>(p\xi')/\alpha$ throughout a sufficiently
small right neighborhood of zero.

\noindent{\bf Point 6.} 
If $a'=0$, the drift is independent of $u$ and the claim is
immediate.  Otherwise, applying the envelope formula 
to Eq.~\eqref{eq:drift-L}, and using
 $u\ge u_1$ gives
 \begin{align}
 0\ge \partial_u \drift(t,u;p,\dot p) = 
 -\frac{1}{\sqrt{\alpha}}
 \frac{a'}{\sqrt{2u a'+a\omega_*^2}}
 \ge - \sqrt{\frac{a'}{2\alpha u_1}}
 \ge - C\big(1+\sqrt{p'(t)}\big),
 \end{align}
where we used $a'=p'\xi'+p\xi''\le C(1+p')$.
\end{proof}

\begin{lemma}\label{lemma:Sqrt-t}
	For $p\in\cuD$, $u\geq0$, and $t\in(0,1]$, set
	\begin{align}
	\Lambda(p,u,t):=\frac{1}{\sqrt\alpha}
	\sqrt{\xi'(t)\big(2u+p(t)\xi(t)\big)}.
	\end{align}
	For every $M<\infty$ there exists a constant
	$C=C(M,\xi,\alpha)<\infty$ such that, for any
	$p\in \Ccont$, $u\in[0,M]$, and every $t\in (0,1]$ at which
	$\dot p(t)$ exists, $0\le L(p,\dot p,u,t)$
	\begin{align}\label{eq:BasicLineq}
	-C+\Lambda(p,u,t)\sqrt{\dot p(t)}
	\le L(p,\dot p,u,t)&\le C+\Lambda(p,u,t)\sqrt{\dot p(t)}\\
	&\le C\big(1+\sqrt{\dot p(t)}\big).
	\end{align}
	Further, if $p\in \Cccv$, then, for almost every $t\in(0,1]$,
	\begin{align}\label{eq:BasicLineq-2}
	L(p,\dot p,u,t)\le
	\frac{C}{\sqrt{t-\theta(p)}}
	\mathbbm{1}_{\{t>\theta(p)\}}.
	\end{align}
	Here and below, the right-hand side is understood to be zero when
	$t\leq\theta(p)$. 

	Finally, the constant $C$ above depends uniquely on 
	$\xi,\alpha$ if $u$ solves Eq.~\eqref{eq:ODE-integral}.
\end{lemma}
\begin{proof}
Fix a point at which $\dot p(t)$ exists, and abbreviate
$r=\dot p(t)$.  Since $p$ is non-decreasing, $r\geq0$.  With the
notation of Eq.~\eqref{eq:ab-def},
\begin{align}
2ua'+ab
&=r\xi'(t)\big(2u+p(t)\xi(t)\big)
+2up(t)\xi''(t)+p(t)^2\xi'(t)^2.
\label{eq:2uaplusab}
\end{align}
In particular, all the terms on the right-hand side are
non-negative.  Choosing $\omega=\sqrt b$ in
Eq.~\eqref{eq:R-omega} gives
\begin{align}
L(p,r,u,t)\le A:=\frac{1}{\sqrt\alpha}\sqrt{2ua'+ab}.
\label{eq:L-upper-A}
\end{align}

We next compare this upper bound with the value at the minimizer.
The map
\begin{align}
\omega\longmapsto \frac{1}{\sqrt\alpha}
\sqrt{2ua'+a\omega^2}
\end{align}
is $\sqrt{a/\alpha}$-Lipschitz on $\reals_{\geq0}$.  Hence, for
every $\omega\geq0$,
\begin{align}
\mathcal R(\omega;p,r,u,t)
&\geq A-\sqrt{\frac a\alpha}\,|\omega-\sqrt b|
+\frac12|\omega-\sqrt b|^2
\geq A-\frac{a}{2\alpha}.
\end{align}
Taking the infimum over $\omega$ and using $a\leq\|\xi'\|_\infty$
yields
\begin{align}
A-\frac{\|\xi'\|_\infty}{2\alpha}
\leq L(p,r,u,t)\leq A.
\label{eq:L-A-comparison}
\end{align}
Non-negativity also follows directly from Eq.~\eqref{eq:R-omega}.

By Eq.~\eqref{eq:2uaplusab},
\begin{align}
\Lambda(p,u,t)\sqrt r
\leq A
&\leq \Lambda(p,u,t)\sqrt r
+\frac{1}{\sqrt\alpha}
\sqrt{2M\|\xi''\|_\infty+\|\xi'\|_\infty^2}.
\end{align}
Moreover,
\begin{align}
\Lambda(p,u,t)
\leq \frac{1}{\sqrt\alpha}
\sqrt{\|\xi'\|_\infty
\big(2M+\|\xi\|_\infty\big)}.
\end{align}
Together with Eq.~\eqref{eq:L-A-comparison}, these bounds prove
Eq.~\eqref{eq:BasicLineq} after increasing $C$.

Finally, suppose that $p\in\Cccv$ and write $\theta=\theta(p)$.
For almost every $t\in(\theta,1]$, Lemma
\ref{lemma:ConcavityLB} gives $\dot p(t)\leq1/(t-\theta)$.
Since $t-\theta\leq1$, Eq.~\eqref{eq:BasicLineq} therefore gives
\begin{align}
L(p,\dot p,u,t)\leq \frac{C}{\sqrt{t-\theta}}.
\end{align}
For $t<\theta$, we have $p(t)=\dot p(t)=0$ almost everywhere, and
Eq.~\eqref{eq:R-omega} gives $L=0$.  This proves
Eq.~\eqref{eq:BasicLineq-2}.

Finally, if $u$ solves Eq.~\eqref{eq:ODE-integral}, 
then $u(t;p)\le p(t)\xi(t)/2$ and therefore 
we can take $M= \xi(1)/2$,
whence $C$ depends uniquely on $\xi,\alpha$.
\end{proof}

\begin{lemma}\label{lemma:Concavity-in-pdot}
	Assume $\alpha>0$, $u\ge 0$, $p\ge 0$.
	Then the map
	\[
	r\mapsto L(p,r,u,t)
	\]
	is concave on $[0,\infty)$. Further, for $u>0$ it is strictly concave
	on $(0,\infty)$.
	\end{lemma}
	
	\begin{proof}
	By Eq.~\eqref{eq:L-inf-s}, $L(p,r,u,t)$ is the infimum of affine
	functions of $r$ and hence
	$r\mapsto L(p,r,u,t)$ is concave on $[0,\infty)$.

	Note that, for $u,r>0$ the infimum is always achieved for
	some $s_*>0$ because the objective function in
	Eq.~\eqref{eq:L-inf-s} diverges at $s\downarrow 0$ and
	$s\uparrow\infty$. Further, we can rewrite 
the above as 
\begin{align}
L(p,r,u,t)
	&=
	\inf_{b>0} \big\{br+ L^*(b;p,u,t)\}\, ,\label{eq:L-inf-b}\\
 L^*(b;p,u,t) &:=
 \frac{1}{2\alpha s(b)}
 +
 up\xiddot s(b)
 +
 \frac{(p\xidot)^2s(b)}{2(1+p\xidot s(b))}\, ,
\end{align}
where $s(b)\in (0,\infty)$ is the unique solution to the equation
\begin{align}
	b = u\xidot s
	+
	\frac{p\xi\xidot s}{2(1+p\xidot s)} =: f(s)\, .
\end{align}
The solution $s(b)$ is unique because the map $s\mapsto f(s)$
is strictly increasing and smooth with $f(0+) = 0$ and $f(\infty) = \infty$.
Further, since $f'(s)>0$ strictly for $s>0$, 
it follows that $s(b)$ is smooth, and therefore 
$b\mapsto L^*(b;p,u,t)$ is smooth. Since the infimum in 
Eq.~\eqref{eq:L-inf-b} is achieved for some $b_*=f(s_*)$,
it follows from standard properties of Legendre duality that
$r\mapsto L(p,r,u,t)$ is strictly concave.
	\end{proof}

\subsection{Existence and uniqueness of the ODE \eqref{eq:ODE-integral}}

\begin{lemma}\label{lem:U-continuous}
		For any $p\in \cuD$, the integral equation~\eqref{eq:ODE-integral}
		has a unique solution $u(\, \cdot\, ;p)$ on $[\theta(p),1]$,
		which we extend to $[0,1]$ by setting $u(t;p)=0$ for $t<\theta(p)$.
		This solution is such that, letting
		\begin{align}
		u(t;p)= \frac{1}{2}p(t)\xi(t)-v(t;p),
		\end{align}
		we have $v(\,\cdot\,;p)\ge 0$ and $v(\,\cdot\,;p)\in \AC([0,1])$.
		If $p\in\Ccont$, then $u(\,\cdot\,;p)$ is continuous and of
		bounded variation on $[\theta(p),1]$.
		Consequently, the function $u(1;\,\cdot\,): \cuD\to \reals_{\ge 0}$ 
	is well-defined, with $u(1;p)\le \xi(1)/2$.
\end{lemma}
\begin{proof}
Set $\theta=\theta(p)$.
Since $p$ is non-decreasing, its derivative $r=\dot p$ exists and is
non-negative almost everywhere, and $r\in L^1([0,1])$.  We assign
an arbitrary value, say zero, to $r$ on the remaining null set.  Put
\begin{align}
g(t):=\frac12p(t)\xi(t),\qquad
F(t,z):=L\big(p(t),r(t),g(t)-z,t\big).
\label{eq:v-vector-field}
\end{align}
Then Eq.~\eqref{eq:ODE-integral} is equivalent, under the change of
variables $u=g-v$, to
\begin{align}
v(t)=\int_\theta^t F(s,v(s))\,\de s,
\qquad v(\theta)=0.
\label{eq:v-integral-equation}
\end{align}

We verify the hypotheses for this Carath\'eodory equation.  The
representation \eqref{eq:R-omega} shows that $z\mapsto F(t,z)$ is
continuous for almost every $t$, while $t\mapsto F(t,z)$ is
measurable.  Moreover, $L$ is non-decreasing in its third argument,
because every function being minimized in \eqref{eq:R-omega} has
this property.  Consequently, $F(t,\,\cdot\,)$ is non-increasing.
Further, by Lemma \ref{lemma:Sqrt-t}, we have
\begin{align}
0\le F(t,z)\le C\big(1+\sqrt{r(t)}\big).
\label{eq:v-integrable-envelope}
\end{align}
Here $C$ may depend on the bounded interval but not on $t$.
The right-hand side is integrable, since $r\in L^1([0,1])$.
Carath\'eodory's existence theorem (see
\cite[Section 10 and Theorem XVIII]{walter2013ordinary}) therefore gives a local
absolutely continuous solution of \eqref{eq:v-integral-equation}.

We next show that every such solution satisfies
\begin{align}
0\le v(t)\le g(t),\qquad t\in[\theta,1].
\label{eq:v-bounds}
\end{align}
The first inequality follows from $F\ge0$.  For the second, note that
$g$ is non-decreasing and right-continuous, and its derivative exists
and equals $b/2$ almost everywhere, whence
\begin{align}
\int_s^t\frac{b(q)}2\,\de q\le g(t)-g(s),
\qquad \theta\le s\le t\le1.
\label{eq:g-ac-part}
\end{align}
By taking $\omega=0$ in Eq.~\eqref{eq:R-omega},
\begin{align}
L(p,r,0,t)\le \frac b2.
\label{eq:L-zero-upper}
\end{align}
If $v(t)>g(t)$ at some time, set
\begin{align*}
s:=\sup\big\{q\in[\theta,t]:v(q)\le g(q)\big\}.
\end{align*}
The set is nonempty because $v(\theta)=0\le g(\theta)$.  Continuity
of $v$ and monotonicity and right-continuity of $g$ imply that
$s<t$, $v(s)\le g(s)$, and $v(q)>g(q)$ for $q\in(s,t]$.  
For $q\in (s,t]$, the positive-part convention
in Eq.~\eqref{eq:R-omega} gives
$F(q,v(q))=L(p(q),r(q),0,q)$.  Equations
\eqref{eq:g-ac-part} and \eqref{eq:L-zero-upper} therefore imply
\begin{align*}
v(t)-v(s)\le \int_s^t\frac{b(q)}2\,\de q
\le g(t)-g(s),
\end{align*}
contradicting $v(t)>g(t)$.  This proves \eqref{eq:v-bounds}.
Equivalently,
\begin{align}
0\le u(t;p)=g(t)-v(t)\le \frac12p(t)\xi(t).
\label{eq:u-upper-bound}
\end{align}
In particular, the solution remains bounded, since $F$
has an integrable bound \eqref{eq:v-integrable-envelope},
the solution can be extended
to all of $[\theta,1]$ by \cite[Theorem XVIII]{walter2013ordinary}).

If $v_1,v_2$ are two solutions, the monotonicity of $F$ gives, almost
everywhere,
\begin{align*}
\frac{\de}{\de t}|v_1(t)-v_2(t)|^2
=2(v_1-v_2)\big(F(t,v_1)-F(t,v_2)\big)\le0.
\end{align*}
Since both solutions vanish at $\theta$, they coincide.  Thus
$u=g-v$ is the unique solution of Eq.~\eqref{eq:ODE-integral}; by
Eq.~\eqref{eq:v-bounds}, it also satisfies $u\ge0$.

For $t<\theta$, we have $p(t)=u(t;p)=0$, and hence $v(t;p)=0$.
Together with $v(\theta;p)=0$ and
$v\in\AC([\theta,1])$, this shows that the extension of $v$ belongs
to $\AC([0,1])$. If $p\in\Ccont$, then $g=p\xi/2$ is continuous
and non-decreasing on $[\theta,1]$, so $u=g-v$ is continuous and of
bounded variation there. Finally, $u(1;p)\ge0$, so in particular
$u(1;\,\cdot\,):\Ccont\to\reals_{\ge0}$ is well-defined.

The inequality $u(1;p)\le \xi(1)/2$ follows from $v(\,\cdot\,;p)\ge 0$
and $p(1)=1$.
\end{proof}

\begin{corollary}\label{coro:scaling}
Let $p,\ovp\in\Ccont$ and assume that
$\ovp(t)= c\,p(t)$ for some constant $c>0$ and all
$t\in [0,t_1]$.
Then the corresponding solutions of the integral equation~\eqref{eq:ODE-integral} are related by
 $u(t;\ovp) = c\cdot u(t;p)$ for all $t\in [0,t_1]$.
\end{corollary}
\begin{proof}
By Eq.~\eqref{eq:R-omega}, we have
\begin{align}
	\mathcal R(\sqrt{c}\omega;\ovp,\dot{\ovp},c u,t)
	=c\mathcal R(\omega;p,\dot p,u,t),
\end{align}
and therefore
$L(\ovp,\dot{\ovp},cu,t)=cL(p,\dot p,u,t)$.
If $p$ vanishes throughout $[0,t_1]$, so does $\ovp$, and the claim
is immediate. Otherwise their activation times coincide on this
interval; write their common value as $\theta$. Multiplying
Eq.~\eqref{eq:ODE-integral} by $c$ gives, for $\theta\le t\le t_1$,
\[
c u(t;p)=\frac12\ovp(t)\xi(t)
-\int_\theta^t L(\ovp(s),\dot{\ovp}(s),c u(s;p),s)\,\de s.
\]
Uniqueness therefore gives $u(t;\ovp)=c u(t;p)$, and both sides
vanish before $\theta$. This argument also includes the initial
jump and any singular continuous variation of the controls.
\end{proof}

\begin{corollary}\label{coro:U-monotone}
The functions $\alpha \mapsto U(p,\alpha)$ and $\alpha \mapsto U_*(\cuG,\alpha)$ 
are non-decreasing for any $p\in\cuD$ and $\cuG\subseteq \cuD$.
Further, for any $p\in\cuD$ 
\begin{align}
U(p,\alpha)\ge \frac{1}{2}p(t)\xi(t)-\frac{1}{\sqrt\alpha}\sqrt{\xi(1)\xi'(1)}\, .
\end{align}
\end{corollary}
\begin{proof}
Again by 
Eq.~\eqref{eq:L-omega} with $\omega=\sqrt{b}$,
we have 
\begin{align}
	L(p,\dot p,u,t) \le \frac{1}{\sqrt\alpha}\sqrt{2a'u+ab}\le  \frac{1}{\sqrt\alpha}
	\sqrt{\frac{\de\phantom{t}}{\de t}(p^2\xi\xi')}\, ,
\end{align} 
where the last inequality follows from the bound 
$u(t;p) \le p(t)\xi(t)/2$
(for $u(t;p)$ solution of Eq.~\eqref{eq:ODE-integral}).

By Cauchy-Schwartz,
\begin{align*}
	u(1;p) &\ge \frac{1}{2}p(t)\xi(t)-\frac{1}{\sqrt\alpha}\int_0^1
	\sqrt{\frac{\de\phantom{t}}{\de t}(p^2\xi\xi')} \de t\\
	&\ge \frac{1}{2}p(t)\xi(t)-\frac{1}{\sqrt\alpha}
\sqrt{p^2\xi\xi'(1)-p^2\xi\xi'(0)}\\
&\ge \frac{1}{2}p(t)\xi(t)-\frac{1}{\sqrt\alpha}\sqrt{\xi(1)\xi'(1)}\, .
\end{align*}
\end{proof}

\begin{lemma}\label{eq:ContInalpha}
	Recalling that $U_*(\cuG;\alpha)=\inf_{p\in\cuG}U(p,\alpha)$, $\cuG\subseteq \cuD$,
	we have $\alpha \mapsto U_*(\cuG;\alpha)$ is non-decreasing and continuous.
	\end{lemma}
	\begin{proof}
	For any $p\in\cuD$, $u_1(t ;p)$, $u_2(t ;p)$ the solutions of 
	Eq.~\eqref{eq:ODE-integral}  corresponding to parameters $\alpha=\alpha_1$, $\alpha=\alpha_2$, respectively,
	with $\alpha_1<\alpha_2$. With an abuse of notation, we denote by
	$\drift(t,u;\alpha)$, $L(u,t;\alpha)$, $\cR(\omega,u,t;\alpha)$ the drift functions defined in 
	Section \ref{app:UsefulIdentities}, where weomit the dependence on $p$ and $\pd$ 
	because they are fixed, and instead emphasize the dependence on $\alpha$.
	
	It follows from Eq.~\eqref{eq:drift-sup} that
	$\drift(t,u;\alpha_1)\le \drift(t,u;\alpha_2)$ for all $t\in[0,1]$ and $u$.
	It follows therefore that $u_1(t)\le u_2(t)$ for all $t\in[0,1]$.
	
	To prove a bound in the other direction, note that, for any 
	$s\in [0,1]$ and $u\in [0,\xi(t)p(t)/2]$,  letting $\omega_2\in
	\argmin_{\omega\in\reals_{\geq0}} \cR(\omega,u,s;\alpha_2)$, we have
	using $\omega_2\le \sqrt{b}$,
	\begin{align}
	L(s,u;\alpha_1) - L(s,u;\alpha_2)& \le \cR(\omega_2,u,s;\alpha_1) - \cR(\omega_2,u,s;\alpha_2)\\
	&= \left(\alpha_1^{-1/2}-\alpha_2^{-1/2} \right)
	\sqrt{2ua'+ a\omega_2^2}\\
	&\le \left(\alpha_1^{-1/2}-\alpha_2^{-1/2} \right)
	\sqrt{(p\xi)a'+ a(p\xi)'} = \,
	\end{align}
	Further, since
	$\drift$ is non-increasing in $u$, we havr
	\begin{align}
	u_2(t) - u_1(t) &=\int_0^t \big\{\drift(s,u_2(s);\alpha_2) - \drift(s,u_1(s);\alpha_1)\big\} \de s\\
	&\le \int_0^t \big\{\drift(s,u_1(s);\alpha_2) - \drift(s,u_1(s);\alpha_1)\big\} \de s\\
	&= \int_0^t \big\{L(s,u_1(s);\alpha_1) - L(s,u_1(s);\alpha_2)\big\} \de s\\
	&\le \left(\alpha_1^{-1/2}-\alpha_2^{-1/2} \right) 
	\int_0^t \sqrt{(a p\xi)'} \de s \,.
	\end{align}
	We thus concluded that
	\begin{align}
	0\le U(p;\alpha_2) - U(p;\alpha_1) \le C(\xi)\left(\alpha_1^{-1/2}-\alpha_2^{-1/2} \right) \, .
	\end{align}
	Since these bounds hold uniformly in $p$, they hold also for
	$U_*(\cuG;\alpha_2) - U_*(\cuG;\alpha_1)$.
	\end{proof}

\subsection{The infima $\inf_{p\in\cuD}U(p)$, $\inf_{p\in\Cccv}U(p)$ are achieved }

\begin{lemma}\label{lemma:U-achieved-D}
Recall that $U(p):=u(1;p)$.
There exists $p_0\in\cuD$ such that 
\begin{align}
\inf_{p\in\cuD}U(p)=U(p_0)\, .
\end{align}
\end{lemma}
\begin{proof}
Let $(p_n)$ be a minimizing sequence for $\inf_{p\in\cuD}U(p)$ and
write
\begin{align*}
r_n(t)&:=\dot p_n(t),&
g_n(t)&:=\frac12p_n(t)\xi(t),&
u_n(t)&:=u(t;p_n),\\
v_n(t)&:=g_n(t)-u_n(t),&
f_n(t)&:=L(p_n(t),r_n(t),u_n(t),t).
\end{align*}
As usual, $r_n$ is defined arbitrarily on the null set on which
$p_n$ is not differentiable.  
Further, we set $f_(t) = 0$ for $t\in[0,\theta_n)$,
where $\theta_n=\inf\{t\in[0,1]: p_n(t)>0\}$.
Lemma \ref{lem:U-continuous} gives
\begin{align}
v_n(t)=\int_0^t f_n(s)\,\de s,\qquad
0\le u_n(t)\le g_n(t).
\label{eq:D-vn-definition}
\end{align}
Lemma \ref{lemma:Sqrt-t}  and
$\int_0^1r_n\le1$ imply that there exists 
a constant $C(\xi,\alpha)$ such that, for every measurable
$E\subseteq[0,1]$,
\begin{align}
\int_E f_n(t)\,\de t\le \sqrt{|E|}\sqrt{\int_0^1f_n(t)^2\,\de t}
\le C(\xi,\alpha)\sqrt{|E|}.
\label{eq:D-uniform-integrability}
\end{align}
In particular, $(f_n)$ is uniformly integrable. As a consequence of the last inequality,
\begin{align*}
|v_n(t)-v_n(s)|
\le C\sqrt{|t-s|}.
\end{align*}
Thus, after passing to a subsequence, $v_n$ converges uniformly to
some function $v$ by the Arzel\'a-Ascoli theorem.  By \eqref{eq:D-uniform-integrability} and the
Dunford--Pettis theorem the family $(f_n)$ is relatively weakly compact in $L^1([0,1])$, 
which means that (further refining the above subsequence),  we can assume 
that $\int_0^1f_n(t) \varphi(t)\,\de t\to \int_0^1f(t) \varphi(t)\,\de t$ 
for some   $f\in L^1([0,1])$ and all $\varphi\in L^{\infty}([0,1])$.  
Passing to the limit in
\eqref{eq:D-vn-definition} gives
\begin{align}
v(t)=\int_0^t f(s)\,\de s.
\label{eq:D-v-limit}
\end{align}
In particular, $v$ is absolutely continuous.

By Helly's selection theorem, after taking a further subsequence,
$p_n\to p_{\infty}$ at every continuity point of $p_{\infty}$ 
for some non-decreasing right continuous $p_{\infty}$,
with $p_{\infty}(1)=1$.
 In particular, $p_n\to p_{\infty}$ almost everywhere
and in $L^1([0,1])$.

Further, recalling that $r_n=\dot{p}_n$, the measures $r_n(t)\,\de t$ have mass
at most one. Therefore, by further refining the subsequence, 
we can assume they have a weak limit $\nu$.
We can decompose $\nu$ into a singular part $\nu_{\rm s}$ and an absolutely continuous part 
$\rho(t)\de t$. Summarizing (recall that we denote by $\Rightarrow$ the weak convergence of measures)
\begin{align}
r_n(t)\,\de t\ \Rightarrow \nu
=\rho(t)\,\de t+\nu_{\rm s}\, .\label{eq:Rn-WeakLimit}
\end{align}
Since $r_n(t)\,\de t\le \de p_n$ and $\de p_n\Rightarrow \de p_{\infty}$,
we have $\nu\le \de p_{\infty}$.  Consequently,
\begin{align}
0\le \rho(t)\le \dot p_{\infty}(t) 
\quad\text{for almost every }t.
\label{eq:D-density-bound}
\end{align}

Set  $u:=p_{\infty}\xi/2-v$, so that
$u_n\to u\ge0$ almost everywhere. We claim that
\begin{align}
f(t)\le L(p_{\infty}(t),\rho(t), u(t),t)
\quad\text{for almost every }t.
\label{eq:L-upper-semicontinuity}
\end{align}
In order to prove this, for  $K\ge1$, put $r_n^K=\min\{r_n,K\}$. Lemma
\ref{lemma:Sqrt-t} gives
$f_n\le L(p_n,r_n^K,u_n,\cdot)+Cr_n/\sqrt K$.
Fix a continuity point $t$ of $p_{\infty}$;
$u$ is also continuous there.
For $I_h=(t-h,t+h)$, letting $\nu_n(I):=\int_Ir_n(s)\,\de s$, we have
\begin{align}
\frac{1}{|I_h|}\int_{I_h}f_n(s)\,\de s 
& \stackrel{(a)}{\le} \frac{1}{|I_h|} \int_{I_h}L(p_n(s),r_n^K(s),u_n(s),s)\, 
\de s+\frac{C\nu_n(I_h)}{|I_h|\sqrt K}\nonumber\\
& \stackrel{(b)}{\le}  \frac{1}{|I_h|} \int_{I_h}L(p_{\infty}(s),r_n^K(s),u(s),s)\, 
\de s+\omega_K(n^{-1})+\frac{C\nu_n(I_h)}{|I_h|\sqrt K}\nonumber\\
& \stackrel{(c)}{\le}  \frac{1}{|I_h|} \int_{I_h}L(p_{\infty}(t),r_n^K(s),u(t),t)\, 
\de s+\omega_K(h)+\omega_K(n^{-1})+\frac{C\nu_n(I_h)}{|I_h|\sqrt K}\nonumber\\
& \stackrel{(d)}{\le}  L\left(p_{\infty}(t),\frac{\nu_n(I_h)}{|I_h|},u(t),t\right)\, 
\de s+\omega_K(h)+\omega_K(n^{-1})+\frac{C\nu_n(I_h)}{|I_h|\sqrt K}\, ,
\label{eq:IntF_UpperBound}
\end{align}
where $\omega_{C}(\eps)$ denote a function such that $\lim_{\eps\to0}\omega_{C}(\eps)=0$
at $C$ fixed.
Here  $(a)$ follows from the bound $f_n\le L(p_n,r_n^K,u_n,\cdot)+Cr_n/\sqrt K$,
$(b)$ from convergence of $p_n$ to $p_{\infty}$ and $u_n$ to $u$, 
by dominated convergence, $(c)$ from continuity of $L$ and $u$, $p_{\infty}$,
$(d)$ from concavity of $L$ in its second argument,
cf. Lemma \ref{lemma:Concavity-in-pdot}.

By Eq.~\eqref{eq:Rn-WeakLimit}, we have $\nu_n(I_h)\to \nu(I_h)$
whenever $\nu(\partial I_h)=0$. Further $|I_h|^{-1}\int_{I_h}f_n(q)\,\de q
\to |I_h|^{-1}\int_{I_h}f(q)\,\de q$ as $n\to\infty$, by weak convergence
of $f_n$ to $f$. Therefore, taking the $n\to\infty$ limit in Eq.~\eqref{eq:IntF_UpperBound},
and using the continuity of $L$ in its second argument (see Lemma \ref{lemma:drift-continuous}), 
we get
\begin{align*}
\frac{1}{|I_h|}\int_{I_h}f(q)\,\de q
\le L\left(p_{\infty}(t),\frac{\nu(I_h)}{|I_h|},u(t),t\right)
+\omega_K(h)+\frac{C\nu(I_h)}{\sqrt K\,|I_h|}.
\end{align*}
For almost every $t$, Lebesgue differentiation gives
$|I_h|^{-1}\int_{I_h}f\to f(t)$ and
$\nu(I_h)/|I_h|\to\rho(t)$, since $\nu_{\rm s}$ has zero density.
Letting first $h\downarrow0$ and then $K\to\infty$ proves
Eq.~\eqref{eq:L-upper-semicontinuity}. Together with
Eq.~\eqref{eq:D-density-bound} and monotonicity in the second argument,
this yields
\begin{align}
v'(t)=f(t)
&\le L(p_{\infty}(t),\rho(t),u(t),t)\nonumber\\
&\le L\Big(p_{\infty}(t),\dot p_{\infty}(t),\frac{1}{2}p_{\infty}(t)\xi(t)/2-v(t),t\Big)
\quad\text{for almost every }t.
\label{eq:D-subsolution}
\end{align}

Let $w=p_{\infty}\xi/2-u(\,\cdot\,;p_{\infty})$.  By Lemma
\ref{lem:U-continuous}, $w$ is the unique solution of
\begin{align*}
w'(t)=L\Big(p_{\infty}(t),\dot p_{\infty}(t),\frac{1}{2}p_{\infty}(t)\xi(t)-w(t),t\Big),
\qquad w(0)=0.
\end{align*}
The right-hand side is non-increasing in $w$ by Eqs.~\eqref{eq:L-omega}, 
\eqref{eq:R-omega}.  Hence the usual scalar
comparison argument applied to \eqref{eq:D-subsolution} gives
\begin{align*}
\frac{\de}{\de t}(v-w)_+^2\le0
\end{align*}
almost everywhere.  Thus $v(t)\le w(t)$ for every $t$.  Since
$p_n(1)=p_{\infty}(1)=1$ and
$v_n(1)\to v(1)$,
\begin{align*}
U(p_{\infty})
=\frac12\xi(1)-w(1)
\le\frac12\xi(1)-v(1)
=\lim_{n\to\infty}U(p_n)
=\inf_{p\in\cuD}U(p).
\end{align*}
The reverse inequality follows from $p_{\infty}\in\cuD$, proving the claim.
\end{proof}

\begin{lemma}\label{lemma:U-achieved}
Recall that $U(p):=u(1;p)$.
There exists $p_*\in\Cccv$ such that
\begin{align}
\inf_{p\in\Cccv}U(p)=U(p_*)\, .
\end{align}
\end{lemma}
\begin{proof}
Let $q_n\in\Cccv$ be a minimizing sequence for
$\inf_{p\in\Cccv}U(p)$, and write
$\theta_n=\theta(q_n)$.  After passing to a subsequence,
$\theta_n\to\theta\in[0,1]$.  Helly's selection theorem gives a
further subsequence converging at every continuity point of a
non-decreasing limit $q_*$.  We take $q_*$ to be right-continuous.
For $t<\theta$,
 eventually $t<\theta_n-\eps$ for some $\eps>0$, and hence 
$q_*(t)=0$. 
For every $\delta>0$, the functions $q_n$ are concave on
$[\theta+\delta,1]$ for all sufficiently large $n$.  Their limit is
therefore concave on $(\theta,1]$.  Defining $q_*(\theta)$ by right
continuity, the restriction of $q_*$ to $[\theta,1]$ is concave and
absolutely continuous, while $q_*(1)=1$.  Thus $q_*\in\Cccv$.

We next claim that $q_n'\to q_*'$ almost everywhere.  If $t<\theta$,
then $q_n'(t)=0$ for all sufficiently large $n$.  If
$t\in(\theta,1)$, use right derivatives as representatives.  For
$0<h<\min(t-\theta,1-t)$, concavity gives, for all sufficiently
large $n$,
\begin{align}
\frac{q_n(t+h)-q_n(t)}{h}
\le q_n'(t+)
\le \frac{q_n(t)-q_n(t-h)}{h}.
\label{eq:derivative-squeeze}
\end{align}
Taking first $n\to\infty$ and then $h\downarrow0$ proves the claim
at every differentiability point of $q_*$.

Set
\begin{align*}
\widehat u_n(t)&:=u(t;q_n),\\
\widehat f_n(t)&:=L(q_n(t),q_n'(t),\widehat u_n(t),t),\\
\widehat v_n(t)&:=\frac12q_n(t)\xi(t)-\widehat u_n(t)
=\int_0^t\widehat f_n(s)\,\de s.
\end{align*}
By Eq.~\eqref{eq:u-upper-bound},
$0\le\widehat u_n(t)\le q_n(t)\xi(t)/2$.  Lemma
\ref{lemma:Sqrt-t} gives
\begin{align}
0\le\widehat f_n(s)\le
\frac{C(\xi,\alpha)}{\sqrt{s-\theta_n}}
\mathbbm{1}_{\{s>\theta_n\}}
\quad\text{for almost every }s.
\label{eq:Ln-dominating-bound}
\end{align}
Consequently, $(\widehat v_n)$ is uniformly bounded and
$1/2$-Hölder continuous.  Passing to a further subsequence, it
converges uniformly to a function $\widehat v_*$, by
the Arzel\'a-Ascoli theorem.  Define
\begin{align*}
\widehat u_*(t):=\frac12q_*(t)\xi(t)-\widehat v_*(t).
\end{align*}
At almost every $s$, the convergence of $q_n(s)$, $q_n'(s)$, and
$\widehat u_n(s)$, together with Eq.~\eqref{eq:R-omega}, gives
\begin{align*}
\widehat f_n(s)\longrightarrow
L(q_*(s),q_*'(s),\widehat u_*(s),s).
\end{align*}
The estimate \eqref{eq:Ln-dominating-bound} also implies uniform
integrability: for every measurable $E\subseteq[0,1]$,
\begin{align*}
\sup_n\int_E\widehat f_n(s)\,\de s\le 2C\sqrt{|E|}.
\end{align*}
Vitali's convergence theorem therefore allows us to pass to the
limit in the integral equation, obtaining
\begin{align*}
\widehat u_*(t)=\frac12q_*(t)\xi(t)
-\int_0^t
L(q_*(s),q_*'(s),\widehat u_*(s),s)\,\de s.
\end{align*}
By Lemma \ref{lem:U-continuous}, this is the unique solution
associated with $q_*$.  It follows that
\begin{align*}
\inf_{p\in\Cccv}U(p)
=\lim_{n\to\infty}U(q_n)
=U(q_*).
\end{align*}
Taking $p_*=q_*$ proves the claim.
\end{proof}

\begin{lemma}
\label{lem:U-continuous-in-D}
If $p\in \cuD$ is such that $U(p)=\inf_{\overline{p}\in\cuD}U(\overline{p})>0$, then $p\in \Ccont$.
\end{lemma}
\begin{proof}
Write $u(t)=u(t;p)$, $\theta=\theta(p)$, and
$M=U(p)=\inf_{\overline{p}\in\cuD}U(\overline{p})>0$.
We first show that 
\begin{align}
u(s)\ge M p(s),\qquad s\in[0,1].
\label{eq:positive-minimizer-bound}
\end{align}
The claim is immediate when $p(s)=0$. Otherwise, define
\[
q_s(t)=
\begin{cases}
p(t)/p(s),&t\le s,\\
1,&t>s.
\end{cases}
\]
Then $q_s\in\cuD$. The homogeneity argument in the proof of
Corollary~\ref{coro:scaling} applies directly to the integral equation
for controls in $\cuD$ as well, and gives
$u(s;q_s)=u(s)/p(s)$. On $[s,1]$ the control $q_s$ is constant.
Since the drift is non-increasing in $u$ and $\alpha\le1$,
Eq.~\eqref{eq:drift-0} gives
\[
\drift(t,u;1,0)\le\drift(t,0;1,0)
=\frac12\left(\sqrt{\xi'(t)}-
\sqrt{\xi'(t)/\alpha}\right)_+^2=0
\qquad (u\ge0).
\]
Thus $M\le U(q_s)\le u(s)/p(s)$, proving
Eq.~\eqref{eq:positive-minimizer-bound}.

Suppose now that $p$ has a jump at some $t_0\in(\theta,1]$, of size
$J=p(t_0)-p(t_0-)>0$. Fix $t_*\in(\theta,t_0)$; then $p(t_*)>0$.
For sufficiently small $\eps>0$, put
\[
t_{-\eps}=t_0-\eps>t_*,\qquad
t_{+\eps}=\min\{t_0+\eps,1\},\qquad
\ell_\eps=t_{+\eps}-t_{-\eps},\qquad
\Delta_\eps=p(t_{+\eps})-p(t_{-\eps}).
\]
Replace $p$ by its chord on this interval:
\[
p_\eps(t)=
\begin{cases}
p(t),&t\notin(t_{-\eps},t_{+\eps}),\\[1ex]
p(t_{-\eps})+\dfrac{\Delta_\eps}{\ell_\eps}(t-t_{-\eps}),
&t\in(t_{-\eps},t_{+\eps}).
\end{cases}
\]
Note that $p_{\eps}\in\cuD$ and $\theta(p_{\eps})=\theta(p)$.
Moreover, $\lim_{\eps\to 0}\Delta_\eps= J$, including when $t_0=1$.
Write $u_\eps(t)=u(t;p_\eps)$ and $r=\dot p$.

Recall the  bounds in Eq.~\eqref{eq:BasicLineq}.
Both solutions $u$ and $u_{\eps}$ are bounded
between zero and $\xi(1)/2$ by
Lemma \ref{lem:U-continuous}. Since $r\in L^1([0,1])$,
Cauchy--Schwarz therefore gives, for $C=C(\xi,\alpha)$,
\begin{align*}
\int_{t_{-\eps}}^{t_{+\eps}}L(p(t),r(t),u(t),t)\,\de t
&\le C\left(\ell_\eps+
\sqrt{\ell_\eps\int_{t_{-\eps}}^{t_{+\eps}}r(t)\,\de t}\right)
=o(\sqrt{\ell_\eps}).
\end{align*}
On the other hand, $p_\eps(t)\ge p(t_*)>0$ for $t\in[t_{-\eps},t_{+\eps}]$.
Since $\xi(0)>0$ (we assumed this to be the case without 
loss of generality by the remark at the beginning of Section
\ref{sec:SolOptControl})
and $\xi'(t)>0$ for $t>0$, for 
$\Lambda(p_\eps,u_\eps,t)$ of Lemma~\ref{lemma:Sqrt-t} 
we have $\Lambda(p_\eps,u_\eps,t)\ge \lambda$ for all 
$t\in[t_{-\eps},t_{+\eps}]$, where $\lambda>0$ is independent of $\eps$. Using
$\dot p_\eps=\Delta_\eps/\ell_\eps$, we obtain
\begin{align*}
		\int_{t_{-\eps}}^{t_{+\eps}}L(p_{\eps}(t),\dot p_{\eps}(t),u_{\eps}(t),t)\,\de t
&\ge\lambda\sqrt{\Delta_\eps\ell_\eps}-C\ell_\eps.
\end{align*}
Thus for
all sufficiently small $\eps$, we have
\begin{align*}
	d_{+\eps}:=\int_{t_{-\eps}}^{t_{+\eps}}L(p(t),r(t),u(t),t)\,\de t-
	\int_{t_{-\eps}}^{t_{+\eps}}L(p_{\eps}(t),\dot p_{\eps}(t),u_{\eps}(t),t)\,\de t
	<0\, .
\end{align*}
Subtracting the
integral equations \eqref{eq:ODE-integral} for $u$ and $u_{\eps}$ at $t_{+\eps}$ gives
\[
u_\eps(t_{+\eps})-u(t_{+\eps})=d_{+\eps} <0.
\]
The upper bound in Eq.~\eqref{eq:BasicLineq} also gives
$|d_{+\eps}|=O(\sqrt{\ell_\eps})$.
If $t_{+\eps}=1$, this already contradicts optimality.

Otherwise, define 
\begin{align*}
d_\eps(t):=u_\eps(t)-u(t) \mbox{ for  }t\in[t_{+\eps},1].
\end{align*}
Of course $d_\eps(t_{+\eps})=d_{+\eps}$.
Since $p_\eps=p$ on $[t_{+\eps},1]$,  it follows that
$d_\eps$
is absolutely continuous on $[t_{+\eps},1]$ 
(the term
$p\xi/2$ cancels when the integral equations are subtracted.)

Monotonicity of $L$ in $u$ gives
$0\le d_\eps(t)\le d_\eps(t_{+\eps})$. By
Eq.~\eqref{eq:positive-minimizer-bound}, $u(t)\ge M p(t_*)$ on
this interval. For small $\eps$, both solutions are consequently
bounded below by $M p(t_*)/2$. The pointwise derivative estimate
in Lemma~\ref{lemma:drift-continuous}, Point~$(6)$
 yields
\[
\dot d_\eps(t)=-A_\eps(t)d_\eps(t),\qquad
0\le A_\eps(t)\le C(1+\sqrt{r(t)})\in L^1([t_{+\eps},1]).
\]
It follows that
\[
U(p_\eps)-U(p)=d_\eps(1)
=d_\eps(t_{+\eps})
\exp\left\{-\int_{t_{+\eps}}^1 A_\eps(t)\,\de t\right\}<0,
\]
again contradicting optimality.

Hence $p$ has no jumps in $(\theta,1]$. Since a monotone function
can only have jump discontinuities and $p$ is right-continuous at
$\theta$, its restriction to $[\theta,1]$ is continuous. 
\end{proof}
\subsection{Continuity of $\dot{u}(t)$}

\begin{lemma}
	\label{lem:phidot-Lip}For any $\xi$ and $\delta,\eta>0$ there exists
	$a=a(\xi,\delta,\eta)>0$ such that the following holds. Suppose
	that $p\in\cuC_*$ and let $u$ be the solution to
	the ODE \eqref{eq:ODE-repeated} with initial condition $u(t_0) = 0$ 
	at $t_{0}=\theta$ 
	(interpreted in integral sense, cf. 
	Eq.~\eqref{eq:ODE-integral}). Define $t_l:=\max\{t:p(t)\leq1-\eta\}$.
	For any $t_{1},t_{2}\in[\theta+\eta,t_l]$, 
	\[
	|t_{1}-t_{2}|\leq a\implies|\dot{u}(t_{1})-\dot{u}(t_{2})|\leq\delta\, .
	\]
\end{lemma}

\begin{proof}
	Since $p\in\cuC_*$ is $C^{\infty}$ on $I(\eta):=[\theta+\eta,t_l]$, 	
	the unique solution of \eqref{eq:ODE-integral}
	 is in fact differentiable (because $u(t;p)-p\xi(t)/2$
	 is the integral of a continuous function) with 
\begin{align}
	\dot{u}(t) = \cF(t,u(t),p(t),\dot{p}(t)) 
\end{align}
Since by Eq.~\eqref{eq:drift-L}, $\cF(t,u(t),p(t),\dot{p}(t))\le (p\xi)'(t)/2$, we have
$u(t) \le C(\xi)$ for a suitable constant $C$ and all $t$. 
By Lemma~\ref{lemma:drift-continuous}, Point $(2)$,
 $(t,v)\mapsto \cF(t,v):=\cF(t,v,p(t),\dot{p}(t))$ is continuous 
on the compact set $I(\eta)\times [0,C(\xi)]$, and therefore uniformly continuous.
In other words, there exists a function $\Delta_0:\reals_{\ge0}\to\reals_{\ge0}$ (depending on $\xi,\eta$, since the modulus of continuity of $\dot{p}$ only
depends on $\eta$ by Lemma \ref{lem:SecondDerivative}) such that $\Delta_0(t)\to0$ as $t\to0$ and 
\begin{align}
|\dot{u}(t_1)-\dot{u}(t_2)|&= \big|\cF(t_1,u(t_1))-\cF(t_2,u(t_2))\big|\\
 &\le \Delta_0\Big(|t_1-t_2|+|u(t_1)-u(t_2)|\Big).
\end{align}
By uniform continuity, $\cF$ is bounded on  $I(\eta)\times [0,C(\xi)]$, 
and therefore $|\dot{u}(t)|$ is also bounded, whence 
$|u(t_1)-u(t_2)|\le C'(\xi,\eta)|t_1-t_2|$.
Substituting above, we see that there exists a second 
modulus of continuity 
$\Delta_1$ (depending on $\xi,\eta$ and such that 
$\lim_{s\to 0}\Delta_1(s)=0$) such that
\begin{align}
	|\dot{u}(t_1)-\dot{u}(t_2)|&\le  \Delta_1\Big(|t_1-t_2|\Big).
	\end{align}
This completes the proof.
\end{proof}

\begin{lemma}
	\label{lem:phidot-Lip-constp}
	For any $\xi$ and $\delta>0$ there exists
	$a=a(\xi,C,\delta)>0$ such that the following holds. Suppose
	that $p\equiv c$ is the constant function on $[t_l,1]$ for some 
	$c\in[0,1]$ and let $u$ be the solution to
	\eqref{eq:ODE-repeated} with $t_{0}=t_l$ 
	(for $t_l$ as defined in Lemma \ref{lem:phidot-Lip}) and $u_0\in[0,C]$. 
	For any $t_{1},t_{2}\in[f_l,1]$, 
	\[
	|t_{1}-t_{2}|\leq a \implies|\dot{u}(t_{1})-\dot{u}(t_{2})|\leq\delta\, .
	\]
\end{lemma}
\begin{proof}
	As in the proof of Lemma \ref{lem:phidot-Lip}, 
	we have that $u(t) \le C_1(C,\xi)$ for a suitable constant $C_1$ and all 
	$t\in[t_l,1]$. 
	By Lemma~\ref{lemma:drift-continuous}, Point $(2)$,
	$(t,v)\mapsto \cF_c(t,v):=\cF(t,v,c,0)$ is continuous 
	on the compact set $[t_l,1]\times [0,C_1]$, and therefore
	 uniformly continuous. Hence there exists a modulus of continuity
	 $\Delta_0$ such that:
\begin{align}
	|\dot{u}(t_1)-\dot{u}(t_2)|&\le \big|\cF_c(t_1,u(t_1))-\cF_c(t_2,u(t_2))\big|\\
	 &\le \Delta_0\Big(|t_1-t_2|+|u(t_1)-u(t_2)|\Big).
	\end{align}
 Further, by continuity $|\cF_c(t,u)|\le C_2(C,\xi)$ on the bounded set
  $[t_l,1]\times [0,C_1]$, whence $|u(t_1)-u(t_2)|\le C_2|t_1-t_2|$,
  and therefore, for a new modulus of continuity $\Delta_1$,
\begin{align}
	|\dot{u}(t_1)-\dot{u}(t_2)|&\le\Delta_1\Big(|t_1-t_2|\Big).
\end{align}
\end{proof}

\section{Equivalence of expressions for the algorithmic threshold:\\
Proof of Theorem \ref{thm:Equivalence}}
\label{app:Equivalence}

In this appendix we use the notational conventions 
introduced in Appendix \ref{app:Properties}.
We assume $\xi(0)>0$, and $t\mapsto \xi(t)$
is non-constant.

\subsection{An optimal control for $t\le q_{\sRS}$}
\label{sec:OptimalControl_qrs}

Recall the definition of `replica symmetric' energy 
from Eq.~\eqref{eq:uRS_def}
\begin{align}
u_{\sRS}(t,\alpha,\xi):=
		\frac{1}{2}
		\Big(\sqrt{\xi(t)}-\sqrt{\frac{1}{\alpha}t\xi'(t)}\Big)_+^2\, .
\end{align}
and, from Eq.~\eqref{eq:qrsDef}, that 
	\begin{align}
		q_{\sRS}:= \arg\min_{q\ge 0}
		\frac{\xi(q)\xi'(q)}{q}\, .
	\end{align}
We refer to the discussion below Eq.~\eqref{eq:qrsDef}
for the case $\xi(t)=\xi_0+\xi_1 t$ and the behavior
at $q=0$.

%
%

The next theorem establishes that this energy 
coincides with the solution of the ODE
\eqref{eq:ODE-integral} for $t\le q_{\sRS}$. 
Recall that we are theoughout assuming, without loss of
generality $\xi(0)>0$.
The assumption $\xi'(0)>0$ below is necessary 
and sufficient to ensure $q_{\sRS}>0$ 
(because otherwise $\xi(q)\xi'(q)/q$ is 
non-negative coefficients of degree at least one.)
We also point out that the statement below holds for
any $t\ge 0$, anltough we will use it only for $t\in[0,1]$.
\begin{theorem}\label{thm:MatchingSol}
	Assume $\xi$ is not affine, and $\xi(0),\xi'(0)>0$. 
	Define, for $C_0>0$ an arbitrary 
	constant:
	\begin{align}
		p(t) &= C_0\sqrt{\frac{t}{\xi(t)\xi'(t)}}\, .
	\end{align}
	Then, for any $t\in [0,q_{\sRS}]$, the solution to the 
	ODE~\eqref{eq:ulb} (equivalently, \eqref{eq:ODE-integral}) is given by
	\begin{align}\label{eq:usol}
	u(t;p) &= p(t) \cdot u_{\sRS}(t,\alpha,\xi)\, .
	\end{align}
	\end{theorem}
	\begin{proof}
	Recall that the definition of 
	 $\drift(t,u) = \drift(t,u;p,\dot{p})$ 
	in Eq.~\eqref{eq:F}. Since the given $p$ is differentiable,
	integral equation~\eqref{eq:ODE-integral} can be written in the classical form
	\begin{align}
	\dot{u} = \drift(t,u;p,\dot{p}) \, .
	\end{align}
    By Corollary \ref{coro:scaling}, $u(t;c\cdot p)=c u(t;p)$ for all $t\ge 0$.
	It is therefore sufficient to prove the claim for $C_0=1$.
	We will also write
	\begin{align}
		\hcL(\omega;p,\dot{p},u,t)&:=
		\cR(\omega;u,t)-\frac{b}{2}\\
		&=\frac{1}{\sqrt{\alpha}}\sqrt{2u(p\xi')'+p\xi'\omega^2}
		-\sqrt{(p\xi)'}\,\omega+
		\frac{1}{2}\omega^2\, ,
	\end{align}
	so that $\drift(t,u;p,\dot{p})= -\min_{\omega\ge 0}
	 \hcL(\omega;p,\dot{p},u,t)$, see Eqs.~\eqref{eq:L-omega}
	 and \eqref{eq:drift-L}.
	We will also use the notation 
	$\hcL(\omega,t)=\hcL(\omega;p,\dot{p},u,t)$.
	
	Recall the notations of Eq.~\eqref{eq:ab-def}
	(to which we add the definition of $d(t)$):
	\begin{align}
	a(t) \coloneqq (p\xi')(t), \quad
	d(t) \coloneqq (p\xi)(t), \quad
	b(t) \coloneqq d'(t)=(p\xi)'(t).
	\end{align}
	From $p(t)=\sqrt{t/(\xi(t)\xi'(t))}$ we have the identity
	\begin{equation}
	a(t)d(t)=p^2\xi\xi'=t,
	\end{equation}
	hence
	\begin{equation}
	d(t)=\frac{t}{a(t)}.
	\end{equation}
	Differentiating yields
	\begin{equation}
	b(t)=\Bigl(\frac{t}{a(t)}\Bigr)'
	=\frac{1}{a(t)}-\frac{t a'(t)}{a(t)^2}.
	\label{eq:bidentity}
	\end{equation}

	In what follows we can assume 
	$t<q_0$ (with $q_0$ defined as per 
	Eq.~\eqref{eq:q0Def}, since for  $t\in (q_0,q_{\sRS})$,
	 we have 
	$u_{\sRS}(t,\alpha,\xi)=0$ and $\drift(t,0;p,\dot{p})=0$.)

	\paragraph{Expression for $u(t)$.}
	We want to check that \eqref{eq:usol} is a solution, hence,
	for $t<q_0$,
	\begin{equation}
	u(t)=\frac12 p(t)\Bigl(\sqrt{\xi(t)}-\sqrt{\tfrac{1}{\alpha}t\xi'(t)}\Bigr)^2 .
	\end{equation}
	Since $a(t)=p(t)\xi'(t)=\sqrt{t\xi'(t)/\xi(t)}$, we obtain
	\[
	\sqrt{\tfrac{1}{\alpha}t\xi'}=\sqrt{\xi}\,\frac{a}{\sqrt{\alpha}},
	\]
	and therefore
	\begin{equation}
	u(t)=\frac12 (p\xi)\Bigl(1-\frac{a}{\sqrt{\alpha}}\Bigr)^2
	=\frac12\,\frac{t}{a}\Bigl(1-\frac{a}{\sqrt{\alpha}}\Bigr)^2.
	\end{equation}
	Expanding,
	\begin{equation}
	u(t)=\frac{t}{2a}-\frac{t}{\sqrt{\alpha}}+\frac{t a}{2\alpha}.
	\label{eq:uexpanded}
	\end{equation}
	
	\paragraph{Derivative of $u(t)$.}
	
	Differentiating \eqref{eq:uexpanded},
	\begin{align}
	u'(t)
	&=\frac12\Bigl(\frac{1}{a}-\frac{t a'}{a^2}\Bigr)
	-\frac{1}{\sqrt{\alpha}}
	+\frac{a+t a'}{2\alpha}.
	\end{align}
	Using \eqref{eq:bidentity},
	\begin{equation}
	u'(t)=\frac12 b(t)-\frac{1}{\sqrt{\alpha}}+\frac{a(t)+t a'(t)}{2\alpha}.
	\label{eq:uprime}
	\end{equation}
	
	\paragraph{Simplified form of $\mathcal L$.}
	
	The function $\hcL(\omega,t)$ reads
	\begin{equation}
	\hcL(\omega,t)
	=\frac{1}{\sqrt{\alpha}}\sqrt{2u a'+a\omega^2}
	+\frac12\omega^2-\omega\sqrt b .
	\label{eq:Lsimplified}
	\end{equation}
	We will denote by $\omega_\star(t)$ the unique minimizer
	of $\hcL(\,\cdot\,,t)$ so that $-\drift(t,u) = \hcL(\omega_\star(t),t)$.
	
	\paragraph{Optimality condition.}
	
	From \eqref{eq:Lsimplified},
	\begin{equation}
	0=\partial_\omega\hcL(\omega^\star,t)
	=\frac{1}{\sqrt{\alpha}}\frac{a\omega^\star}{\sqrt{2u a'+a\omega^{\star2}}}
	+\omega^\star-\sqrt b .
	\end{equation}
	Hence
	\begin{equation}
	\sqrt b-\omega^\star
	=\frac{1}{\sqrt{\alpha}}\frac{a\omega^\star}{\sqrt{2u a'+a\omega^{\star2}}}.
	\label{eq:FOC-star}
	\end{equation}
	
	\paragraph{Evaluation of $L(t)$.}
	
	At the minimizer,
	\begin{align}
	-\drift(t,u) =
	&=\frac{1}{\sqrt{\alpha}}\sqrt{2u a'+a\omega^{\star2}}
	+\frac12\omega^{\star2}-\omega^\star\sqrt b \\
	&=\frac{2u a'}{\sqrt{\alpha}\sqrt{2u a'+a\omega^{\star2}}}
	-\frac12\omega^{\star2}.
	\label{eq:Lintermediate}
	\end{align}
	
	Next observe that from \eqref{eq:uexpanded},
	\begin{equation}
	2u a'
	=\frac{t a'}{a}\Bigl(1-\frac{a}{\sqrt{\alpha}}\Bigr)^2.
	\end{equation}
	Using \eqref{eq:bidentity}, $\frac{t a'}{a}=1-a b$, hence
	\begin{equation}
	2u a'=(1-a b)\Bigl(1-\frac{a}{\sqrt{\alpha}}\Bigr)^2.
	\label{eq:keycollapse}
	\end{equation}
	
	\paragraph{More details: removing $\omega^\star\sqrt b$ and isolating the square root.}
	
	Let
	\begin{equation}
	S \coloneqq \sqrt{2u(t)a'(t)+a(t)\,\omega^{\star}(t)^2}\, .
	\end{equation}
	Then \eqref{eq:FOC-star} is equivalently
	\begin{equation}
	\sqrt b
	=\omega^\star + \frac{1}{\sqrt{\alpha}}\frac{a\,\omega^\star}{S}
	=\omega^\star\Bigl(1+\frac{a}{\sqrt{\alpha}\,S}\Bigr).
	\label{eq:FOC_S_correct}
	\end{equation}
	In particular,
	\begin{equation}
	\omega^\star \sqrt b
	=\omega^{\star 2} + \frac{1}{\sqrt{\alpha}}\frac{a\,\omega^{\star 2}}{S}.
	\label{eq:omegasqrtb}
	\end{equation}
	Substituting \eqref{eq:omegasqrtb} into \eqref{eq:Lsimplified} evaluated at $\omega^\star$ gives
	\begin{align}
	-\drift(t,u) =
	&=\frac{1}{\sqrt{\alpha}}S+\frac12\omega^{\star 2}-\omega^\star\sqrt b \\
	&=\frac{1}{\sqrt{\alpha}}\Bigl(S-\frac{a\,\omega^{\star 2}}{S}\Bigr)-\frac12\omega^{\star 2}.
	\end{align}
	Using $S^2=2u a' + a\omega^{\star 2}$, we have
	\begin{equation}
	S-\frac{a\,\omega^{\star 2}}{S}
	=\frac{S^2-a\omega^{\star 2}}{S}
	=\frac{2u a'}{S},
	\end{equation}
	hence we recover \eqref{eq:Lintermediate}:
	\begin{equation}
	-\drift(t,u)=\frac{2u a'}{\sqrt{\alpha}\,S}-\frac12\omega^{\star 2}.
	\label{eq:Lintermediate_rederived}
	\end{equation}
	
	\paragraph{Reducing the minimization to a scalar equation.}
	
	Introduce the dimensionless ratio
	\begin{equation}
	r \coloneqq \frac{a}{\sqrt{\alpha}\,S}\qquad\Bigl(\text{so that }S=\frac{a}{\sqrt{\alpha}\,r}\Bigr).
	\end{equation}
	Then \eqref{eq:FOC_S_correct} yields
	\begin{equation}
	\omega^\star = \frac{\sqrt b}{1+r},
	\qquad\text{and hence}\qquad
	\omega^{\star 2}=\frac{b}{(1+r)^2}.
	\label{eq:omega_in_r}
	\end{equation}
	Plugging \eqref{eq:omega_in_r} into $S^2=2u a'+a\omega^{\star 2}$ gives a scalar equation for $r$:
	\begin{equation}
	\frac{a^2}{\alpha r^2}
	=2u a' + a\,\frac{b}{(1+r)^2}.
	\label{eq:r_equation}
	\end{equation}
	Moreover \eqref{eq:Lintermediate_rederived} becomes
	\begin{equation}
	-\drift(t,u)=\frac{2u a'}{a}\,r-\frac12\,\frac{b}{(1+r)^2}.
	\label{eq:L_in_r}
	\end{equation}
	
	\paragraph{Solving for $\drift(t,u)$ and final simplification.}
	Substituting the identity \eqref{eq:keycollapse},
	into \eqref{eq:r_equation} yields
	\begin{equation}
	0=-\frac{a^2}{\alpha r^2}
	+(1-a b)\Bigl(1-\frac{a}{\sqrt{\alpha}}\Bigr)^2
	+\frac{a b}{(1+r)^2} =: \frac{1}{r^2(1+r)^2} \, P(r).
	\label{eq:r_expanded}
	\end{equation}
	Here $P(r)$ is a degree four polynomial which factorizes as follows:
	\begin{align}
	P(r) & = (1+r)^2\left[r^2\Bigl(1-\frac{a}{\sqrt{\alpha}}\Bigr)^2-
	\frac{a^2}{\alpha}\right]
	+abr^2\left[1-(1+r)^2\Bigl(1-\frac{a}{\sqrt{\alpha}}\Bigr)^2\right]\\
	& = \Bigl(1-\frac{a}{\sqrt{\alpha}}\Bigr)^2(r-r_0)P_0(r)\, ,
	\end{align}
	where
	\begin{align}
		r_0&:=\frac{a}{\sqrt{\alpha}-a}\, ,\\
		P_0(r) & := (1+r^2)(r+r_0) +abr^2(r+r_0+2)\, .
	\end{align}
	Note that $a(t) = p(t)\xi'(t) = \sqrt{t\xi'(t)/\xi(t)}
	<\sqrt{\alpha}$ for $t<q_0$ fixed.
	It follows that $r_0>0$ is the only non-negative root of $P(r)$ 
	hence we have
	\begin{equation}
	r=\frac{a}{\sqrt{\alpha}-a}.
	\label{eq:r_solution}
	\end{equation}
	Substituting \eqref{eq:r_solution} into \eqref{eq:L_in_r} gives
	\begin{align}
	-\drift(t,u)
	&=\frac{2u a'}{a}\,\frac{a}{\sqrt{\alpha}-a}
	-\frac12\,\frac{b(\sqrt{\alpha}-a)^2}{\alpha} \\
	& = \frac{1}{\sqrt{\alpha}-1}(1-ab) \Big(1-\frac{a}{\sqrt{\alpha}}\Big)^2
	-\frac12 b  \Big(1-\frac{a}{\sqrt{\alpha}}\Big)^2\\
	& = -\frac{1}{2}b+\frac{1}{\sqrt{\alpha}}
	-\frac{a}{2\alpha}\big(2-ab\big)\\
	&=\frac{1}{\sqrt{\alpha}}-\frac12 b-\frac{a+t a'}{2\alpha},
	\label{eq:LastStepL}
	\end{align}
	where in the first step we again used \eqref{eq:keycollapse} 
	and in the last we used
	$\frac{t a'}{a}=1-a b$.

	\paragraph{Conclusion.}
	
	Comparing \eqref{eq:uprime} and \eqref{eq:LastStepL} we obtain the 
	identity $u'(t)=\drift(t,u)$, thus completing the proof.
	
	\end{proof}

	\begin{lemma}\label{thm:RS-lower-bound}
		For any $p\in\Ccont$ and $t\in[0,1]$,
		\begin{align}
			u(t;p) \ge p(t)u_{\sRS}(t,\alpha,\xi)\,.
		\end{align}
		\end{lemma}
		\begin{proof}
		By Lemma~\ref{lem:U-continuous}, $u(t;p)\ge0$, so the claim is
		immediate when $u_{\sRS}(t,\alpha,\xi)=0$. We therefore work on
		$0<t<q_0$. If $\theta(p)\ge q_0$, the claim also follows
		immediately, since $p=0$ before $\theta(p)$. Hence assume
		$\theta(p)<q_0$, and put
		\begin{align*}
		\beta(t)&:=\sqrt{\frac{t\xid}{\alpha \xi}}\in(0,1),\\
		u_{\sRS}(t)&=u_{\sRS}(t,\alpha,\xi)
		=\frac{\xi}{2}(1-\beta)^2,\\
		K(t)&:=\frac{\xi}{2}-u_{\sRS}(t)
		=\xi \left(\beta-\frac{\beta^2}{2}\right)\ge0.
		\end{align*}
		We claim that, at almost every $t\in(\theta(p),q_0\wedge1)$,
		\begin{equation}\label{eq:RS-barrier}
		L\big(p(t),\dot p(t),p(t)u_{\sRS}(t),t\big)
		\le \dot p(t)K(t)+p(t)K'(t).
		\end{equation}
		Indeed, $p(t)>0$ and $\xid >0$ there. In the variational
		formula~\eqref{eq:L-s}, choose
		\[
		s=\frac{\beta}{p \xid (1-\beta)},
		\qquad r:=p\xid s=\frac{\beta}{1-\beta}.
		\]
		Since $(p\xi')'=\dot p\,\xid+p\xidd$ and
		$(p\xi)'=\dot p\,\xi+p\xid$, this gives
		\begin{align*}
		L(p,\dot p,pu_{\sRS}(t),t)
		&\le \dot p\left(ru_{\sRS}(t)+\frac{r\xi}{2(1+r)}\right)
		+p\left(\frac{u_{\sRS}(t)\xidd r}{\xid}+\frac{\xid}{2\alpha r}
		+\frac{\xid r}{2(1+r)}\right)\\
		&=\dot p K+pK'.
		\end{align*}
		For the last equality, use $t\xid =\alpha \xi \beta^2$ and
		$\beta'/\beta=(1/t+\xidd/\xid -\xid/\xi)/2$. This proves
		\eqref{eq:RS-barrier}.

		Set $g(t)=p(t)\xi/2$, $v(t)=g(t)-u(t;p)$, and
		$h(t)=p(t)K(t)$. Then $v$ is absolutely continuous by
		Lemma~\ref{lem:U-continuous}, while $h$ is continuous and of
		bounded variation on $[\theta(p),q_0\wedge1]$. At the initial
		time, $v(\theta(p))=0\le h(\theta(p))$.
		The absolutely continuous part of $\de h$ is
		$(\dot p K+pK')\,\de t$, and its singular part is
		$K\,\de p_{\rm s}\ge0$. On the set $v>h$ we have
		$u=g-v<pu_{\sRS}(t)$. Since $L$ is non-decreasing in its energy
		argument, Eq.~\eqref{eq:RS-barrier} implies
		\[
		v'(t)=L(p,\dot p,u,t)
		\le L(p,\dot p,pu_{\sRS}(t),t)
		\le \dot p K+pK'
		\]
		for almost every such $t$. Thus the positive part of $v-h$
		cannot increase: its absolutely continuous derivative is
		non-positive wherever it is positive, and its singular
		variation is $-K\,\de p_{\rm s}\le0$.
		As it starts at zero, $v\le h$, or equivalently
		$u(t;p)\ge p(t)u_{\sRS}(t)$, for $\theta(p)\le t<q_0\wedge1$.
		For $t<\theta(p)$ both sides vanish, and continuity gives
		the claim at $q_0\wedge1$.
		\end{proof}

\subsection{Concavity of the drift and optimal value}
\label{sec:ConcDrift}

Recall the definition of $L(p,\dot p,u,t)$ from 
Eq.~\eqref{eq:L-omega}, which we rewrite, setting 
$p=p(t)$ and $\dot p(t)=r$, as
	\[
	L(p,r,u,t)
	=
	\min_{\omega\ge 0}
	\left\{
	\frac{1}{\sqrt\alpha}
	\sqrt{
	2u\,(p\xi''+r\xi') + p\xi'\,\omega^2
	}
	+
	\frac12\left(
	\omega-\sqrt{p\xi'+r\xi}
	\right)^2
	\right\}.
	\]
 Using the same identities as in 
 Appendix \ref{app:UsefulIdentities}, we have
	\begin{align}
		L(p,r,u,t)
		=
		\inf_{s>0} 
		\left\{
		\frac{1}{2\alpha s}
		+
		u(r\xidot+p\xiddot)s
		+
		\frac{(r\xi+p\xidot)p\xidot s}{2(1+p\xidot s)}
		\right\}.\label{eq:L-inf-s}
		\end{align}

We also recall the integral equation
\eqref{eq:ODE-integral}. Extending $\dot p$ and $L$ by zero before
$\theta(p)$ and using the identity \eqref{eq:drift-L}, it can be
written for every $t\in[0,1]$ as
\begin{align}\label{eq:ODE-integral-2}
	u(t;p) = \frac{1}{2}p(t)\xi(t)-\int_0^t 
	L(p(s),\dot{p}(s),u(s;p),s) \, \de s\, .
	\end{align}

	\begin{lemma}\label{lemma:U-strict-inequality-gen}
		Assume that $\xi$ is non-constant and
		$U_*(\Ccont) >0$.
		Let $p\in\Ccont$ and 
		$t_0\in (\theta(p),1)$ be such that $u(1;p)>0$, 
		and there exists an open interval $I\subset (0,1)$, $t_0\in I$,
		 and  such that $p$ is convex on $I$ and $p(t)\ge 
		 p(t_0)+c_1(t-t_0)+c_2(t-t_0)^2$ for some $c_1\in\reals$,
		 $c_2>0$ and all $t\in I$.

		Then
		\begin{align}
		u(1;p)> U_*(\Ccont) = \inf_{\ovp\in\Ccont} u(1;\ovp)\, .
		\end{align}
		\end{lemma}
		\begin{proof}
Suppose, for a contradiction, that
$u(1;p)=U_*(\Ccont)>0$. Since $\xi$ is non-constant and has
non-negative coefficients, $\xi'(t)>0$ for $t>0$.
Write $u(t)=u(t;p)$.

First, $u(s)>0$ whenever $p(s)>0$. Indeed, if $u(s)=0$ and
$p(s)>0$, define
\[
\widehat p(t)=
\begin{cases}
p(t)/p(s),&t\le s,\\
1,&t>s.
\end{cases}
\]
Then $\widehat p\in\Ccont$ and, by
Corollary~\ref{coro:scaling}, $u(s;\widehat p)=0$.
For the constant control $\widehat p=1$, the variational formula
\eqref{eq:L-omega} gives
\[
\drift(t,0;1,0)
=\frac12\left(\sqrt{\xi'(t)}-
\sqrt{\xi'(t)/\alpha}\right)_+^2=0,
\]
since $\alpha\in(0,1]$. Uniqueness therefore gives
$u(1;\widehat p)=0$. We conclude that $u(s)>0$ as claimed.

We next choose a closed interval $J=[a,b]\subset I$ on which
$0<p<1$ and
\[
0<m\le k(t):=\frac{\dot p(t)}{p(t)}\le M<\infty
\qquad\text{for almost every }t\in J.
\]
Such an interval exists: for small $h>0$, monotonicity and the
quadratic lower bound imply
$0\ge p(t_0-h)-p(t_0)\ge-c_1h+c_2h^2$, so $c_1>0$.
The same lower bound on the right of $t_0$ implies
$p'_+(t_0)\ge c_1$. Convexity then bounds $\dot p$ away from zero
to the right of $t_0$, and bounds it above on compact subintervals
of $I$. In particular, $p$ is Lipschitz on $J$, so it has no singular
variation there. Equation~\eqref{eq:ODE-integral} then implies that
$u$ is absolutely continuous on $J$.

On $J$, set $e=u/p$, which is also absolutely continuous.
Note that, for any $c>0$,
$L(cp,cr,cu,t)=cL(p,r,u,t)$.
This follows from Eq.~\eqref{eq:L-omega} by rescaling 
the minimization variable $\omega=\sqrt c\,\widetilde\omega$.
In particular, $L(p,\dot p,u,t)/p=L(1,k,e,t)$ on $J$.
Using this identity and Eq.~\eqref{eq:L-inf-s}, with $z=p\xi's$, gives
\begin{align}
\dot e&=G(t,e,k):=k\left(\frac\xi2-e\right)+\frac{\xi'}2
-\ell(t,e,k),\label{eq:ODE_for_e}\\
\ell(t,e,k)&=\min_{z>0}\left\{
\frac{\xi'}{2\alpha z}
+e\left(k+\frac{\xi''}{\xi'}\right)z
+\frac{(k\xi+\xi')z}{2(1+z)}\right\}.
\end{align}
Here and below the derivatives of $\xi$ are evaluated at $t$.
For $e,k>0$, the minimizer $z$ is unique, interior, and smooth
in $(t,e,k)$. Indeed, its equation is
\[
\frac{\xi'}\alpha
=2e\left(k+\frac{\xi''}{\xi'}\right)z^2
+(k\xi+\xi')\frac{z^2}{(1+z)^2},
\]
whose right-hand side is strictly increasing from zero to infinity.

For $h\in C_c^\infty((a,b))$, consider the perturbation
\[
p_\eps(t)=p(t)\exp\left\{-\eps\int_t^1h(s)\,\de s\right\},
\]
where $h$ is extended by zero outside $J$.
Note that $p_{\eps}\in\Ccont$ for $\eps$
sufficiently small. Indeed,  for either sign of
$\eps$ sufficiently small, its logarithmic derivative on
 $J$ is, allmost everywhere, 
$k+\eps h\ge0$. Write $e_\eps=u(\,\cdot\,;p_\eps)/p_\eps$ on $J$.
Fot $t<a$, using Corollary \ref{coro:scaling}, 
both $p$ and $u$ are multiplied by
the same constant, so $e_\eps(a)=e(a)$.
For  $t\ge b$, we have $p_\eps(t)=p(t)$.

We next prove that $u(t;p_\eps)<u(t)$ for all 
$t\in (b,1]$ as well. Continuity and strict
positivity of $u$ give $m_0:=\min_{[a,1]}u>0$. Since $p_\eps=p$
on $[b,1]$, subtracting the two integral equations shows that
$d_\eps:=u(\,\cdot\,;p_\eps)-u$ is absolutely continuous there,
although the individual solutions need not be. Almost everywhere,
\[
\dot d_\eps(t)=-\big\{L(p(t),\dot p(t),u(t)+d_\eps(t),t)
-L(p(t),\dot p(t),u(t),t)\big\}.
\]
Since  $L$ is non-decreasing in $u$, the right-had side is positive
whence
$0\ge d_\eps(t)\ge  d_\eps(b)$. Continuous dependence 
of  $L(p,r,u,t)$  with respect to its arguments
for $t\in J$,
implies $d_\eps(b)\to0$, so for small $\eps$ both solutions
are at least $m_0/2$ on $[b,1]$. 
Define $A_{\eps}(t)$ by
\begin{align}
L(p(t),\dot p(t),u(t)+d_\eps(t),t)
-L(p(t),\dot p(t),u(t),t) = A_\eps d_\eps\, .
\end{align}
Lemma~\ref{lemma:drift-continuous},
Point~$(6)$, and the mean value formula now give
\[
\dot d_\eps=-A_\eps d_\eps,\qquad
0\le A_\eps(t)\le C\big(1+\sqrt{\dot p(t)}\big)\in L^1([b,1]).
\]
Consequently, for all $t\in (b,1]$,
\[
d_\eps(t)=d_\eps(b)
\exp\left\{-\int_b^t A_\eps(s)\,\de s\right\},
\]
so $d_\eps(b)<0$ 
implies $d_\eps(t)<0$ for all $t\in (b,1]$, as claimed.

Optimality of $p$  and $p_{\eps}(b)=p(b)$ 
thus implies stationarity of $e_\eps(b)$.
Differentiating the scalar ODE \eqref{eq:ODE_for_e} 
gives 
\begin{align}
	\frac{\de\phantom{t}}{\de t}\partial_{\eps}e_\eps(t)
	= G_e(t,e_{\eps}(t),k_{\eps}(t))\partial_{\eps}e_\eps(t)
+G_k(t,e_{\eps}(t),k_{\eps}(t))h(t)\, ,\;\;\;\;
\partial_{\eps}e_\eps(a)=0\, .
\end{align}
where $k_{\eps}(t) = \pd_{\eps}(t)/p_{\eps}(t)$,
and $G_e$ and $G_k$ are the partial
derivatives of $G$ with respect to $e$ and $k$, respectively.
Setting $\eps=0$ and integrating this equation, we get
for $t\ge a$,
\begin{align}
\left.\partial_{\eps}e_\eps(t)\right|_{\eps=0}
=\int_a^t
\exp\left\{\int_{t'}^bG_e(s,e(s),k(s))\,\de s\right\}
G_k(t',e(t'),k(t'))h(t')\,\de t'.
\end{align}
Setting $\left.\partial_{\eps}e_\eps(b)\right|_{\eps=0}=0$
yields
\begin{align}
\int_a^b
\exp\left\{\int_s^bG_e(t,e(t),k(t))\,\de t\right\}
G_k(s,e(s),k(s))h(s)\,\de s=0.
\end{align}
Since $h$ is arbitrary, this implies 
$G_k=0$ almost everywhere on $J$.
The envelope theorem yields
\[
0=G_k=\frac{\xi}{2(1+z)}-e(1+z),
\qquad\text{hence}\qquad
 e=\frac{\xi}{2(1+z)^2}.
\]

Introduce
\[
A=\frac{\xi'}\xi,\qquad B=\frac{\xi''}{\xi'},\qquad
\rho=\sqrt{\frac{2e}\xi},\qquad
w=\alpha(1-\rho)^2.
\]
Then $0<\rho<1$. Combining the stationarity equation with the
minimizing equation for $z$ gives
\[
k=\frac12\left(\frac A w-A-B\right).
\]
This expresses $k$ as a continuous function of $e,t$; the ODE
then implies that $e$ is $C^1$, and subsequently that $k$ is
$C^1$ and $p$ is $C^2$ in the interior of $J$.
Substituting the stationarity equation into $\dot e=G$ gives
\[
\dot\rho=\frac{A-B}{2}(1-\rho)
-\frac{A}{2\alpha(1-\rho)},
\qquad
\dot w=(B-A)w+A.
\]
Consequently, writing $Q=A/w$, we have
$\dot Q=-Q^2$ and $Q=2k+A+B$. Since
\[
\dot A=AB-A^2,\qquad
\dot B=\frac{\xi'''}{\xi'}-B^2,
\]
we obtain
\[
\frac{\ddot p}{p}=\dot k+k^2
=-k^2-2k(A+B)-\frac32AB-\frac12\frac{\xi'''}{\xi'}<0.
\]
The strict inequality uses $k>0$, $A>0$, and the non-negative
coefficients of $\xi$. This contradicts convexity of $p$ on $J$,
and proves the lemma.
\end{proof}

\begin{corollary}\label{coro:Achieved_on_CCV}
	Recall that $U_*(S) := \inf_{p\in S} u(1;p)$, and assume $U_*(\cuD) > 0$. 
	Then there exists $p_*\in\Cccv$ such that
\begin{align}
	U_*(\cuD) = U_*(\Cccv) =  u(1;p_*)\, .
\end{align}
\end{corollary}
\begin{proof}
By Lemma \ref{lemma:U-achieved-D}, and Lemma \ref{lem:U-continuous-in-D},
there exists $p_*\in \Ccont$ such that $U_*(\cuD)=U_*(\Ccont) = u(1;p_*)$.
Let $\theta=\theta(p_*)$ and assume by contradiction that $p_*$ is not concave on 
$(\theta,1]$. Then we claim that there exists $t_0\in (\theta,1)$ such
that $p_*(t)\ge 
p_*(t_0)+c_1(t-t_0)+c_2(t-t_0)^2$ for some $c_1,c_2\in\reals$,
$c_2>0$, and all $t$ in an open interval containing $t_0$. 
However, by Lemma \ref{lemma:U-strict-inequality-gen}, 
$u(1;p_*)> U_*(\Ccont)$, thus implying a contradiction. Therefore, $p_*$ is concave on $(\theta,1]$.

We are left with proving the claim $p_*(t)\ge 
p_*(t_0)+c_1(t-t_0)+c_2(t-t_0)^2$. If $p_*\in \Ccont$ is not concave on $(\theta,1]$,
then there exists $\theta<t_1<t_2<1$ such that 
$p_*(s)<B(s)-\eps(s-t_1)(t_2-s)$ for some $s\in (t_1,t_2)$ and $\eps>0$ small enough, 
where $B(t)$ is the affine function
interpolating $B(t_1) = p_*(t_1)$ and $B(t_2) = p_*(t_2)$. Define $g(t) = p_*(t)-B(t)+\eps(t-t_1)(t_2-t)$,
and let $t_0\in (t_1,t_2)$ be a minimizer of $g(t)$ over $[t_1,t_2]$.
Then the claimed bound holds on $(t_1,t_2)$ by the inequality $g(t)\ge g(t_0)$.
\end{proof}

\begin{lemma}
\label{lemma:Monotonicity-Ualg}
Assume $\xi(0),\xi'(0)>0$ and define
		$p_{*}:[0,1]\to [0,1]$ by
	\begin{align}
		p_{*}(t) := 
		\begin{cases}
			C_0(\xi)\sqrt{\frac{t}{\xi(t)\xi'(t)}} & \mbox{ if } t\in [0,q_{\sRS}]\, ,\\
			1 & \mbox{ if } t\in (q_{\sRS},1]\, ,
		\end{cases}
		\label{eq:p-star}
	\end{align}
	where $q_{\sRS}$ is defined in Eq.~\eqref{eq:qrsDef}, 
	and, for $\uq_{\sRS} :=\min(q_{\sRS},1)$,
	\begin{align}
	C_0(\xi) = \min_{t\in [0,1]}\sqrt{\frac{\xi(t)\xi'(t)}{t}} 
	= \sqrt{\frac{\xi(\uq_{\sRS})\xi'(\uq_{\sRS})}{\uq_{\sRS}}} \, .
	\end{align}
	Then $p_{*}\in\Cccv$ and $u(1;p_*)= u_{\salg}(\alpha,\xi)$.
	Further $\alpha\mapsto u_{\salg}(\alpha,\xi)$ is monotone non-decreasing.
\end{lemma}
\begin{proof}
	In order to prove that $p_*\in\ Cccv$,
	it is sufficient to show that $t\mapsto \sqrt{t/g(t)}$, $g(t):=
	\xi(t)\xi'(t)$ is concave on $[0,q_{\sRS}]$. To see this, note that
    $g$ is a polynomial with non-negative coeffigient and hence $g(t)$, $g'(t)$, $g''(t)\ge 0$
	for all $t\ge 0$. Further
	\begin{align}
	tg'(t) -g(t) =  -t^2\frac{\de \phantom{t}}{\de t}\frac{\xi(t)\xi'(t)}{t} \le 0\,\;\;
	\forall t\in [0,q_{\sRS}]\, .
	\end{align} 
    Therefore, for $t\in [0,q_{\sRS}]$,
	\begin{align}
	\frac{\de^2 \phantom{t}}{\de t^2}\sqrt{\frac{t}{g(t)}}
	= \frac{\bigl(tg'(t)-g(t)\bigr)\bigl(3tg'(t)+g(t)\bigr)
	-2t^2g(t)g''(t)}{4t^{3/2}g(t)^{5/2}} \le 0\, .
	\end{align}

	The identity $u(1;p_*)= u_{\salg}(\alpha,\xi)$ follows from Theorem \ref{thm:MatchingSol},
		and comparing the ODE \eqref{eq:ODE-integral} with the definition of $u_{\salg}(\alpha,\xi)$.
    
	Monotonicity follows simply from Corollary \ref{coro:U-monotone}, noting that $p_*$ 
	does not depend on $\alpha$. 
\end{proof}

	\begin{theorem}
		We have
	\begin{align}
		U_*(\cuD,\alpha)  = U_*(\Cccv,\alpha) =  u_{\salg}(\alpha,\xi)\, .\label{eq:IdentityUU}
	\end{align}
	\end{theorem}
	
	\begin{proof}
		Clearly $U_*(\cuD,\alpha)  \le  U_*(\Cccv,\alpha) \le  u(1;p_*)=u_{\salg}(\alpha,\xi)$ 
		for  $p_*\in\Cccv$ defined in Lemma \ref{lemma:Monotonicity-Ualg}.
		It is therefore sufficient to prove that $U_*(\cuD,\alpha) \ge  u(1;p_*)$.
		Since both $\alpha\mapsto U_*(\cuD,\alpha)$ and $\alpha\mapsto u(1;p_*)$ are monotone 
		non-decreasing,
		$\alpha \mapsto U_*(\cuD,\alpha)$ is continuous by 
		(by Corollary \ref{coro:U-monotone} and 
		Lemma~\ref{lemma:Monotonicity-Ualg},
		Lemma \ref{eq:ContInalpha}),
		and $U_*(\cuD,\alpha)>0$ for all $\alpha$
large enough, it is sufficient to prove the claim $U_*(\cuD,\alpha) \ge  u(1;p_*)$ 
(and hence $U_*(\cuD,\alpha) = u(1;p_*)$) for 
$U_*(\cuD,\alpha)>0$.

	From here on, we fix $\alpha$ andassume 
		 $U_*(\cuD)=U_*(\cuD,\alpha)>0$.
		Assume by contradiction that $U_*(\cuD,\alpha) < u(1;p_*)$.
		Let $p\in\Cccv$ be such that $U_*(\cuD) = U_*(\Cccv) =  u(1;p)<u(1;p_*)$ 
		(existence follows from Corollary \ref{coro:Achieved_on_CCV}),
		and set
		$t_*(p):=\inf\{t\in [0,1]: p(t)=1\}$.
		We will distinguish the two cases $t_*(p)<q_{\sRS}$ and $t_*(p)\ge q_{\sRS}$.

		\paragraph{Case $t_*(p)\le q_{\sRS}$.}
		At time $t=q_{\sRS}$, we have, by Theorem \ref{thm:MatchingSol},
		\begin{align}
		u(q_{\sRS};p_*) = p_*(q_{\sRS})\, u_{\sRS}(q_{\sRS};\xi)
		=u_{\sRS}(q_{\sRS};\xi), .
		\end{align}
		(We omit the argument $\alpha$ from $u_{\sRS}$, because it is fixed throughout the 
		rest of the proof.)
		On the other hand, by Lemma \ref{thm:RS-lower-bound},
		for any $p$
		\begin{align}
			u(q_{\sRS};p) \ge u_{\sRS}(q_{\sRS};\xi)\, .
		\end{align}
		In particular $u(q_{\sRS};p)\ge u(q_{\sRS};p_*)$.
		However, for $t\in [q_{\sRS},1]$, we have 
		$p(t)=p_*(t)=1$, and hence $u(t;p)$, $u(t;p_*)$ obey the
		integral equation~\eqref{eq:ODE-integral} with the different initial conditions.
		By existence and uniqueness of the solution, 
		we have $u(t;p)\ge u(t;p_*)$ for all $t\in[q_{\sRS},1]$,
		and in particular $u(1;p_*)\le u(1;p)=U_*(\Cccv)$, in contradiction with our assumption.

		\paragraph{Case $t_*(p)>q_{\sRS}$ and $p(q_{\sRS})>0$.}
		Let $b:=p(q_{\sRS})\in(0,1)$
		and define
		\begin{align}
			p_{+}(t) = 
			\begin{cases}
				b\cdot p_*(t) & \mbox{ if } t\in [0,q_{\sRS}]\, ,\\
				p(t)& \mbox{ if } t\in (q_{\sRS},1]\, ,
			\end{cases}
		\end{align}
		By Lemma \ref{thm:RS-lower-bound},
		$u(q_{\sRS};p)\ge b\cdot u_{\sRS}(q_{\sRS};\xi) = u(q_{\sRS};p_+)$
		and hence by existence and uniqueness 
        $u(1;p)\ge  u(1;p_+)$.
		However, by definition of $q_{\sRS}$, we have 
		$\dot p_+(q_{\sRS}-)=b\dot p_*(q_{\sRS}-)=0$,
		while by concavity of $p$, $\pd_+(q_{\sRS}+)=\pd(q_{\sRS}+)\ge (1-b)/(1-q_{\sRS})>0$.
		In other words $\pd_+(q_{\sRS}-)<\pd_+(q_{\sRS}+)$, and therefore 
		by Lemma \ref{lemma:U-strict-inequality-gen},
		$u(1;p_+)>U_*(\Ccont)$. Corollary \ref{coro:Achieved_on_CCV} and the
		comparison above now give
		\begin{align}
		u(1;p)\ge u(1;p_+)>U_*(\Ccont)=U_*(\Cccv),
		\end{align}
		contradicting the choice of $p$.

		\paragraph{Case $t_*(p)>q_{\sRS}$ and $p(q_{\sRS})>0$.}
		In this case, $\theta(p)\ge q_{\sRS}$. We then define, for $\eps>0$
		\begin{align}
			p_{\eps}(t) = 
			\begin{cases}
				\eps\cdot p_*(t) & \mbox{ if } t\in [0,q_{\sRS}]\, ,\\
				p(t)& \mbox{ if } t\in (q_{\sRS},1]\, ,
			\end{cases}
		\end{align}
		By the ODE \eqref{eq:ODE-integral}, we have, for $t\le \theta$
		\begin{align}
		u(\theta;p_{\eps})-u(\theta;p) &=-\int_{0}^{\theta} 
		L(p_{\eps}(s),\dot p_{\eps}(s),u(s;p_{\eps}),s)\, \de s\nonumber\\
		& \le  -\int_{0}^{q_{\sRS}} 
		L(p_{\eps}(s),\dot p_{\eps}(s),u(s;p_{\eps}),s)\, \de s\nonumber\\
		& = \eps \Big\{u_{\sRS}(q_{\sRS},\alpha,\xi)-\frac{1}{2}\xi(q_{\sRS})\big\}<0\, .
		\end{align}
		Again, by existence and uniquencess of the solution
		for $u>0$, we obtain $u(1;p_{\eps})<u(1;p)$, yielding a contradiction. 
	\end{proof}

\bibliographystyle{amsalpha}

\providecommand{\bysame}{\leavevmode\hbox to3em{\hrulefill}\thinspace}
\providecommand{\MR}{\relax\ifhmode\unskip\space\fi MR }
\providecommand{\MRhref}[2]{%
  \href{http://www.ams.org/mathscinet-getitem?mr=#1}{#2}
}
\providecommand{\href}[2]{#2}

\end{document}